\documentclass[a4paper,10pt,oneside]{article} 
\usepackage{geometry}
\usepackage{standalone}
\usepackage[T1] {fontenc}
\usepackage [utf8] {inputenc}
\usepackage[english] {babel}
\usepackage{xcolor}
\usepackage{graphicx}
\usepackage{authblk}
\usepackage{bbm}
\usepackage{amsmath,amssymb,amsthm, xcolor, listings}
\usepackage{mathrsfs}
\usepackage{tikz}
\usetikzlibrary{patterns}
\usetikzlibrary{decorations.pathreplacing}
\usepackage{paralist}
\usepackage{fancyhdr}
\usepackage{mathtools}
\usepackage{titlesec}
\usepackage{enumitem}
\usepackage{float,placeins}
\usepackage{url}
\usepackage{hyperref}
\usepackage{bm}
\usepackage{mathdots}
\usepackage{microtype} 
\usepackage{eurosym}

\newtheorem{theorem}{Theorem}[section]
\newtheorem{proposition}{Proposition}[section]
\newtheorem{lemma}{Lemma}[section]
\newtheorem{cor}{Corollary}[section]

\theoremstyle{definition}
\newtheorem{definition}{Definition}[section]
\newtheorem{remark}{Remark}[section]

\numberwithin{equation}{section}
\newtheorem{assumption}{Assumption}[section]

\newcommand{\mycircleblack}{\tikz\draw[black, fill=black!88!white] (0,0) circle (3pt);}
\newcommand{\mycircleintgray}{\tikz\draw[black, fill=black!50!white] (0,0) circle (3pt);}
\newcommand{\mycirclegray}{\tikz\draw[black, fill=black!28!white] (0,0) circle (3pt);}
\newcommand{\mycirclelightgray}{\tikz\draw[black, fill=black!12!white] (0,0) circle (3pt);}
\newcommand{\mycircleerreprimo}{\tikz\draw[black, fill=black!70!white] (0,0) circle (3pt);}

\begin{document}
\title{Metastability for the degenerate Potts Model with\\ positive external magnetic field under Glauber dynamics.}

\author[
         {}\hspace{0.5pt}\protect\hyperlink{hyp:email1}{1},\protect\hyperlink{hyp:affil1}{a}
        ]
        {\protect\hypertarget{hyp:author1}{Gianmarco Bet}}

\author[
         {}\hspace{0.5pt}\protect\hyperlink{hyp:email2}{2},\protect\hyperlink{hyp:affil1}{a},\protect\hyperlink{hyp:corresponding}{$\dagger$}
        ]
        {\protect\hypertarget{hyp:author2}{Anna Gallo}}

\author[
         {}\hspace{0.5pt}\protect\hyperlink{hyp:email3}{3},\protect\hyperlink{hyp:affil1}{a},\protect\hyperlink{hyp:affil2}{b}
        ]
        {\protect\hypertarget{hyp:author3}{Francesca R.~Nardi}}

\affil[ ]{
          \small\parbox{365pt}{
             \parbox{5pt}{\textsuperscript{\protect\hypertarget{hyp:affil2}{a}}}Università degli Studi di Firenze,
            \enspace
             \parbox{5pt}{\textsuperscript{\protect\hypertarget{hyp:affil1}{b}}}Eindhoven University of Technology
            }
          }

\affil[ ]{
          \small\parbox{365pt}{
             \parbox{5pt}{\textsuperscript{\protect\hypertarget{hyp:email1}{1}}}\texttt{\footnotesize\href{mailto:gianmarco.bet@unifi.it}{gianmarco.bet@unifi.it}},
             \parbox{5pt}{\textsuperscript{\protect\hypertarget{hyp:email2}{2}}}\texttt{\footnotesize\href{mailto:anna.gallo1@stud.unifi.it}{anna.gallo1@stud.unifi.it}},
             \parbox{5pt}{\textsuperscript{\protect\hypertarget{hyp:email3}{3}}}\texttt{\footnotesize\href{mailto:francescaromana.nardi@unifi.it}{francescaromana.nardi@unifi.it}}
            }
          }

\affil[ ]{
          \small\parbox{365pt}{
             \parbox{5pt}{\textsuperscript{\protect\hypertarget{hyp:corresponding}{$\dagger$}}}Corresponding author
            }
          }

\date{\today}

\maketitle

\begin{abstract}
We consider the ferromagnetic q-state Potts model on a finite grid graph with non-zero external field and periodic boundary conditions. The system evolves according to Glauber-type dynamics described by the Metropolis algorithm, and we focus on the low temperature asymptotic regime. We analyze the case of positive external magnetic field. In this energy landscape there are $1$ stable configuration and $q-1$ metastable states. We study the asymptotic behavior of the first hitting time from any metastable state to the stable configuration as $\beta\to\infty$ in probability, in expectation, and in distribution. We also identify the exponent of the mixing time and find an upper and a lower bound for the spectral gap. We also geometrically identify the union of all minimal gates and the tube of typical trajectories for the transition from any metastable state to the unique stable configuration.

\\
 \emph{Acknowledgment:} The research of Francesca R.~Nardi was partially supported by the NWO Gravitation Grant 024.002.003--NETWORKS and by the PRIN Grant 20155PAWZB ``Large Scale Random Structures''.

\end{abstract}

\maketitle
\tableofcontents
\section{Introduction}
Metastability is a phenomenon that is observed when a physical system is close to a first order phase transition. More precisely, this phenomen takes place when the physical system, for some specific values of the parameters, is imprisoned for a long time in a state which is different from the equilibrium state. The former is known as the \textit{metastable state}, the latter is the \textit{stable state}. After a long (random) time, the system may exhibit the so-called \textit{metastable behavior} and this happens when the system performs a sudden transition from the metastable state to the stable state. On the other hand, when the system lies on the phase coexistence line, it is of interest to understand precisely the transition between two (or more) stable states. This is the so-called \textit{tunneling behavior}.

The phenomenon of metastability occurs in several physical situations, such as supercooled liquids, supersaturated gases, ferromagnets in the hysteresis loop and wireless networks. For this reason, many models for the metastable behavior have been developed throughout the years. In these models 
a suitable stochastic dynamics is chosen and typically three main issues are investigated. The first is the study of the \textit{first hitting time} at which the process starting from a metastable state visits some stable states. The second issue is the study of the so-called set of \textit{critical configurations}, which are those configurations that the process visits during the transition from the metastable state to some stable states. The third issue is the study of the \textit{tube of typical paths},  i.e., the set of the typical trajectories followed by the process during the transition from the metastable state to some stable states. On the other hand, when a system exhibits tunneling behavior the same three issues above are investigated for the transition between any two stable states.\\

In this paper we study the metastable behavior of the $q$-state Potts model with non-zero external magnetic field on a finite two-dimensional discrete torus $\Lambda$. We will refer to $\Lambda$ as a \textit{grid-graph}. The $q$-state Potts model is an extension of the classical Ising model from $q=2$ to an arbitrary number $q$ of spins with $q>2$. The state space $\mathcal X$ is given by all possible configurations $\sigma$ such that at each site $i$ of $\Lambda$ lies a spin with value $\sigma(i)\in\{1,\dots,q\}$.
To each configuration $\sigma\in\mathcal X$ we associate an energy $H(\sigma)$ that depends on the local ferromagnetic interaction, $J=1$, between nearest-neighbor spins, and on the external magnetic field $h$ related only to a specific spin value. Without loss of generality, we choose this spin equal to the spin $1$. We study the $q$-state ferromagnetic Potts model with Hamiltonian $H(\sigma)$ in the limit of large inverse temperature $\beta\to\infty$. The stochastic evolution is described by a \textit{Glauber-type dynamics}, that is a Markov chain on the finite state space $\mathcal X$ with transition probabilities that allow single spin-flip updates and that is given by the Metropolis algorithm.  This dynamics is reversible with respect to the stationary distribution that is the \textit{Gibbs measure} $\mu_\beta$, see \eqref{gibbs}.\\

Our analysis focuses on the model to which we will refer as \textit{$q$-Potts model with positive external magnetic field}. In this energy landscape there are $q-1$ degenerate-metastable states and only one stable state. In the metastable configurations all spins are equal to some $m$, for $m\in\{2,\dots,q\}$, while the stable state is the configuration in which all spins are equal to $1$. In this case, we focus our attention on the transition from one of the metastable states to the stable configuration.

The goal of this paper is to investigate all the three issues of metastability for the metastable behavior of the $q$-Potts model with positive external magnetic field. More precisely, we investigate the asymptotic behavior of the transition time, and identify the set of critical configurations and the tube of typical trajectories for the transition from a metastable state to the unique stable state. Furthermore, we also identify the union of all the critical configurations for the transition from a metastable configuration to the other metastable states.

Let us now briefly describe the strategy that we adopt. First we show that the metastable set contains the $q-1$ configurations where all spins are equal to some $m\in\{2,\dots,q\}$. We give asymptotic bounds in probability for the first hitting time from any metastable state to the stable configuration and we identify the order of magnitude of the expected hitting time. Moreover, we characterize the behavior of the mixing time in the low-temperature regime and give a lower and un upper bound of the spectral gap. Finally, we find the set of all minimal gates for the transition from a metastable state to the stable configuration. For any $m\in\{2,\dots,q\}$, if the starting configuration is the one with all spins equal to $m$, we prove that the minimal gate contains those configurations in which all spins are $m$ except those, which are $1$, in a quasi-square with a unit protuberance on one of the longest sides. Furthermore, we prove that during the metastable transition the process almost surely does not visit any metastable state different from the initial one, and we exploit this result to identify the union of all minimal gates for the transition from a metastable state to the other metastable configurations. 
Finally, we identify geometrically those configurations that belong to the tube of typical paths for the transition from any metastable state to the stable state.
 \\

The Potts model is one of the most studied statistical physics models, as the vast literature on the subject attests, both on the mathematics side and the physics side. The study of the equilibrium properties of the Potts model and their dependence on $q$, have been investigated on the square lattice $\mathbb Z^d$ in  \cite{baxter1973potts,baxter1982critical}, on the triangular lattice in \cite{baxter1978triangular,enting1982triangular} and on the Bethe lattice in \cite{ananikyan1995phase,de1991metastability,di1987potts}. The mean-field version of the Potts model has been studied in \cite{costeniuc2005complete,ellis1990limit,ellis1992limit,gandolfo2010limit,wang1994solutions}. Furthermore, the tunneling behaviour for the Potts model with zero external magnetic field has been studied in \cite{nardi2019tunneling,bet2021critical,kim2021metastability}. In this energy landscape there are $q$ stable states and there is not any relevant metastable state.
In \cite{nardi2019tunneling}, the authors derive the asymptotic behavior of the \textit{first hitting time} for the transition between stable configurations, and give results in probability, in expectation and in distribution. They also characterize the behavior of the \textit{mixing time} and give a lower and an upper bound for the \textit{spectral gap}. In \cite{bet2021critical}, the authors study the tunneling from a stable state to the other stable configurations and between two stable states. In both cases, they geometrically identify the union of all minimal gates and the tube of typical trajectories. Finally, in \cite{kim2021metastability}, the authors study the model in two and three dimensions. In both cases, they give a description of \textit{gateway configurations} that is suitable to allow them to prove sharp estimate for the tunneling time by computing the so-called \textit{prefactor}. These \textit{gateway configurations} are quite different from the states belonging to the \textit{minimal gates} identified by \cite{bet2021critical}.
The $q$-Potts model with non-zero external magnetic field has been studied in \cite{bet2021metastabilityneg}, where the authors study the energy landscape defined by a Hamiltonian function with \textit{negative} external magnetic field. In this scenario there are a unique metastable configuration and $q-1$ stable states, and the authors answer to all the three issues of the metastability introduced above for the transition from the metastable state to the set of the stable states and also to any fixed stable state. Furthermore, they give sharp estimates on the expected transition time by computing the \textit{prefactor}.

\paragraph{State of art} All grouped citations here and henceforth are in chronological order of publication. In this paper we adopt the framework known as \textit{pathwise approach}, that was initiated in 1984 by Cassandro, Galves, Olivieri, Vares in \cite{cassandro1984metastable} and it was further developed in \cite{olivieri1995markov,olivieri1996markov,olivieri2005large} and independtly in \cite{catoni1997exit}. The pathwise approach is based on a detailed knowledge of the energy landscape and, thanks to ad hoc large deviations estimates it gives a quantitative answer to the three issues of metastability which we described above. This approach was further developed in \cite{manzo2004essential,cirillo2013relaxation,cirillo2015metastability,nardi2016hitting,fernandez2015asymptotically,fernandez2016conditioned} to distinguish the study of the transition time and of the critical configurations from the study of the third issue. This is achieved proving the recurrence property and identifying the communication height between the metastable and the stable state that are the only two model-dependent inputs need for the results concerning the first issue of metastability. In particular, in \cite{manzo2004essential,cirillo2013relaxation,cirillo2015metastability,nardi2016hitting,fernandez2015asymptotically,fernandez2016conditioned} this method has been exploited to find answers valid increasing generality in order to reduce as much as possible the number of model dependent inputs necessary to study the metastable and tunneling behaviour, and to consider situations in which the energy landscapes has multiple stable and/or metastable states. 
For this reason, the pathwise approach has been used to study metastability in statistical mechanics lattice models.
The pathwise approach has been also applied in \cite{arous1996metastability,cirillo1998metastability,cirillo1996metastability,kotecky1994shapes,nardi1996low,neves1991critical,neves1992behavior,olivieri2005large} with the aim of answering to the three issues for Ising-like models with Glauber dynamics. Moreover, it was also applied in \cite{hollander2000metastability,den2003droplet,gaudilliere2005nucleation,apollonio2021metastability,nardi2016hitting,zocca2019tunneling} to study the transition time and the gates for Ising-like and hard-core models with Kawasaki and Glauber dynamics. Furthermore, this method was applied to probabilistic cellular automata (parallel dynamics) in \cite{cirillo2003metastability,cirillo2008competitive,cirillo2008metastability,procacci2016probabilistic,dai2015fast}.

The pathwise approach is not the only method which is applied to study the physical systems that approximate a phenomenon of metastability. For instance, the so-called \emph{potential-theoretical approach} exploits a suitable Dirichlet form and spectral properties of the transition matrix to investigate the sharp asymptotic of the hitting time. An interesting aspect of this method is that it allows to estimate the expected value of the transition time including the so-called \textit{prefactor}, i.e., the
coefficient that multiplies the leading-order exponential factor. To find these results, it is necessary to prove the recurrence property, the communication height between the metastable and the stable state and a detailed knowledge of the critical configurations as well of those configurations connected with them by one step of the dynamics, see \cite{bovier2002metastability,bovier2004metastability,bovier2016metastability,cirillo2017sum}. In particular, the potential theoretical approach was applied to find the prefactor for Ising-like models and the hard-core model in \cite{bashiri2019on,boviermanzo2002metastability,cirillo2017sum,bovier2006sharp,den2012metastability,jovanovski2017metastability,den2018metastability} for Glauber and Kawasaki dynamics and in \cite{nardi2012sharp,bet2020effect} for parallel dynamics.
Recently, other approaches have been formulated in \cite{beltran2010tunneling,beltran2012tunneling,gaudillierelandim2014} and in \cite{bianchi2016metastable} and they are particularly adapted to estimate the pre-factor when dealing with the tunnelling between two or more stable states.\\

\paragraph{Outline}The outline of the paper is as follows. In Section \ref{moddescr} we define the ferromagnetic $q$-state Potts model and the Hamiltonian that we associate to each Potts configuration. In Section \ref{defnot} we give a list of both model-independent and model dependent definitions that occur to state our main results in Section \ref{mainres}. In Section \ref{energylandscape} we analyse the energy landscape and give the explicit proofs of the main results stated in Subsections \ref{mainresenergylandscape} and \ref{mainrestime}. Subsection \ref{subsectionnew} is devoted to the study on the transition from a metastable state to the other metastable configurations. In Subsections \ref{proofgates} and \ref{prooftube} we give the explicit proofs of the main results on the critical configurations and on the tube of typical paths, respectively. 
\section{Model description}\label{moddescr}
In the $q$-state Potts model each spin lies on a vertex of a finite two-dimensional rectangular lattice $\Lambda=(V,E)$, where $V=\{0,\dots,K-1\}\times\{0,\dots,L-1\}$ is the vertex set and $E$ is the edge set, namely the set of the pairs of vertices whose spins interact with each other. We consider periodic boundary conditions. More precisely, we identify
each pair of vertices lying on opposite sides of the rectangular lattice, so that we obtain a two-dimensional torus. Two vertices $v,w\in V$ are said to be nearest-neighbors when they share an edge of $\Lambda$. 
We denote by $S$ the set of spin values, i.e., $S:= \{1,\dots,q\}$ and assume $q>2$. To each vertex $v\in V$ is associated a spin value $\sigma(v)\in S$, and $\mathcal X := S^V$ denotes the set of spin
configurations. We denote by $\textbf{1},\dots,\textbf{q} \in \mathcal X$ those configurations in which all the vertices have spin value $1,\dots,q$, respectively.  

\noindent To each configuration $\sigma\in\mathcal X$ we associate the energy $H(\sigma)$ given by
\begin{align}\label{hamiltoniangeneral}
H(\sigma):=-J \sum_{(v,w)\in E} \mathbbm{1}_{\{\sigma(v)=\sigma(w)\}}-h\sum_{u\in V} \mathbbm{1}_{\{\sigma(u)=1\}}, \ \ \sigma\in\mathcal X,
\end{align}
where $J$ is the \textit{coupling} or \textit{interation constant} and $h$ is the \textit{external magnetic field}. The function $H:\mathcal X\to\mathbb R$ is called \textit{Hamiltonian} or \textit{energy function}. The Potts model is said to be \textit{ferromagnetic} when $J>0$, and \textit{antiferromagnetic} otherwise. In this paper we set $J=1$ without loss of generality and, we focus on the ferromagnetic $q$-state Potts model with non-zero external magnetic field. More precisely, we study the model with positive external magnetic field, i.e., we rewrite \eqref{hamiltoniangeneral} by considering the magnetic field $h_{\text{pos}}:=h$, 
\begin{align}\label{hamiltonianpos}
H_{\text{pos}}(\sigma) :&= - \sum_{(v,w)\in E} \mathbbm{1}_{\{\sigma(v)=\sigma(w)\}}-h_{\text{pos}}\sum_{u\in V} \mathbbm{1}_{\{\sigma(u)=1\}}\notag\\
&= - \sum_{(v,w)\in E} \mathbbm{1}_{\{\sigma(v)=\sigma(w)\}}-h\sum_{u\in V} \mathbbm{1}_{\{\sigma(u)=1\}}.
\end{align}

The \textit{Gibbs measure} for the $q$-state Potts model on $\Lambda$ is a probability distribution on the state space $\mathcal X$ given by
\begin{align}\label{gibbs}
\mu_\beta(\sigma):=\frac{e^{-\beta H(\sigma)}}{Z},
\end{align}
where $\beta>0$ is the inverse temperature and where $Z:=\sum_{\sigma'\in\mathcal X}e^{-\beta H_\text{pos}(\sigma')}$.

The spin system evolves according to a Glauber-type dynamics. This dynamics is described by a single-spin update Markov chain $\{X_t^\beta\}_{t\in\mathbb{N}}$ on the state space $\mathcal X$ with the following transition probabilities: for $\sigma, \sigma' \in \mathcal X$,
\begin{align}\label{metropolisTP}
P_\beta(\sigma,\sigma'):=
\begin{cases}
Q(\sigma,\sigma')e^{-\beta [H_\text{pos}(\sigma')-H_\text{pos}(\sigma)]^+}, &\text{if}\ \sigma \neq \sigma',\\
1-\sum_{\eta \neq \sigma} P_\beta (\sigma, \eta), &\text{if}\ \sigma=\sigma',
\end{cases}
\end{align}
where $[n]^+:=\max\{0,n\}$ is the positive part of $n$ and 
\begin{align}\label{Qmatrix}
Q(\sigma,\sigma'):=
\begin{cases}
\frac{1}{q|V|}, &\text{if}\ |\{v\in V: \sigma(v) \neq \sigma'(v)\}|=1,\\
0, &\text{if}\ |\{v\in V: \sigma(v) \neq \sigma'(v)\}|>1,
\end{cases}
\end{align}
for any $\sigma, \sigma' \in \mathcal X$. $Q$ is the so-called \textit{connectivity matrix} and it is symmetric and irreducible, i.e., for all $\sigma, \sigma' \in \mathcal X$, there exists a finite sequence of configurations $\omega_1,\dots,\omega_n \in \mathcal X$ such that $\omega_1=\sigma$, $\omega_n=\sigma'$ and $Q(\omega_i,\omega_{i+1})>0$ for $i=1,\dots,n-1$.  Hence, the resulting stochastic dynamics defined by \eqref{metropolisTP} is reversible with respect to the Gibbs measure \eqref{gibbs}. The triplet $(\mathcal X,H,Q)$ is the so-called \textit{energy landscape}. \\

The dynamics defined above belongs to the class of Metropolis dynamics. Given a configuration $\sigma$ in $\mathcal X$, at each step
\begin{itemize}
\item[1.] a vertex $v \in V$ and a spin value $s \in S$ are selected independently and uniformly at random;
\item[2.] the spin at $v$ is updated to spin $s$ with probability
\begin{align}\begin{cases}
1, &\text{if}\ H_\text{pos}(\sigma^{v,s)})-H_\text{pos}(\sigma) \le 0,\\
e^{-\beta[H_\text{pos}(\sigma^{v,s})-H_\text{pos}(\sigma)]},&\text{if}\ H_\text{pos}(\sigma^{v,s})-H_\text{pos}(\sigma) > 0,
\end{cases}\end{align}
\end{itemize}
where $\sigma^{v,s}$ is the configuration obtained from $\sigma$ by updating the spin in the vertex $v$ to $s$, i.e., 
\begin{align}\label{confspinflip}
\sigma^{v,s}(w):=
\begin{cases}
\sigma(w)\ &\text{if}\ w\neq v,\\
s\ &\text{if}\ w=v.
\end{cases}
\end{align}
Hence, at each step the update of vertex $v$ depends on the neighboring spins of $v$ and on the energy difference
\begin{align}\label{energydifference3}
H_{\text{pos}}(\sigma^{v,s})-H_{\text{pos}}(\sigma)= 
\begin{cases}
\sum_{w \sim v} (\mathbbm{1}_{\{\sigma(v)=\sigma(w)\}}-\mathbbm{1}_{\{\sigma(w)=s\}})+h, &\text{if}\ \sigma(v)=1,\ s\neq1,\\
\sum_{w \sim v} (\mathbbm{1}_{\{\sigma(v)=\sigma(w)\}}-\mathbbm{1}_{\{\sigma(w)=s\}}), &\text{if}\ \sigma(v)\neq1,\ s\neq1,\\
\sum_{w \sim v} (\mathbbm{1}_{\{\sigma(v)=\sigma(w)\}}-\mathbbm{1}_{\{\sigma(w)=s\}})-h, &\text{if}\ \sigma(v)\neq1,\ s=1.
\end{cases}
\end{align}

\section{Definitions and notations}\label{defnot}
In order to state our main results, we need to give some definitions and notations which are used throughout the next sections.
\subsection{Model-independent definitions and notations}\label{modinddef}
We now give a list of model-independent definitions and notations that will be useful in formulating our main results.
\begin{itemize}
\item[-] We call \textit{path} a finite sequence $\omega$ of configurations $\omega_0,\dots,\omega_n \in \mathcal X$, $n \in \mathbb{N}$, such that $Q(\omega_i,\omega_{i+1})>0$ for $i=0,\dots,n-1$. Given $\sigma, \sigma' \in \mathcal X$, if $\omega_1=\sigma$ and $\omega_n=\sigma'$, we denote a path from $\sigma$ to $\sigma'$ as $\omega: \sigma \to \sigma'$. 
\item[-] Let $\Omega_{\sigma,\sigma'}$ be the set of all paths between $\sigma$ and $\sigma'$.
\item[-] Given a path $\omega=(\omega_0,\dots,\omega_n)$, we define the \textit{height} of $\omega$ as
\begin{align}\label{height}
\Phi_\omega:=\max_{i=0,\dots,n} H(\omega_i).
\end{align}
\item[-] We say that a path $\omega\in\Omega_{\sigma,\sigma'}$ is the \textit{concatenation} of the $L$ paths \[\omega^{(i)}=(\omega^{(i)}_0,\dots,\omega^{(i)}_{n_i}),\ \text{for some}\ n_i\in\mathbb{N},\ i=1,\dots,L\] if \[\omega=(\omega^{(1)}_0=\sigma,\dots,\omega^{(1)}_{n_1},\omega^{(2)}_0,\dots,\omega^{(2)}_{n_2},\dots,\omega^{(L)}_0,\dots,\omega^{(L)}_{n_L}=\sigma').\] 
\item[-] A path $\omega=(\omega_0,\dots,\omega_n)$ is said to be \textit{downhill} (\textit{strictly downhill}) if $H(\omega_{i+1})\le H(\omega_i)$ ($H(\omega_{i+1})<H(\omega_i)$) for $i=0,\dots,n-1$.
\item[-] For any pair $\sigma, \sigma' \in \mathcal X$, the \textit{communication height} or \textit{communication energy} $\Phi(\sigma,\sigma')$ between $\sigma$ and $\sigma'$ is the minimal energy across all paths $\omega:\sigma \to \sigma'$, i.e., 
\begin{align}\label{comheight}
\Phi(\sigma,\sigma'):=\min_{\omega:\sigma \to \sigma'} \Phi_\omega = \min_{\omega:\sigma \to \sigma'} \max_{\eta \in \omega} H(\eta). 
\end{align}
More generally, the communication energy between any pair of non-empty disjoint subsets $\mathcal A,\mathcal B \subset \mathcal X$ is $\Phi(\mathcal A,\mathcal B):=\min_{\sigma \in \mathcal A,\ \sigma' \in \mathcal{B}} \Phi(\sigma,\sigma').$
\item[-] We define \textit{optimal paths} those paths that realize the min-max in \eqref{comheight} between $\sigma$ and $\sigma'$. Formally, we define the set of \textit{optimal paths} between $\sigma, \sigma' \in\mathcal X$ as
\begin{align}\label{optpaths}
\Omega_{\sigma,\sigma'}^{opt}:=\{\omega\in\Omega_{\sigma,\sigma'}:\ \max_{\eta\in\omega} H(\eta)=\Phi(\sigma,\sigma')\}.
\end{align}

\item[-] For any $\sigma \in \mathcal X$, let
\begin{align}
\mathcal{I}_\sigma:=\{\eta \in \mathcal X:\ H(\eta)<H(\sigma)\}
\end{align}
be the set of states with energy strictly smaller than $H(\sigma)$.
We define \textit{stability level} of $\sigma$ the energy barrier
\begin{align}\label{stabilitylevel}
V_\sigma:=\Phi(\sigma,\mathcal{I}_\sigma)-H(\sigma).
\end{align}
If $\mathcal{I}_\sigma = \varnothing$, we set $V_\sigma:=\infty$.
\item[-] The bottom $\mathscr{F}(\mathcal A)$ of a non-empty set $\mathcal A\subset\mathcal X$ is the set of the \textit{global minima} of $H$ in $\mathcal A$, i.e.,
\begin{align}\label{bottom}
\mathscr{F}(\mathcal A):=\{\eta \in \mathcal A:H(\eta)=\min_{\sigma \in \mathcal A}H(\sigma)\}. 
\end{align}
In particular, $\mathcal X^s:=\mathscr{F}(\mathcal X)$ is the set of the \textit{stable states}.
\item[-] For any $\sigma\in\mathcal X$ and any $\mathcal A\subset\mathcal X$, $\mathcal A\neq\varnothing$, we set
\begin{align}\label{gammatwo}
\Gamma(\sigma,\mathcal A):=\Phi(\sigma,\mathcal A)-H(\sigma).
\end{align} 
\item[-] We define the set of \textit{metastable states} as 
\begin{align}\label{metastableset}
\mathcal X^m:=\{\eta \in \mathcal X: V_\eta=\max_{\sigma\in\mathcal X\backslash\mathcal X^s} V_\sigma\}.
\end{align}
We denote by $\Gamma^m$ the stability level of a metastable state.
\item[-] We define \textit{metastable set at level $V$} the set of all the configurations with stability level larger tha $V$, i.e.,
\begin{align}\label{metasetV}
\mathcal X_V:=\{\sigma\in\mathcal X:V_\sigma>V\}.
\end{align}
\item[-] The set of \textit{minimal saddles} between $\sigma, \sigma' \in \mathcal X$ is defined as
\begin{align}\label{saddles}
\mathcal S(\sigma,\sigma'):=\{\xi\in\mathcal X:\exists\omega\in\Omega_{\sigma,\sigma'}^{opt},\ \xi\in\omega:\ \max_{\eta\in\omega} H(\eta)=H(\xi)\}.
\end{align}
\item[-] We say that $\eta\in\mathcal S(\sigma,\sigma')$ is an \textit{essential saddle} if there exists $\omega\in\Omega_{\sigma,\sigma'}^{opt}$ such that either
\begin{itemize}
\item  $\{\text{arg max}_\omega H\}=\{\eta\}$ or
\item $\{\text{arg max}_\omega H\}\supset\{\eta\}$ and $\{\text{arg max}_{\omega'} H\}\not\subseteq\{\text{arg max}_\omega H\}\backslash \{\eta\}$ for all $\omega'\in\Omega_{\sigma,\sigma'}^{opt}$.
\end{itemize}
\item[-] A saddle $\eta\in\mathcal S(\sigma,\sigma')$ that is not essential is said to be \textit{unessential}. 
\item[-] Given $\sigma, \sigma' \in \mathcal X$, we say that $\mathcal W(\sigma,\sigma')$ is a \textit{gate} for the transition from $\sigma$ to $\sigma'$ if  $\mathcal W(\sigma,\sigma')\subseteq\mathcal S(\sigma,\sigma')$ and $\omega\cap\mathcal W(\sigma,\sigma')\neq\varnothing$ for all $\omega\in\Omega_{\sigma,\sigma'}^{opt}$.
\item[-] We say that  $\mathcal W(\sigma,\sigma')$ is a \textit{minimal gate} for the transition from $\sigma$ to $\sigma'$ if it is a minimal (by inclusion) subset of $\mathcal S(\sigma,\sigma')$ that is visited by all optimal paths. More in detail, it is a gate and for any $\mathcal W'\subset\mathcal W(\sigma,\sigma')$ there exists $\omega'\in\Omega_{\sigma,\sigma'}^{opt}$ such that $\omega'\cap\mathcal W'=\varnothing$. We denote by $\mathcal{G}=\mathcal{G}(\sigma,\sigma')$ the union of all minimal gates for the transition $\sigma\to\sigma'$.
\item[-] Given a non-empty subset $\mathcal A \subset \mathcal X$ and a configuration $\sigma \in \mathcal X$, we define
\begin{align}\label{firsthittingtime}
\tau_\mathcal A^\sigma := \text{inf}\{t>0: \ X_t^\beta \in \mathcal A\}
\end{align}
as the \textit{first hitting time} of the subset $\mathcal A$ for the Markov chain $\{X_t^\beta\}_{t \in \mathbb{N}}$ starting from $\sigma$ at time $t=0$. 
\item[-] Let $\{X_t^\beta\}_{t\in\mathbb N}$ be the Markov chain with transition probabilities \eqref{metropolisTP} and stationary distribution \eqref{gibbs}. For every $\epsilon\in(0,1)$, we define the \textit{mixing time} $t^{\text{mix}}_\beta(\epsilon)$ by
\begin{align}\label{mixingtimedef}
t^{\text{mix}}_\beta(\epsilon):=\min\{n\ge 0|\ \max_{\sigma\in\mathcal X}||P_\beta^n(\sigma,\cdot)-\mu_\beta(\cdot)||_{\text{TV}}\le\epsilon\},
\end{align}
where the total variance distance is defined by $||\nu-\nu'||_{\text{TV}}:=\frac 1 2 \sum_{\sigma\in\mathcal X}|\nu(\sigma)-\nu'(\sigma)|$ for every two probability distribution $\nu,\nu'$ on $\mathcal X$. Furthermore, we define \textit{spectral gap} as
\begin{align}\label{spectralgapdef}
\rho_\beta:=1-\lambda_\beta^{(2)},
\end{align}
where $1=\lambda_\beta^{(1)}>\lambda_\beta^{(2)}\ge\dots\ge\lambda_\beta^{(|\mathcal X|)}\ge-1$ are the eigenvalues of the matrix $P_\beta(\sigma,\eta))_{\sigma,\eta\in\mathcal X}$.
\item[-] Given a non-empty subset $\mathcal A\subseteq\mathcal X$, it is said to be \textit{connected} if for any $\sigma,\eta\in\mathcal A$ there exists a path $\omega:\sigma\to\eta$ totally contained in $\mathcal A$. Moreover, we define $\partial A$ as the \textit{external boundary} of $\mathcal A$, i.e., the set
\begin{align}\label{boundaryset}
\partial\mathcal A:=\{\eta\not\in\mathcal A:\ P(\sigma,\eta)>0\ \text{for some}\ \sigma\in\mathcal A\}.
\end{align}
%
%
%
%
\item[-] A non-empty subset $\mathcal C\subset\mathcal X$ is called \textit{cycle} if it is either a singleton or a connected set such that
\begin{align}\label{cycle}
\max_{\sigma\in\mathcal C} H(\sigma)<H(\mathscr{F}(\partial\mathcal C)).
\end{align}
When $\mathcal C$ is a singleton, it is said to be a \textit{trivial cycle}. Let $\mathscr C(\mathcal X)$ be the set of cycles of $\mathcal X$. 
\item[-] The \textit{depth} of a cycle $\mathcal C$ is given by
\begin{align}\label{defdepcycle}
\Gamma(\mathcal C):=H(\mathscr F(\partial\mathcal C))-H(\mathscr F(\mathcal C)).
\end{align}
If $\mathcal C$ is a trivial cycle we set $\Gamma(\mathcal C)=0$.
\item[-] Given a non-empty set $\mathcal A\subset\mathcal X$, we denote by $\mathcal M(\mathcal A)$ the \textit{collection of maximal cycles} $\mathcal A$, i.e., 
\begin{align}
\mathcal M(\mathcal A):=\{\mathcal C\in\mathscr C(\mathcal X)|\ \mathcal C\ \text{maximal by inclusion under constraint}\ \mathcal C\subseteq\mathcal A\}.\notag
\end{align}
\item[-] For any $\mathcal A\subset\mathcal X$, we define the \textit{maximum depth of $\mathcal A$} as the maximum depth of a cycle contained in $\mathcal A$, i.e., 
\begin{align}\label{defgammatildenzb}
\widetilde{\Gamma}(\mathcal A):=\max_{\mathcal C\in\mathcal M(\mathcal A)} \Gamma(\mathcal C).
\end{align}
In \cite[Lemma 3.6]{nardi2016hitting} the authors give an alternative characterization of \eqref{defgammatildenzb} as the maximum initial energy barrier that the process started from a configuration $\eta\in\mathcal A$ possibly has to overcome to exit from $\mathcal A$, i.e.,
\begin{align}
\widetilde{\Gamma}(\mathcal A)=\max_{\eta\in\mathcal A} \Gamma(\eta,\mathcal X\backslash\mathcal A).
\end{align}
\item[-] For any $\sigma\in\mathcal X$, if $\mathcal{A}$ is a non-empty target set, we define the \textit{initial cycle} for the transition from $\sigma$ to $\mathcal A$ as
\begin{align}\label{initialcycle}
\mathcal C_{\mathcal{A}}^\sigma(\Gamma):=\{\sigma\}\cup\{\eta\in\mathcal X:\ \Phi(\sigma,\eta)-H(\sigma)<\Gamma=\Phi(\sigma,\mathcal A)-H(\sigma)\}.
\end{align}
If $\sigma\in\mathcal A$, then $\mathcal C_{\mathcal{A}}^\sigma(\Gamma)=\{\sigma\}$ and it is a trivial cycle. Otherwise, $\mathcal C_{\mathcal{A}}^\sigma(\Gamma)$ is either a trivial cycle (when $\Phi(\sigma, \mathcal{A})=H(\sigma))$ or a non-trivial cycle containing $\sigma$ (when $\Phi(\sigma, \mathcal{A})>H(\sigma)$). In any case, if $\sigma\notin\mathcal{A}$, then $C_{\mathcal{A}}^\sigma(\Gamma)\cap{\mathcal{A}}=\varnothing$. Note that \eqref{initialcycle} coincides with \cite[Equation (2.25)]{nardi2016hitting}.
\end{itemize}
\begin{itemize}
\item[-] A $\beta\to f(\beta)$ is said to be \textit{super-exponentially small} (SES) if 
\begin{align}
\lim_{\beta\to\infty} \frac 1 \beta \log f(\beta)=-\infty.
\end{align}
%
\end{itemize}
\subsection{Model-dependent definitions and notations}\label{moddepdef}
In this section we give some further model-dependent notations, which hold for any fixed $q$-Potts configuration $\sigma\in\mathcal X$.
\begin{itemize}
\item[-] For any $v,w\in V$, we write $w\sim v$ if there exists an edge $e\in E$ that links the vertices $v$ and $w$.
\item[-] We denote the edge that links the vertices $v$ and $w$ as $(v,w)\in\ E$. Each $v\in V$ is naturally identified by its coordinates $(i,j)$, where $i$ and $j$ denote respectively the number of the row and of the column where $v$ lies. Moreover, the collection of vertices with first coordinate equal to $i=0,\dots,K-1$ is denoted as $r_i$, which is the $i$-th row of $\Lambda$. The collection of those vertices with second coordinate equal to $j=0,\dots,L-1$ is denoted as $c_j$, which is the $j$-th column of $\Lambda$.
\item[-] We define the set $C^s(\sigma)\subseteq\mathbb R^2$ as the union of unit closed squares centered at the vertices $v\in V$ such that $\sigma(v)=s$. We define $s$-\textit{clusters} the maximal connected components $C^s_1,\dots,C^s_n,\ n\in\mathbb N$, of $C^s(\sigma)$. 
\item[-] For any $s\in S$, we say that a configuration $\sigma\in\mathcal X$ has an $s$-rectangle if it has a rectangular cluster in which all the vertices have spin $s$. 
\item[-] Let $R_1$ an $r$-rectangle and $R_2$ an $s$-rectangle. They are said to be \textit{interacting} if either they intersect (when $r=s$) or are disjoint but there exists a site $v\notin R_1\cup R_2$ such that $\sigma(v)\neq r,s$ and $v$ has two nearest-neighbor $w, u$ lying inside $R_1, R_2$ respectively. For instance, in Figure \ref{figurerecurrencepositive}(c) the gray rectangles are not interacting. Furthermore, we say that $R_1$ and $R_2$ are \textit{adjacent} when they are at lattice distance one from each other, see for instance Figure \ref{figurerecurrencepositive}(d). Since in our scenario there are $q$ types of spin, note that here we give a geometric characterization of interacting rectangles. This geometric definition coincides with the usual one, in which flipping the spin on vertex $v$ decreases the energy, only if the two rectangles have the same spin value. 
\item[-] We set $R(C^s(\sigma))$ as the smallest rectangle containing $C^s(\sigma)$.
\item[-] Let $R_{\ell_1\times\ell_2}$ be a rectangle in $\mathbb R^2$ with  sides of length $\ell_1$ and $\ell_2$.
\item[-] Let $s\in S$. If $\sigma$ has a cluster of spins $s$ which is a rectangle that wraps around $\Lambda$, we say that $\sigma$ has an \textit{$s$-strip}. For any $r,s\in S$, we say that an $s$-strip is \textit{adjacent} to an $r$-strip if they are at lattice distance one from each other. For instance, in Figure \ref{figurerecurrencepositive}(a)-(b) there are depicted vertical and horizontal adjacent strips, respectively.
%
%
%
\item[-] $\bar R_{a,b}(r, s)$ denotes the set of those configurations in which all the vertices have spins equal to $r$, except those, which have spins $s$, in a rectangle $a\times b$, see Figure \ref{esempioBRbar}(a). Note that when either $a=L$ or $b=K$,  $\bar R_{a,b}(r, s)$ contains those configurations which have an $r$-strip and an $s$-strip.
\item[-] $\bar B_{a,b}^l(r, s)$ denotes the set of those configurations in which all the vertices have spins $r$, except those, which have spins $s$, in a rectangle $a\times b$ with a bar $1\times l$ adjacent to one of the sides of length $b$, with $1\le l\le b-1$, see Figure \ref{esempioBRbar}(b).
\begin{figure}[h!]
\centering
\begin{tikzpicture}[scale=0.7, transform shape]
\draw [fill=gray,lightgray] (0.6,0.6) rectangle (1.5,2.4);
\draw[step=0.3cm,color=black] (0,0) grid (3.6,2.7);
\draw (1.8,-0.1) node[below] {\LARGE(a)};
\end{tikzpicture}\ \ \ \ \ \ \ \
\begin{tikzpicture}[scale=0.7, transform shape]
\draw [fill=gray,lightgray] (0.9,0.3) rectangle (1.5,2.1);
\draw [fill=gray,lightgray] (1.5,0.6) rectangle (1.8,1.8);
\draw[step=0.3cm,color=black] (0,0) grid (3.6,2.7);
\draw (1.8,-0.1) node[below] {\LARGE(b)};
\end{tikzpicture}
\caption{\label{esempioBRbar} Examples of configurations which belong to $\bar R_{3,6}(r, s)$ (a), $\bar B_{2,6}^4(r,s)$ (b). We color white the vertices whose spin is $r$ and we color gray the vertices whose spin is $s$.}
\end{figure}
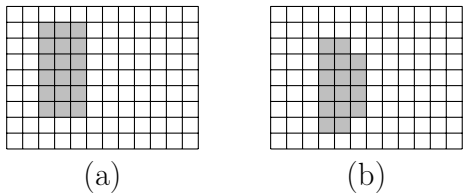
\item[-] We define 
\begin{align}\label{ellestar}
\ell^*:=\left\lceil \frac{2}{h} \right\rceil
\end{align} as the \textit{critical length}.
\end{itemize}

\section{Main results on the $q$-state Potts model with positive external magnetic field}\label{mainres}
This section is devoted to the statement of our main results on the $q$-state Potts model with positive external magnetic field.  Note that we give the proof of the main results by considering the condition $L\ge K\ge3\ell^*$, where $\ell^*$ is defined in \eqref{ellestar}. It is possible to extend the results to the case $K>L$ by interchanging the role of rows and columns in the proof.\\

In this scenario related to the Hamiltonian $H_\text{pos}$, we add either a subscript or a superscript ``pos'' to the notation of the model-independent quantities (defined in general in Subsection \ref{modinddef}) in order to remind the reader that these quantities are computed in the case of positive external magnetic field. For example, we denote the set of the metastable configurations of the Hamiltonian $H_\text{pos}$ by $\mathcal X^m_\text{pos}:=\{\eta\in\mathcal X: V_\eta^\text{pos}=\max_{\sigma\in\mathcal X\backslash\mathcal X^s_\text{pos}} V_\sigma^\text{pos}\}$, where $V_\xi^\text{pos}$ is the stability level related to $H_\text{pos}$ of any configuration $\xi\in\mathcal X$.  

\noindent In order to give the main results on the $q$-Potts model with energy function defined in \eqref{hamiltonianpos}, we have the following assumption.
\begin{assumption}\label{remarkconditionpos}
We assume that the following conditions are verified:
\begin{itemize}
\item[(i)]  the magnetic field $h_{\text{pos}}:=h$ is such that $0<h<1$;
\item[(ii)] $2/h$ is not integer.
\end{itemize}
\end{assumption}

\subsection{Energy landscape}\label{mainresenergylandscape}
Using the definition \eqref{hamiltonianpos} and by simple algebraic calculations, in the following proposition we identificate the set of the global minima of $H_\text{pos}$.
\begin{proposition}[Identification of the stable configuration]\label{stablesetposprop}
Consider the $q$-state Potts model on a $K\times L$ grid $\Lambda$, with periodic boundary conditions and with positive external magnetic field. Then, the set of global minima of the Hamiltonian \eqref{hamiltonianpos} is given by $\mathcal X^s_\emph{pos}:=\{\bold 1\}$.
\end{proposition}
In the next theorem we define the configurations that belong to $\mathcal X^m_{\text{pos}}$ and give an estimate of the stability level $\Gamma_{\text{pos}}^m$. We refer to Figure \ref{figurepositive} for a pictorial representation of the $4$-state Potts model related to the Hamiltonian $H_\text{pos}$.
\begin{theorem}[Identification of the metastable states]\label{teometastablepos}
Consider the $q$-state Potts model on a $K\times L$ grid $\Lambda$, with periodic boundary conditions and with positive external magnetic field. Then,
\begin{align}
\mathcal X^m_{\emph{pos}}=\{\bold 2,\dots,\bold q\}
\end{align}
and, for any $\bold m\in\mathcal X^m_{\emph{pos}}$,
\begin{align}\label{estimatestablevelmetapos}
\Gamma^m_{\emph{pos}}=\Gamma_{\emph{pos}}(\bold m,\mathcal X^s_{\emph{pos}})=4\ell^*-h(\ell^*(\ell^*-1)+1).
\end{align}
\end{theorem}
\textit{Proof.} The theorem follows by \cite[Theorem 2.4]{cirillo2013relaxation} since the first assumption follows by 
  Propositions \ref{refpathpos} and \ref{lowerboundpos} and the second assumption is satisfied thanks to Proposition \ref{ricorrenzapos}. $\qed$

\noindent Using \eqref{estimatestablevelmetapos}, in Subsection \ref{proofenergylandscape} we prove the following corollary. 
\begin{cor}[Maximum depth of a cycle in $\mathcal X\backslash\mathcal X^s_\text{pos}$]\label{corollarygammatildepos}
Consider the $q$-state Potts model on a $K\times L$ grid $\Lambda$, with periodic boundary conditions and with positive external magnetic field. Then, 
\begin{align}\label{aligngammatildepos}
\widetilde{\Gamma}_\emph{pos}(\mathcal X\backslash\mathcal X^s_\emph{pos})=\Gamma^m_\emph{pos}.
\end{align}
\end{cor}

\begin{figure}[h!]
\centering
\begin{tikzpicture}[scale=1,transform shape]
\draw[thick,white] (0.4,0.62) parabola (1.32,2.92);
\draw[dotted] (0.4,0.62) parabola (1.31,2.92);
\draw[thick,white] (0.4,0.62) parabola (-0.51,2.92);
\draw[dotted] (0.4,0.62) parabola (-0.51,2.92);
\draw[white,thick] (-0.51,2.92)--(-0.466,2.7);
\draw (0.4,2.92) ellipse (0.91cm and 0.16cm);
\fill[black!5!white] (0.4,2.92) ellipse (0.91cm and 0.16cm);
\draw (1.31,2.92) -- (1.22,2.48);
\draw (-0.51,2.92)--(-0.474,2.73);
\fill (0.4,0.62) circle (1.3pt);  \draw (0.4,0.62) node [below] {{$\bold 3$}};
\draw[thick,white] (0.72,0.9)--(0.915,1.34);
\draw (0.7,0.85)--(0.912,1.34);
\draw[black] (0,-0.5) parabola (-1.05,2.55);
\draw[black] (0,-0.5) parabola (1.05,2.55);
\draw[white,thick] (1.05,2.55) -- (0.86,1.55);
\draw[black,dotted] (1.05,2.55) -- (0.86,1.55);
\draw[dotted,thick] (0,2.55) ellipse (1.05cm and 0.21cm);
\fill[black!5!white] (0,2.55) ellipse (1.05cm and 0.21cm);
\fill[black] (0,-0.5) circle (1.3pt);  \draw (0,-0.5)  node [below] {{$\bold 1$}};
\draw (1.6,0) parabola (0.69,2.3);\draw (1.6,0) parabola (2.51,2.3);
\draw (1.6,2.3) ellipse (0.91cm and 0.2cm);
\fill[black!5!white]  (1.6,2.3) ellipse (0.91cm and 0.2cm);
\fill (1.6,0) circle (1.3pt);  \draw (1.6,0) node [below] {{$\bold 4$}};
\draw (-1.95,0.2) parabola (-2.86,2.5);\draw (-1.95,0.2) parabola (-1.04,2.5);
\draw (-1.95,2.5) ellipse (0.91cm and 0.2cm);
\fill[black!5!white]  (-1.95,2.5) ellipse (0.91cm and 0.2cm);
\fill (-1.95,0.2) circle (1.3pt); \draw (-1.93,0.2) node [below] {{$\bold 2$}};
\fill (0.85,2.43) circle (1.5pt);
\fill (0.22,2.75) circle (1.5pt);
\fill (-1.05,2.52) circle (1.5pt);
\draw [->,ultra thin] (0.22,2.75)--(-1,3.6) node[above] {{$\bar B_{\ell^*-1,\ell^*}^1(3,1)$}};
\draw [->,ultra thin] (-1.05,2.52) --(-3.1,3); \draw (-2.7,3.2) node[above,left] {{$\bar B_{\ell^*-1,\ell^*}^1(2,1)$}};
\draw [->,ultra thin] (0.85,2.43)--(2.8,3.3) node[above,right] {{$\bar B_{\ell^*-1,\ell^*}^1(4,1)$}};
\draw (4.25,2.3) node {$\Phi_{\text{pos}}(\bold m,\bold 1)$, $\bold m\in\mathcal X^m_{\text{pos}}$};
\end{tikzpicture}
\caption{\label{figurepositive} Energy landscape in the case of $4$-state Potts model with positive external magnetic field around the unique stable state $\bold 1$ cutting the configurations with the energy bigger than $\Phi_\text{pos}(\bold m,\bold 1)$, $\bold m\in\mathcal X^m_\text{pos}=\{\bold 2,\bold 3,\bold 4\}$. This picture is simplified since there are not represented the cycles (valleys) that contain configurations with stability level smaller than or equal to $2$ (see Proposition \ref{ricorrenzapos}).}
\end{figure}
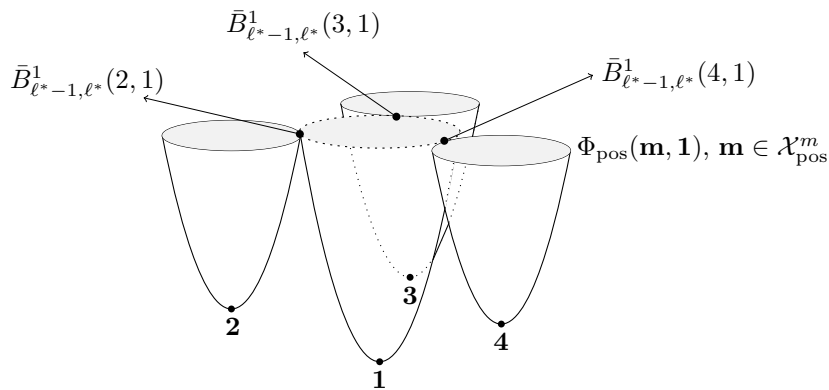\FloatBarrier

In the following proposition, that we prove in Subsection \ref{stablevpos}, we investigate on the stability level of any configuration $\eta\in\mathcal X\backslash\{\bold 2,\dots,\bold q\}$.
\begin{proposition}[Estimate on the stability level]\label{ricorrenzapos}
If the external magnetic field is positive, then for any $\eta\in\mathcal X\backslash \{\bold 1,\dots,\bold q\}$ and $\bold m\in\{\bold 2,\dots,\bold q\}$, $V_\eta^\emph{pos}\le 2<\Gamma_{\emph{pos}}(\bold m,\mathcal X^s_{\emph{pos}})$. 
\end{proposition}
 Exploiting the estimate of the stability level given in Proposition \ref{ricorrenzapos}, we prove the following result on a recurrence property to set of the metastable and stable configurations, i.e., $\{\bold 1,\dots,\bold q\}$.
\begin{theorem}[Recurrence property]\label{teorecproppos}
Consider the $q$-state Potts model on a $K\times L$ grid $\Lambda$, with periodic boundary conditions and with positive external magnetic field. For any  $\sigma\in\mathcal X$ and for any $\epsilon>0$, there exists $k>0$ such that for $\beta$ sufficiently large 
\begin{align}
\mathbb P(\tau^\sigma_{\{\bold 1,\dots,\bold q\}}>e^{\beta(2+\epsilon)})\le e^{-e^{k\beta}}=\emph{SES}.
\end{align}
\end{theorem}
\textit{Proof}. Apply \cite[Theorem 3.1]{manzo2004essential} with $V=2$ and use \eqref{metasetV} and Proposition \ref{ricorrenzapos} to get $\mathcal X_2=\{\bold 1,\dots,\bold q\}=\mathcal X^m_\text{pos}\cup\mathcal X^s_\text{pos}$, where the last equality follows by Proposition \ref{stablesetposprop} and Theorem \ref{teometastablepos}. $\qed$

\subsection{Asymptotic behavior of the first hitting time to the stable state and mixing time}\label{mainrestime}
We are interested in studying the first hitting time $\mathcal X^s_\text{pos}=\{\bold 1\}$ of the stable configuration starting from any $\bold m\in\mathcal X^m_\text{pos}$, i.e., $\tau^\bold m_{\mathcal X^s_\text{pos}}$. In the following theorem we give asymptotic bounds in probability and identify the order of magnitude of the expected value of $\tau^\bold m_{\mathcal X^s_\text{pos}}$. Moreover, we identify the exponent at which the mixing time of the Markov chain $\{X_t^\beta\}_{t\in\mathbb N}$ asymptotically grows as $\beta$ and give an upper and a lower bound for the spectral gap, see  see \eqref{mixingtimedef} and \eqref{spectralgapdef}.
\begin{theorem}[Asymptotic behavior of $\tau_{\mathcal X^s_{\text{pos}}}^\bold m$ and mixing time]\label{timepositive}
Consider the $q$-state Potts model on a $K\times L$ grid $\Lambda$, with periodic boundary conditions and with positive external magnetic field. Then, for any $\bold m\in\mathcal X^m_{\emph{pos}}$, the following statements hold:
\begin{itemize}
\item[\emph{(a)}] for every $\epsilon>0$, $\lim_{\beta\to\infty}\mathbb P_\beta(e^{\beta(\Gamma^m_\emph{pos}-\epsilon)}<\tau_{\mathcal X^s_{\emph{pos}}}^{\bold m}<e^{\beta(\Gamma^m_\emph{pos}+\epsilon)})=1$;
\item[\emph{(b)}] $\lim_{\beta\to\infty} \frac{1}{\beta}\log\mathbb E[\tau_{\mathcal X^s_{\emph{pos}}}^\bold m]=\Gamma^m_\emph{pos}$;
\item[\emph{(c)}] for every $\epsilon\in(0,1)$, $\lim_{\beta\to\infty}\frac 1 \beta \log t^{\text{mix}}_\beta(\epsilon)=\Gamma_\emph{pos}^m$ and there exist two constants $0<c_1\le c_2<\infty$ independent of $\beta$ such that, for any $\beta>0$, $c_1e^{-\beta\Gamma^m_\emph{pos}}\le\rho_\beta\le c_2e^{-\beta\Gamma_\emph{pos}^m}$.
\end{itemize}
\end{theorem}
\textit{Proof.} Items (a) and (b) follow by \cite[Theorem 4.1]{manzo2004essential} and by \cite[Theorem 4.9]{manzo2004essential}, respectively, with $\eta_0=\bold m$ and $\Gamma=\Gamma^m_\text{pos}$ (using Theorem \ref{teometastablepos}). Item (c) follows by Corollary \ref{corollarygammatildepos} and by \cite[Proposition 3.24]{nardi2016hitting}. $\qed$\\

In the literature there exist some model-independent results on the asymptotic rescaled distribution of the first hitting time from some $\eta\in\mathcal X$ to a certain target $G\subset\mathcal X$, see for instance \cite[Theorem 4.15]{manzo2004essential}, \cite[Theorem 2.3]{fernandez2015asymptotically}, \cite[Theorem 3.19]{nardi2016hitting}. Unfortunately none of these results are suitable for our scenario when $\eta\in\mathcal X^m_\text{pos}$ and $G=\mathcal X^s_\text{pos}$. This fact follows by the presence of multiple degenerate metastable states that implies the presence of other deep wells in $\mathcal X$ different from the initial cycle $\mathcal C^\bold m_{\mathcal X^s_\text{pos}}(\Gamma^m_\text{pos})$. Hence, we consider a different target and we investigate the asymptotic rescaled distribution of the first hitting time from a metastable state to the subset $G\subset\mathcal X$ setting
\begin{align}\label{deftargetG}
G=\mathcal X\backslash\mathcal C^\bold m_{\mathcal X^s_\text{pos}}(\Gamma^m_\text{pos}).
\end{align}
We defer to Subsection \ref{proofenergylandscape} for the proof of the following theorem.
\begin{theorem}\label{teotimetargetG}
Consider the $q$-state Potts model on a $K\times L$ grid $\Lambda$, with periodic boundary conditions and with positive external magnetic field. Let $\bold m\in\mathcal X^m_\emph{pos}$ and let $G$ as defined in \eqref{deftargetG}. Then, 
\begin{align}
\frac{\tau^\bold m_G}{\mathbb E[\tau^\bold m_G]}\overset{d}\to\emph{Exp}(1),\ as\ \beta\to\infty.
\end{align}
\end{theorem}
Note that by definition \eqref{deftargetG}, by Proposition \ref{stablesetposprop} and by Theorem \ref{teometastablepos}, we have $\mathscr F(G)=\mathcal X^s_\text{pos}$ and that the maximal stability level is $V(G)=\Gamma^m_\text{pos}$.
\subsection{Minimal gates for the metastable transition}\label{mainresgates}
A further goal is to identify the union of all minimal gates for the transition from any metastable state to the unique stable state $\mathcal X^s_\text{pos}=\{\bold 1\}$. In order to do this, for any $m\in S\backslash\{1\}$, let us define 
\begin{align}\label{gatexm1}
&\mathcal W_{\text{pos}}(\bold m,\mathcal X^s_{\text{pos}}):= \bar B_{\ell^*-1,\ell^*}^1(m,1)
\ \ \text{and}\ \ \mathcal{W'}_{\text{pos}}(\bold m,\mathcal X^s_{\text{pos}}):=\bar B_{\ell^*,\ell^*-1}^1(m,1).
\end{align}
We refer to Figure \ref{figureexamplebis}(b)--(c) for an example of configurations belonging respectively to $\mathcal{W'}_{\text{pos}}(\bold m,\mathcal X^s_{\text{pos}})$ and to $\mathcal W_{\text{pos}}(\bold m,\mathcal X^s_{\text{pos}})$. 
These sets are investigated in Subsection \ref{mingatespos}. In particular, in Proposition \ref{propgatepos} we show that $\mathcal W_{\text{pos}}(\bold m,\mathcal X^s_{\text{pos}})$ is a gate for the transition from any $\bold m\in X^m_\text{pos}$ to $\mathcal X^s_\text{pos}$.

 Furthermore, in Subsection \ref{proofgates} we prove the following result.
\begin{theorem}[Minimal gates for the transition from $\bold m\in\mathcal X^m_\text{pos}$ to $\mathcal X^s_\text{pos}$]\label{teogateposset}
Consider the $q$-state Potts model on a $K\times L$ grid $\Lambda$, with periodic boundary conditions and with positive external magnetic field. Then, $\mathcal W_{\emph{pos}}(\bold m,\mathcal X^s_{\emph{pos}})$ is a minimal gate for the transition from any metastable state $\bold m\in\mathcal X^m_{\emph{pos}}$ to $\mathcal X^s_{\emph{pos}}=\{\bold 1\}$. Moreover,
\begin{align}\label{setmingatespos}
\mathcal G_{\emph{pos}}(\bold m,\mathcal X^s_{\emph{pos}})=\mathcal W_{\emph{pos}}(\bold m,\mathcal X^s_{\emph{pos}}).
\end{align}
\end{theorem}
\begin{figure}[h!]
\centering
\begin{tikzpicture}[scale=0.8,transform shape]
\draw[dotted,thick] (0,0.5) circle (0.8cm);
\fill[black!5!white] (0,0.5) circle (0.8cm);
\draw (0,2) circle (0.7cm);\draw (1.31,-0.24) circle (0.7cm);\draw (-1.31,-0.24) circle (0.7cm);
\fill[black!5!white] (0,2) circle (0.7cm);
\fill[black!5!white] (1.31,-0.24) circle (0.7cm);
\fill[black!5!white] (-1.31,-0.24) circle (0.7cm);
\fill  (0,0.5) circle (1.2pt);\draw (0,0.5) node[below] {{$\bold 1$}};
\fill (0,2) circle (1.2pt);\draw (0,2) node[below] {{$\bold 3$}};
\fill (1.31,-0.24) circle (1.2pt);\draw (1.31,-0.24) node[below] {{$\bold 4$}};
\fill (-1.31,-0.24) circle (1.2pt);\draw (-1.31,-0.24) node[below] {{$\bold 2$}};

\fill (-0.7,0.12) circle (1.5pt);
\fill (0.7,0.12) circle (1.5pt);
\fill (0,1.3) circle (1.5pt);
\draw[->,ultra thin] (-0.7,0.12) --(-1.2,-1.2);\draw (-1.4,-1.1) node[below] {$\bar B_{\ell^*-1,\ell^*}^1(2,1)$};
\draw[->,ultra thin] (0.7,0.12) --(1.2,-1.2);\draw (1.4,-1.1) node[below] {$\bar B_{\ell^*-1,\ell^*}^1(4,1)$};
\draw[->,ultra thin] (0,1.3)--(-1,1.3);;\draw (-0.9,1.3) node[left] {$\bar B_{\ell^*-1,\ell^*}^1(3,1)$};
\end{tikzpicture}
\caption{\label{figurepositiveabove} Viewpoint from above of the energy landscape depicted in Figure \ref{figurepositive}. For every $\bold m\in\mathcal X^m_\text{pos}$, the cycle whose bottom is the stable state $\bold 1$ is deeper than the initial cycles $\mathcal C^\bold m_{\mathcal X^s_\text{pos}}(\Gamma_\text{pos}^m)$. These last cycles are depicted with circles whose diameter is smaller than the one related to the stable state $\bold 1$.}
\end{figure}
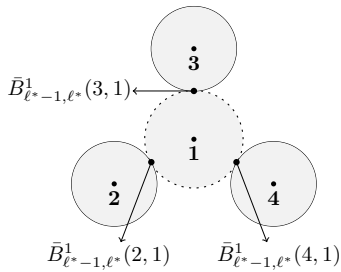\FloatBarrier
We remark that in \cite[Theorems 4.5 and 4.6]{bet2021metastabilityneg} the authors identify the union of all minimal gates for the metastable transitions for the $q$-Potts model with \textit{negative} external magnetic fields. These minimal gates have the same geometric definition of those of our scenario, the main difference is that in the negative case there are $q-1$ possible ``colors'' for the vertices inside the quasi-square with a unit protuberance.\\

In the next corollary we show that the process typically intersects $\mathcal W_{\text{pos}}(\bold m,\mathcal X^s_{\text{pos}})$ during the transition $\bold m\in\mathcal X^m_\text{pos}\to\mathcal X^s_\text{pos}$.
\begin{cor}\label{corpos}
Consider the $q$-state Potts model on a $K\times L$ grid $\Lambda$, with periodic boundary conditions and with positive external magnetic field. Then, for any $\bold m\in\mathcal X^m_{\emph{pos}}$
\begin{align}
\lim_{\beta\to\infty} \mathbb P_\beta(\tau^\bold m_{\mathcal W_{\emph{pos}}(\bold m,\mathcal X^s_{\emph{pos}})}<\tau^\bold m_{\mathcal X^s_{\emph{pos}}})=1.
\end{align}
\end{cor}
\textit{Proof}. The corollary follows from Proposition \ref{propgatepos} and from \cite[Theorem 5.4]{manzo2004essential}. $\qed$
\subsection{Minimal gates for the transition from a metastable state to the other metastable configurations}\label{mainresmetameta}
In Subsection \ref{subsectionnew} we study the transition from a metastable state to the set of the other metastable states. We prove that during the transition from any $\bold m\in\mathcal X^m_\text{pos}$ to $\mathcal X^s_\text{pos}$ almost surely the process does not intersect $\mathcal X^m_\text{pos}\backslash\{\bold m\}$, and we exploit this result to identificate the union of all minimal gates for this type of transition.
\begin{theorem}[Minimal gates for the transition from $\bold m$ to $\mathcal X^m_\text{pos}\backslash\{\bold m\}$]\label{mingatetranmetameta}
Let $\mathbf m\in\mathcal X^m_\emph{pos}$. If the external magnetic field is positive, then the following sets are minimal gates for the transition from $\mathbf m$ to $\mathcal X^m_\emph{pos}\backslash\{\mathbf m\}$
\begin{align}
\emph{(a)}&  \ \mathcal W_\emph{pos}(\mathbf m,\mathcal X^s_\emph{pos}),\hspace{90mm}\\
\emph{(b)}&\bigcup_{\mathbf z\in\mathcal X^m_\emph{pos}\backslash\{\mathbf m\}} \mathcal W_\emph{pos}(\mathbf z,\mathcal X^s_\emph{pos}).
\end{align}
Furthermore,
\begin{align}\label{Gmingatemetameta}
\mathcal G_\emph{pos}(\mathbf m,\mathcal X^m_\emph{pos}\backslash\{\mathbf m\})=\bigcup_{\mathbf z\in\mathcal X^m_\emph{pos}} \mathcal W_\emph{pos}(\mathbf z,\mathcal X^s_\emph{pos}).
\end{align}
\end{theorem}
We defer the proof of the above theorem in Subsection \ref{proofgates}. Note that in the negative scenario \cite{bet2021metastabilityneg} the theorem corresponding to Theorem \ref{mingatetranmetameta} is not present since there is only one metastable state.
\begin{cor}\label{corposmetameta}
If the external magnetic field is positive, then for any $\bold m\in\mathcal X^m_{\emph{pos}}$,
\begin{itemize}
\item[\emph{(a)}] $\lim_{\beta\to\infty} \mathbb P_\beta(\tau^\bold m_{\mathcal W_{\emph{pos}}(\bold m,\mathcal X^s_{\emph{pos}})}<\tau^\bold m_{\mathcal X^m_{\emph{pos}}\backslash\{\bold m\}})=1,$
\item[\emph{(b)}] $\lim_{\beta\to\infty} \mathbb P_\beta(\tau^\bold m_{\bigcup_{\mathbf z\in\mathcal X^m_\emph{pos}\backslash\{\mathbf m\}} \mathcal W_\emph{pos}(\mathbf z,\mathcal X^s_\emph{pos})})<\tau^\bold m_{\mathcal X^m_{\emph{pos}}\backslash\{\bold m\}})=1$.
\end{itemize}
\end{cor}
\textit{Proof}. The corollary follows from Propositions \ref{propgateposmetameta} and \ref{unionWmingatemetameta} and from \cite[Theorem 5.4]{manzo2004essential}. $\qed$

\subsection{Tube of typical trajectories of the metastable transition}\label{mainrestube}

\subsubsection{Further model-independent and model-dependent definitions}
In addition to the list of Subection \ref{modinddef}, in order state the main result concerning the tube of the typical trajectories we give some further definitions that are taken from \cite{nardi2016hitting}, \cite{cirillo2015metastability} and \cite{olivieri2005large}.
\begin{itemize}
\item[-] We call \textit{cycle-path} a finite sequence $(\mathcal C_1,\dots,\mathcal C_m)$ of trivial or non-trivial cycles $\mathcal C_1,\dots,\mathcal C_m\in\mathscr C(\mathcal X)$, such that $\mathcal C_i\cap\mathcal C_{i+1}=\varnothing\ \ \text{and}\ \ \partial\mathcal C_i\cap\mathcal C_{i+1}\neq\varnothing,\ \text{for every}\ i=1,\dots,m-1.$ 
\item[-] A cycle-path  $(\mathcal C_1,\dots,\mathcal C_m)$ is said to be \textit{downhill} (\textit{strictly downhill}) if the cycles $\mathcal C_1,\dots,\mathcal C_m$ are pairwise connected with decreasing height, i.e., when $H(\mathscr F(\partial\mathcal C_i))\ge H(\mathscr F(\partial\mathcal C_{i+1}))$ ($H(\mathscr F(\partial\mathcal C_i))> H(\mathscr F(\partial\mathcal C_{i+1}))$) for any $i=0,\dots,m-1$.
\item[-] For any $\mathcal C\in\mathscr C(\mathcal X)$, we define as
\begin{align}\label{principalboundary}
\mathcal B(\mathcal C):=
\begin{cases}
\mathscr F(\partial\mathcal C) &\text{if}\ \mathcal C\ \text{is a non-trivial cycle},\\
\{\eta\in\partial\mathcal C: H(\eta)< H(\sigma)\} &\text{if}\ \mathcal C=\{\sigma\}\ \text{is a trivial cycle},
\end{cases}
\end{align}
the \textit{principal boundary} of $\mathcal C$. Furthermore, let $\partial^{np}\mathcal C$ be the \textit{non-principal boundary} of $\mathcal C$, i.e., $\partial^{np}\mathcal C:=\partial\mathcal C\backslash\mathcal B(\mathcal C).$
\item[-] The \textit{relevant cycle} $\mathcal C_{\mathcal{A}}^+(\sigma)$ is 
\begin{align}\label{relevantcycle}
\mathcal C_{\mathcal{A}}^+(\sigma):=\{\eta\in\mathcal{X}:\ \Phi(\sigma,\eta)<\Phi(\sigma,\mathcal A)+\delta/2\},
\end{align} 
where $\delta$ is the minimum energy gap between any optimal and any non-optimal path from $\sigma$ to $\mathcal A$.
\item[-]  We denote the set of cycle-paths that lead from $\sigma$ to $\mathcal A$ and consist of maximal cycles  in $\mathcal X\backslash\mathcal A$ as 
\begin{align}
\hspace{0pt}\mathcal {P}_{\sigma,\mathcal A}:=\{\text{cycle-path}\ (\mathcal C_1,\dots,\mathcal C_m)|\ \mathcal C_1,\dots,\mathcal C_m&\in\mathcal M(\mathcal C_{\mathcal{A}}^+(\sigma)\backslash A),\sigma\in\mathcal C_1, \partial\mathcal C_m\cap\mathcal A\neq\varnothing\}.\notag
\end{align}
\item[-]  Given a non-empty set $\mathcal{A}\subset\mathcal X$ and $\sigma\in\mathcal{X}$, we constructively define a mapping $G: \Omega_{\sigma,A}\to\mathcal P_{\sigma,\mathcal A}$ in the following way. Given $\omega=(\omega_1,\dots,\omega_n)\in\Omega_{\sigma,A}$, we set $m_0=1$, $\mathcal C_1=\mathcal C_{\mathcal A}(\sigma)$ and define recursively $m_i:=\min\{k>m_{i-1}|\ \omega_k\notin\mathcal C_i\}$ and $\mathcal C_{i+1}:=\mathcal C_{\mathcal A}(\omega_{m_i})$. We note that $\omega$ is a finite sequence and $\omega_n\in\mathcal A$, so there exists an index $n(\omega)\in\mathbb N$ such that $\omega_{m_{n(\omega)}}=\omega_n\in\mathcal A$ and there the procedure stops. By  $(\mathcal C_1,\dots,\mathcal C_{m_{n(\omega)}})$ is a cycle-path with $\mathcal C_1,\dots,\mathcal C_{m_{n(\omega)}}\subset\mathcal M(\mathcal X\backslash\mathcal A)$. Moreover, the fact that $ \omega\in\Omega_{\sigma,A}$ implies that $\sigma\in\mathcal C_1$ and that $\partial\mathcal C_{n(\omega)}\cap\mathcal A\neq\varnothing$, hence $G(\omega)\in\mathcal P_{\sigma,\mathcal A}$ and the mapping is well-defined.

\item[-] We say that a cycle-path $(\mathcal C_1,\dots,\mathcal C_m)$ is \textit{connected via typical jumps} to $\mathcal A\subset\mathcal X$ or simply $vtj-$\textit{connected} to $\mathcal A$ if
\begin{align}\label{cyclepathvtj}
\mathcal B(\mathcal C_i)\cap\mathcal C_{i+1}\neq\varnothing,\ \ \forall i=1,\dots,m-1,\ \ \text{and}\ \ \mathcal B(\mathcal C_m)\cap\mathcal A\neq\varnothing.\end{align}
Let $J_{\mathcal C,\mathcal A}$ be the collection of all cycle-paths $(\mathcal C_1,\dots,\mathcal C_m)$ that are vtj-connected to $\mathcal A$ and such that $\mathcal C_1=\mathcal C$. 
\item[-] Given a non-empty set $\mathcal{A}$ and $\sigma\in\mathcal{X}$, we define $\omega\in\Omega_{\sigma,A}$ as a \textit{typical path} from $\sigma$ to $\mathcal A$ if its corresponding cycle-path $G(\omega)$ is vtj-connected to $\mathcal A$ and we denote by $\Omega_{\sigma,A}^{\text{vtj}}$ the collection of all typical paths from $\sigma$ to $\mathcal A$, i.e., 
\begin{align}\label{defvtjpaths}
\Omega_{\sigma,A}^{\text{vtj}}:=\{\omega\in\Omega_{\sigma,A}|\ G(\omega)\in J_{\mathcal C_{\mathcal A}(\sigma),\mathcal A}\}.
\end{align}
\item[-] We define the \textit{tube of typical paths} $T_{\mathcal A}(\sigma)$ from $\sigma$ to $\mathcal A$ as the subset of states $\eta\in\mathcal X$ that can be reached from $\sigma$ by means of a typical path which does not enter $\mathcal A$ before visiting $\eta$, i.e.,
\begin{align}
T_{\mathcal A}(\sigma):=\{\eta\in\mathcal X|\ \exists\omega\in\Omega_{\sigma,A}^{\text{vtj}}:\ \eta\in\omega\}.
\end{align}
Moreover, we define $\mathfrak T_{\mathcal A}(\sigma)$ as the set of all maximal cycles that belong to at least one vtj-connected path from $\mathcal C_{\mathcal A}^\sigma(\Gamma)$ to $\mathcal A$, i.e.,
\begin{align}\label{tubostorto}
\mathfrak T_{\mathcal A}(\sigma):=\{\mathcal C\in\mathcal M(\mathcal C_{\mathcal A}^+(\sigma)\backslash\mathcal A)|\ \exists(\mathcal C_1,\dots,\mathcal C_n)\in J_{\mathcal C_{\mathcal A}^\sigma(\Gamma),\mathcal A}\ \text{and}\exists j\in\{1,\dots,n\}:\ \mathcal C_j=\mathcal C\}.
\end{align}
Note that 
\begin{align}\label{reltubodrittoetorto}
\mathfrak T_{\mathcal A}(\sigma)=\mathcal M(T_{\mathcal A}(\sigma)\backslash\mathcal A)
\end{align}
and that the boundary of $T_{\mathcal A}(\sigma)$ consists of states either in $\mathcal A$ or in the non-principal part of the boundary of some $\mathcal C\in\mathfrak T_{\mathcal A}(\sigma)$:
\begin{align}
\partial T_{\mathcal A}(\sigma)\backslash\mathcal A\subseteq \bigcup_{\mathcal C\in\mathfrak T_{\mathcal A}(\sigma)} (\partial\mathcal C\backslash\mathcal B(\mathcal C))=:\partial^{np}\mathfrak T_{\mathcal A}(\sigma).
\end{align}
For the sake of semplicity, we will also refer to $\mathfrak T_{\mathcal A}(\sigma)$ as tube of typical paths from $\sigma$ to $\mathcal A$.
\end{itemize}
Furthermore, in addition to the list given in Subsection \ref{moddepdef}, we give some further model-dependent definitions.
\begin{itemize}
\item[-] For any $m,s\in S, m\neq s$, we define $\mathcal X(m,s)=\{\sigma\in\mathcal X: \sigma(v)\in\{m,s\}\ \text{for any}\ v\in V\}$. 
\item[-] For any $m\in S\backslash\{1\}$, we define
\begin{align}\label{setverticalstrip}
\mathscr S_\text{pos}^v(m,1)&:=\{\sigma\in\mathcal X(m,1): \sigma\ \text{has a vertical $1$-strip of thickness at least $\ell^*$ with }\notag\\&\text{possibly a bar of length $l=1,\dots,K$ on one of the two vertical edges}\},\\ 
\label{sethorizontalstrip}
\mathscr S_\text{pos}^h(m,1)&:=\{\sigma\in\mathcal X(m,1): \sigma\ \text{has a horizontal $1$-strip of thickness at least $\ell^*$ with}\notag\\&\text{possibly a bar of length $l=1,\dots,L$ on one of the two horizontal edges}\}.
\end{align} 
\end{itemize}
\subsubsection{Main results on the tube of typical trajectories}
In this subsection we give our main result concerning the tube of the typical trajectories for the transition $\mathbf m\to\mathcal X^s_\text{pos}$ for any fixed $\mathbf m\in\mathcal X^m_\text{pos}$. The tube of typical paths for this transition turns out to be
\begin{align}\label{tubepos}
\mathfrak T&_{\mathcal X^s_\text{pos}}(\mathbf m):=\bigcup_{\ell=1}^{\ell^*-1} \bar R_{\ell-1,\ell}(m,1)\cup\bigcup_{\ell=1}^{\ell^*} \bar R_{\ell-1,\ell-1}(m,1)\cup\bigcup_{\ell=1}^{\ell^*-1}\bigcup_{l=1}^{\ell-1} \bar B^l_{\ell-1,\ell}(m,1)\notag\\
&\cup\bigcup_{\ell=1}^{\ell^*}\bigcup_{l=1}^{\ell-2} \bar B^l_{\ell-1,\ell-1}(m,1)\cup\bar B_{\ell^*-1,\ell^*}^1(m,1)\cup\bigcup_{\ell_1=\ell^*}^{K-1}\bigcup_{\ell_2=\ell^*}^{K-1} \bar R_{\ell_1,\ell_2}(m,1)\cup\bigcup_{\ell_1=\ell^*}^{K-1}\bigcup_{\ell_2=\ell^*}^{K-1}\bigcup_{l=1}^{\ell_2-1}\bar B_{\ell_1,\ell_2}^l(m,1)\notag\\
&\cup\bigcup_{\ell_1=\ell^*}^{L-1}\bigcup_{\ell_2=\ell^*}^{L-1} \bar R_{\ell_1,\ell_2}(m,1)\cup\bigcup_{\ell_1=\ell^*}^{L-1}\bigcup_{\ell_2=\ell^*}^{L-1}\bigcup_{l=1}^{\ell_2-1}\bar B_{\ell_1,\ell_2}^l(m,1)\cup\mathscr S_\text{pos}^v(m,1)\cup\mathscr S_\text{pos}^h(m,1).
\end{align}
As illustrated in the next result, which we prove in Section \ref{prooftube}, $\mathfrak T_{\mathcal X^s_\text{pos}}(\mathbf m)$ includes those configurations with a positive probability of being visited by the Markov chain $\{X_t\}^\beta_{t\in\mathbb N}$ started in $\mathbf m$ before hitting $\mathcal X^s_\text{pos}$ in the limit $\beta\to\infty$. Note that the relation between $T_{\mathcal X^s_\text{pos}}(\mathbf m)$ and $\mathfrak T_{\mathcal X^s_\text{pos}}(\mathbf m)$ is given by \eqref{reltubodrittoetorto}.

\begin{theorem}\label{teotubesetpos}
If the external magnetic field is positive, then for any $\mathbf m\in\mathcal X^m_\emph{pos}$ the tube of typical trajectories for the transition $\mathbf m\to\mathcal X^s_\emph{pos}$ is \eqref{tubepos}. Furthermore,
\begin{align}\label{ristimetubenzbneg}
\mathbbm P_\beta(\tau_{\partial^{np}\mathfrak T_{\mathcal X^s_\emph{pos}}(\mathbf m)}^\mathbf m\le\tau_{\mathcal X^s_\emph{pos}}^\mathbf m)\le e^{-k\beta}.
\end{align}
\end{theorem}
\section{Energy landscape analysis and asymptotic behavior}\label{energylandscape}
In this section we analyse the enrgy landscape of the $q$-state Potts model with positive external magnetic field. 
\noindent First we recall some useful definitions and lemmas from \cite{nardi2019tunneling}. 
\subsection{Known local geometric properties}\label{subseclocgeo}
In the following list we introduce the notions of \textit{disagreeing edges}, \textit{bridges} and \textit{crosses} of a Potts configuration on a grid-graph $\Lambda$.
\begin{itemize}
\item[-] We call $e=(v,w)\in E$ a \textit{disagreeing edge} if it connects two vertices with different spin values, i.e., $\sigma(v)\neq\sigma(w)$.
\item[-] For any $i=0,\dots,K-1$, let
\begin{align}\label{diseedgesrow}
d_{ r_i}(\sigma):=\sum_{(v,w)\in r_i} \mathbbm{1}_{\{\sigma(v)\neq\sigma(w)\}}
\end{align} 
be the total number of disagreeing edges that $\sigma$ has on row $r_i$. Furthermore, for any $j=0,\dots,L-1$ let
\begin{align}\label{disedgescol}
d_{c_j}(\sigma):=\sum_{(v,w)\in c_j} \mathbbm{1}_{\{\sigma(v)\neq\sigma(w)\}},
\end{align} 
be the total number of disagreeing edges that $\sigma$ has on column $c_j$.
\item[-] We define $d_h(\sigma)$ as the total number of disagreeing horizontal edges and $d_v(\sigma)$ as the total number of disagreeing vertical edges, i.e.,
\begin{align}\label{disedgeshorizontal}
d_h(\sigma):=\sum_{i=0}^{K-1} d_{ r_i}(\sigma),
\end{align} and
\begin{align}\label{disedgesvertical}
d_v(\sigma):=\sum_{j=0}^{L-1} d_{c_j}(\sigma).
\end{align}
Since we may partition the edge set $E$ in the two subsets of horizontal edges $E_h$ and of vertical edges $E_v$, such that $E_h\cap E_v=\varnothing$, the total number of disagreeing edges is given by
\begin{align}\label{totnumberdisedges}
\sum_{(v,w)\in E} \mathbbm{1}_{\{\sigma(v)\neq\sigma(w)\}}&=
\sum_{(v,w)\in E_v} \mathbbm{1}_{\{\sigma(v)\neq\sigma(w)\}}+\sum_{(v,w)\in E_h} \mathbbm{1}_{\{\sigma(v)\neq\sigma(w)\}}=d_v(\sigma)+d_h(\sigma).
\end{align} 
\item[-] We say that $\sigma$ has a \textit{horizontal bridge} on row $r$ if $\sigma(v)=\sigma(w),\ \text{for all}\ v,w\in r$.
\item[-] We say that $\sigma$ has a \textit{vertical bridge} on column $c$ if $\sigma(v)=\sigma(w),\ \text{for all}\ v,w\in c$. 
\item[-] We say that $\sigma\in\mathcal X$ has a \textit{cross} if it has at least one vertical and one horizontal bridge.
\end{itemize}
For sake of semplicity, if $\sigma$ has a bridge of spins $s\in S$, then we say that $\sigma$ has an $s$-bridge. Similarly, if $\sigma$ has a cross of spins $s$, we say that $\sigma$ has an $s$-cross.
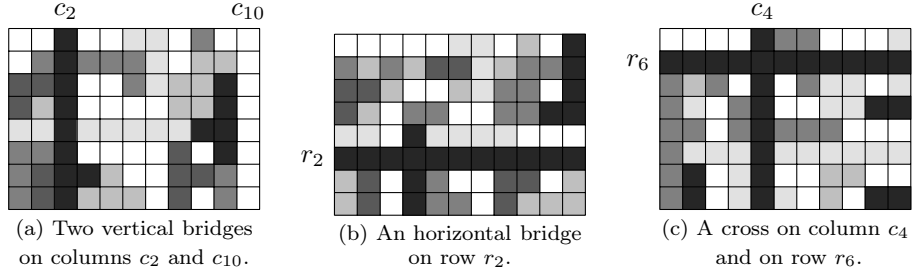
\begin{figure}[h!]
\centering
\begin{tikzpicture}[scale=1, transform shape]
\draw [fill=gray,gray] (0,0) rectangle (0.3,0.9) (0.3,0.6) rectangle (0.6,0.9) (2.7,0) rectangle (3,0.9) (2.4,2.1) rectangle (2.7,2.4) (0.6,0.9) rectangle (0.9,0.6) (0.9,1.8) rectangle (1.8,2.1) rectangle (1.5,1.5) (0.3,1.8)rectangle(0.6,2.1);
\draw [fill=gray,lightgray] (0.9,0) rectangle (1.8,0.3) (1.2,0.3) rectangle (1.5,0.6) (2.1,1.2) rectangle (2.7,1.8) rectangle (3,2.1) (2.4,1.8) rectangle (2.7,2.1) (0.3,1.2)rectangle(0.6,1.5);
\draw [fill=gray,black!85!white] (0.6,0) rectangle (0.9,2.4) (0.9,0.3) rectangle (1.2,0.6) (2.7,0.6) rectangle (3,1.8) (2.4,0.9) rectangle (3,1.2);
\draw [fill=gray,black!12!white] (0,0.9)rectangle(0.6,1.2) (0.9,0.9)rectangle(2.4,1.2)(1.8,1.5)rectangle(2.1,2.4)(1.5,2.1)rectangle(1.8,2.4);
\draw [fill=gray,black!65!white] (0.3,0)rectangle(0.6,0.6) (0,1.2)rectangle(0.3,1.5)(0,1.5)rectangle(0.6,1.8)(2.1,0)rectangle(2.4,0.9)(2.4,0.3)rectangle(2.7,0.6);
\draw[step=0.3cm,color=black] (0,0) grid (3.3,2.4);
\draw (0.75, 2.4) node[above] {$c_2$};\draw (3.15, 2.4) node[above] {$c_{10}$};
\draw (1.65,0) node[below] {\footnotesize{(a) Two vertical bridges}};\draw(1.65,-0.4) node[below]{\footnotesize{on columns $c_2$ and $c_{10}$.}};
\end{tikzpicture}\ \
\begin{tikzpicture}[scale=1, transform shape]
\draw [fill=gray,gray] (0,0) rectangle (0.3,2.1) (0.9,1.2)rectangle(1.5,1.5)(0.6,1.5)rectangle(0.9,2.1)(0.9,0) rectangle (1.8,0.3) (1.2,0.3) rectangle (1.5,0.6) (2.1,1.2) rectangle (2.7,1.8) rectangle (3,2.1) (2.4,1.8) rectangle (2.7,2.1) (0.3,1.2)rectangle(0.6,1.5);
\draw [fill=gray,lightgray] (0.3,0) rectangle (0.6,1.5)(2.1,1.8)rectangle(2.4,2.1)(2.1,0)rectangle(3.3,0.3)(0,1.2)rectangle(0.3,1.5)(0,1.5)rectangle(0.9,1.8)(0,0) rectangle (0.3,0.9) (0.3,0.6) rectangle (0.6,0.9) (3,0) rectangle (3.3,0.9) (2.4,2.1) rectangle (2.7,2.4) (0.6,0.9) rectangle (0.9,0.6) (0.9,1.8) rectangle (1.8,2.1) rectangle (1.5,1.5) (0.3,1.8)rectangle(0.6,2.1);;
\draw [fill=gray,black!12!white] (0,0.9)rectangle(0.6,1.2) (0.9,0.9)rectangle(2.4,1.2)(1.8,1.5)rectangle(2.1,2.4)(1.5,2.1)rectangle(1.8,2.4) (0.3,0)rectangle(0.6,0.6) (0,1.2)rectangle(0.3,1.5);
\draw [fill=gray,black!65!white] (0.3,0)rectangle(0.6,0.9) (0,1.2)rectangle(0.3,1.5)(0,1.5)rectangle(0.6,1.8)(2.1,0)rectangle(2.4,0.9)(2.4,0.3)rectangle(2.7,0.6)(1.2,1.8)rectangle(1.8,2.1);
\draw [fill=gray,black!85!white] (0,0.6)rectangle(3.3,0.9)(2.7,1.2)rectangle(3.3,1.5)rectangle(3,2.4) (0.9,0)rectangle(1.2,1.2);
\draw[step=0.3cm,color=black] (0,0) grid (3.3,2.4);
\draw (0,0.75) node[left] {$r_2$};
\draw (1.65,0) node[below] {\footnotesize{(b) An horizontal bridge}};\draw(1.65,-0.4) node[below]{\footnotesize{on row $r_2$.}};
\end{tikzpicture} \ \
\begin{tikzpicture}[scale=1, transform shape]
\draw [fill=gray,lightgray] (0.3,0)rectangle(0.6,0.6) (0,1.2)rectangle(0.3,1.5)(0,1.5)rectangle(0.6,1.8)(2.1,0)rectangle(2.4,0.9)(2.4,0.3)rectangle(2.7,0.6)(1.2,1.8)rectangle(1.8,2.1);
\draw [fill=gray,black!12!white] (0,0) rectangle (0.3,2.1) (0.9,1.2)rectangle(1.5,1.5)(0.6,1.5)rectangle(0.9,2.1)(0.9,0) rectangle (1.8,0.3) (1.2,0.3) rectangle (1.5,0.6) (2.1,1.2) rectangle (2.7,1.8) rectangle (3,2.1) (2.4,1.8) rectangle (2.7,2.1) (0.3,1.2)rectangle(0.6,1.5) (0,0.6)rectangle(3.3,0.9)(2.7,1.2)rectangle(3.3,1.5)rectangle(3,2.4) (0.9,0)rectangle(1.2,1.2);
\draw [fill=gray,gray] (0,0) rectangle (0.3,2.1) (0.9,1.2)rectangle(1.5,1.5)(0.6,1.5)rectangle(0.9,2.1)(0,0.9)rectangle(0.6,1.2) (0.9,0.9)rectangle(2.4,1.2)(1.8,1.5)rectangle(2.1,2.4)(1.5,2.1)rectangle(1.8,2.4);
\draw [fill=gray,black!85!white] (1.2,0) rectangle (1.5,2.4) (0.3,0.3) rectangle (0.6,0.6) (0,1.8)rectangle(3.3,2.1)(2.7,1.2)rectangle(3.3,1.5) (0.3,0)rectangle(0.6,0.6) (2.7,0)rectangle(3.3,0.3);
\draw[step=0.3cm,color=black] (0,0) grid (3.3,2.4);
\draw (0,1.95) node[left] {$r_6$};\draw (1.35, 2.4) node[above] {$c_4$};
\draw (1.65,0) node[below] {\footnotesize{(c) A cross on column $c_4$}};\draw(1.65,-0.4) node[below]{\footnotesize{ and on row $r_6$.}};
\end{tikzpicture}
\caption{\label{esbridgescross} Example of configurations on a $8\times 11$ grid graph displaying a vertical $s$-bridge (a), a horizontal $s$-bridge (b) and a $s$-cross (c). We color black the spins $s$.}
\end{figure}\FloatBarrier
\begin{itemize}
\item[-] For any $s\in S$, the total number of $s$-bridges of the configuration $\sigma$ is denoted by $B_s(\sigma)$. 
\end{itemize}
Note that if a configuration $\sigma\in\mathcal X$ has an $s$-cross, then $B_s(\sigma)$ is at least $2$ since the presence of an $s$-cross implies the presence of two $s$-bridges, i.e., of a horizontal $s$-bridge and of a vertical $s$-bridge. \\
We conclude this section by recalling the following three useful lemmas from \cite{nardi2019tunneling}. These results give us some geometric properties for the $q$-state Potts model on a grid-graph and they are verified regardless of the definition of the external magnetic field. 
\begin{lemma}\emph{\cite[Lemma 2.2]{nardi2019tunneling}}\label{firstlemmaNZ}
A Potts configuration on a grid-graph $\Lambda$ does not have simultaneously a horizontal bridge and a vertical bridge of different spins.
\end{lemma}
\begin{lemma}\emph{\cite[Lemma 2.6]{nardi2019tunneling}}\label{numberbridges}
Let $\sigma, \sigma'\in\mathcal X$ be two Potts configurations which differ by a single-spin update, that is $|\{v\in V: \sigma(v)\neq\sigma'(v)\}|=1$. Then for every $s\in S$ we have that 
\begin{itemize}
\item[\emph{(i)}] $B_s(\sigma')-B_s(\sigma)\ \in \{-2,-1,0,1,2\}$,
\item[\emph{(ii)}] $B_s(\sigma')-B_s(\sigma)=2$ \ \text{if and only if $\sigma'$ has an $s$-cross that $\sigma$ does not have}.
\end{itemize}
\end{lemma}
\begin{lemma}\emph{\cite[Lemma 2.3]{nardi2019tunneling}}\label{gapbridges}
The following properties hold for every Potts configuration $\sigma\in\mathcal X$ on a grid graph $\Lambda$ with periodic boundary conditions:
\begin{itemize}
\item[\emph{(i)}] $d_{ r}(\sigma)$=0 if and only if $\sigma$ has a horizontal bridge on row $r$;
\item[\emph{(ii)}] $d_c(\sigma)$=0 if and only if $\sigma$ has a vertical bridge on column $c$;
\item[\emph{(iii)}] if $\sigma$ has no horizontal bridge on row $r$, then $d_{ r}(\sigma)\ge 2$;
\item[\emph{(iv)}] if $\sigma$ has no vertical bridge on column $c$, then $d_c(\sigma)\ge 2.$
\end{itemize}
\end{lemma}
\subsection{Metastable states and stability level of the metastable configurations}\label{stablevpos}
By Proposition \ref{stablesetposprop}, $H_{\text{pos}}$ has only one global minimum, $\mathcal X^s_{\text{pos}}=\{\bold 1\}$. Furthermore, the configurations $\bold 2,\dots,\bold q$ are such that $H_{\text{pos}}(\bold 2)=\dots=H_{\text{pos}}(\bold q)$. In this subsection, our aim is to prove that the metastable set $\mathcal X^m_{\text{pos}}$ is the union of these configurations. We are going to prove this claim by steps. We begin by obtaining an upper bound for the stability level of the states $\bold 2,\dots,\bold q$. 

\noindent Given $\bold m\in\{\bold 2,\dots,\bold q\}$, let us compute the following energy gap between any $\sigma\in\mathcal X$ and $\bold m$,
\begin{align}\label{rewritegap1pos}
H_{\text{pos}}(\sigma)-H_{\text{pos}}(\bold m)&=-\sum_{(v,w)\in E} \mathbbm{1}_{\{\sigma(v)=\sigma(w)\}} -h\sum_{u\in V} \mathbbm{1}_{\{\sigma(u)=1\}}-(-|E|) \notag \\
&=\sum_{(v,w)\in E} \mathbbm{1}_{\{\sigma(v)\neq\sigma(w)\}}-h\sum_{u\in V} \mathbbm{1}_{\{\sigma(u)=1\}} \notag \\
&=d_v(\sigma)+d_h(\sigma)-h\sum_{u\in V} \mathbbm{1}_{\{\sigma(u)=1\}},
\end{align}
where in the last equality we used \eqref{totnumberdisedges}.
\begin{definition}\label{refpathmiopos}
For any $\bold m\in\mathcal X^m_{\text{pos}}$, we define a \textit{reference path} $\tilde\omega:\bold m\to\bold 1$, $\tilde\omega=(\tilde\omega_0,\dots,\tilde\omega_{KL})$ as the concatenation of the two paths ${\tilde\omega}^{(1)}:=(\bold 1=\tilde\omega_0,\dots,\tilde\omega_{(K-1)^2})$ and ${\tilde\omega}^{(2)}:=(\tilde\omega_{(K-1)^2},\dots,$ $\bold m=\tilde\omega_{KL})$. The paths ${\tilde\omega}^{(1)}$ and ${\tilde\omega}^{(2)}$ are obtained by replacing $\bold 1$ with $\bold m$ and $\bold s$ with $\bold 1$ in the paths $\hat{\omega}^{(1)}$ and $\hat{\omega}^{(2)}$ of \cite[Definition 5.1]{bet2021metastabilityneg}. See Appendix \ref{appendixproofposdef} for the explicit definition.
\end{definition}

For any fixed $\bold m\in\mathcal X^m_\text{pos}$, let us focus on the transition from $\bold m$ to $\mathcal X^s_\text{pos}=\{\bold 1\}$. Given $m\in S$, let \begin{align}\label{ennekappa} N_m(\sigma):=|\{v\in V:\ \sigma(v)=m\}| \end{align} be the number of vertices with spin $m$ in $\sigma\in\mathcal X$.
\begin{lemma}\label{lemmaposuno}
Let $\bold m\in\mathcal X^m_\emph{pos}$. If the external magnetic field is positive, then for any $\sigma\in\bar R_{\ell^*-1,\ell^*}(m,1)$ there exists a path $\gamma:\sigma\to\bold m$ such that the maximum energy along $\gamma$ is bounded as
\begin{align}
\max_{\xi\in\gamma} H_\emph{pos}(\xi)<4\ell^*-h(\ell^*(\ell^*-1)+1)+H_\emph{pos}(\bold m).
\end{align}
\end{lemma}
\textit{Proof.} The proof proceeds analogously to the proof of \cite[Lemma 5.4]{bet2021metastabilityneg} by replacing $\bold 1$ with $\bold m$, $\bold s$ with $\bold 1$ and $\hat\omega$ with $\tilde\omega$. See Appendix \ref{appendixproofposlemma61} for the explicit proof. $\qed$\\

In the next lemma we show that for any $\bold m\in\mathcal X^m_\text{pos}$, $\bar B_{\ell^*-1,\ell^*}^2(m,1)$ is connected to the stable set $\mathcal X^s_\text{pos}$ by a path that does not overcome the energy value $4\ell^*-h(\ell^*(\ell^*-1)+1)+H_\text{pos}(\bold m)$.
\begin{lemma}\label{lemmaduepositive}
Let $\bold m\in\mathcal X^m_\emph{pos}$. If the external magnetic field is positive, then for any $\sigma\in\bar B_{\ell^*-1,\ell^*}^2(m,1)$ there exists a path $\gamma:\sigma\to\bold m$ such that the maximum energy along $\gamma$ is bounded as
\begin{align}
\max_{\xi\in\gamma} H_\emph{pos}(\xi)<4\ell^*-h(\ell^*(\ell^*-1)+1)+H_\emph{pos}(\bold m).
\end{align}
\end{lemma}
\textit{Proof.} The proof proceeds analogously to the proof of \cite[Lemma 5.5]{bet2021metastabilityneg} by replacing $\bold 1$ with $\bold m$, $\bold s$ with $\bold 1$ and $\hat\omega$ with $\tilde\omega$. See Appendix \ref{appendixproofposlemma62} for the explicit proof. $\qed$\\

We are now able to prove the following propositions, in which we give an upper bound and a lower bound for $\Gamma_\text{pos}(\bold m,\mathcal X^s_\text{pos}):=\Phi_\text{pos}(\bold m,\mathcal X^s_\text{pos})-H_\text{pos}(\bold m)$, for any $\bold m\in\mathcal X^m_\text{pos}$.
\begin{proposition}[Upper bound for the communication height]\label{refpathpos}
If the external magnetic field is positive, then for every $\bold m\in\mathcal X^m_{\emph{pos}}$, 
\begin{align}\label{upperboundposequation}
\Phi_{\text{\emph{pos}}}(\bold m,\mathcal X^s_\text{\emph{pos}})-H_{\text{\emph{pos}}}(\bold m)\le  4\ell^*-h(\ell^*(\ell^*-1)+1).
\end{align}
\end{proposition}
\textit{Proof.} The upper bound \eqref{upperboundposequation} follows by the proof of Lemma \ref{lemmaduepositive}, where we proved that $\max_{\xi\in\tilde\omega} H_\text{pos}(\xi)=H_\text{pos}(\tilde\omega_{k^*})= 4\ell^*-h(\ell^*(\ell^*-1)+1)+H_\text{pos}(\bold m)$.
$\qed$

\begin{proposition}[Lower bound for the communication height]\label{lowerboundpos}
If the external magnetic field is positive, then for every $\bold m\in\{\bold 2,\dots,\bold q\}$, 
\begin{align}\label{lowerboundaligniphipos}
\Phi_{\text{\emph{pos}}}(\bold m,\mathcal X^s_\text{\emph{pos}})-H_{\text{\emph{pos}}}(\bold m)\ge 4\ell^*-h(\ell^*(\ell^*-1)+1).\end{align} 
\end{proposition} 
\textit{Proof.} The proof proceeds analogously to the proof of \cite[Proposition 5.2]{bet2021metastabilityneg} by replacing $\bold 1$ with $\bold m$, $\bold s$ with $\bold 1$ and $\hat\omega$ with $\tilde\omega$. See Appendix \ref{appendixproofposprop62} for the explicit proof. $\qed$\\

The above Propositions \ref{refpathpos} and \ref{lowerboundpos} are used to prove \eqref{estimatestablevelmetapos} in Theorem \ref{teometastablepos}. Note that \eqref{estimatestablevelmetapos} is the min-max energy value reached by any optimal path $\omega\in\Omega_{\bold m,\mathcal X^s_\text{pos}}^{opt}$ for every $\bold m\in\mathcal X^m_\text{pos}$.
\begin{lemma}\label{lemmatrepospre}
If the external magnetic field is positive, then any $\omega\in\Omega_{\bold m,\mathcal X^s_\emph{pos}}^{opt}$ is such that $\omega\cap\bar R_{\ell^*-1,\ell^*}(m,1)\neq\varnothing$.
\end{lemma}
\textit{Proof.} The proof proceeds analogously to the proof of \cite[Lemma 5.6]{bet2021metastabilityneg} by replacing $\bold 1$ with $\bold m$, $\bold s$ with $\bold 1$ and $\hat\omega$ with $\tilde\omega$. See Appendix \ref{appendixproofposlemma63} for the explicit proof. $\qed$\\

Let $\sigma\in\mathcal X$ and let $v\in V$. We define the \textit{tile} centered in $v$, denoted by $v$-tile, as the set of five sites consisting of $v$ and its four nearest neighbors. See for instance Figure \ref{figmattonellepos}.
A $v$-tile is said to be \textit{stable} for $\sigma$ if by flipping the spin on vertex $v$ from $\sigma(v)$ to any $s\in S$ the energy difference $H_\text{pos}(\sigma^{v,s})-H_\text{pos}(\sigma)$ is greater than or equal to zero. 

\noindent In Lemma \ref{stabletilespositive} we define the set of all possible stable tiles induced by the Hamiltonian \eqref{hamiltonianpos} and next we exploit this result to prove Proposition \ref{proplocalminimapos}. For any $\sigma\in\mathcal X$, $v\in V$ and $s\in S$, we define $n_s(v)$ as the number of nearest neighbors to $v$ with spin $s$ in $\sigma$, i.e.,
\begin{align}\label{numberspinneighborhood}
n_s(v):=|\{w\in V:\ w\sim v,\ \sigma(w)=s\}|.
\end{align}
\begin{lemma}[Characterization of stable $v$-tiles in a configuration $\sigma$]\label{stabletilespositive}
Let $\sigma\in\mathcal X$ and let $v\in V$. If the external magnetic field is positive, then the tile centered in $v$ is stable for $\sigma$ if and only if it satisfies one of the following conditions.
\begin{itemize}
\item[(1)] For any $m\in S$, $m\neq 1$, if $\sigma(v)=m$, $v$ has either at least three nearest neighbors with spin $m$ or two nearest neighbors with spin $m$ and two nearest neighbors with spin $r,t\in S\backslash\{m\}$ such that they may be not both equal to $1$, see Figure \ref{figmattonellepos}\emph{(a),(c),(d),(f)--(m)}, or one nearest neighbor $m$ and three nearest neighbors with spin $r,s,t\in S\backslash\{1\}$ different from each other, see Figure \ref{figmattonellepos}\emph{(r)}.
\item[(2)] If $\sigma(v)=1$, $v$ has either at least two nearest neighbors with spin $1$, see Figure \ref{figmattonellepos}\emph{(b),(e),(n)--(q)} or it has one nearest neighbor $1$ and three nearest neighbors with spin $r,s,t\in S\backslash\{1\}$ different from each other, see Figure \ref{figmattonellepos}\emph{(s)}.
\end{itemize}
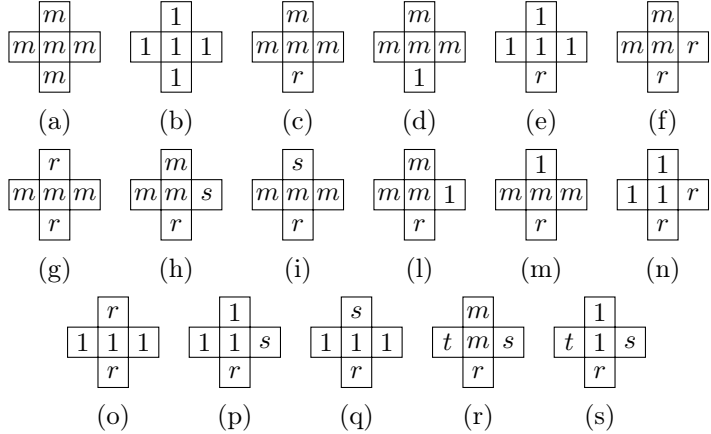
\begin{figure}[h!]
\centering
\begin{tikzpicture}
\foreach \i in {0,1.6,3.2,4.8,6.4,8}
\draw (\i,0)rectangle(\i+1.2,0.4);
\foreach \i in {0,1.6,3.2,4.8,6.4,8}
\draw (\i+0.4,-0.4)rectangle(\i+0.8,0.8);

\draw (0.6,-0.5) node[below] {(a)}(2.2,-0.5) node[below] {(b)}(3.8,-0.5) node[below] {(c)}(5.4,-0.5) node[below] {(d)}(7,-0.5) node[below] {(e)}(8.6,-0.5) node[below] {(f)};

\foreach \i in {0.6,3.8,5.4,8.6}\draw (\i,0.2) node {$m$} (\i,0.6) node {$m$};
\foreach \i in {2.2,7}\draw (\i,0.2) node {$1$}(\i,0.6) node {$1$};
\foreach \i in {2.2,5.4}\draw (\i,-0.2) node {$1$};
\foreach \i in {1.8,2.6,6.6,7.4}\draw  (\i,0.2) node {$1$}; 
\draw (0.6,-0.2) node {$m$}; 
\foreach \i in {0.2,3.4,5,5.8,8.2}\draw  (\i,0.2) node {$m$}; 
\foreach \i in {1,4.2}\draw  (\i,0.2) node {$m$}; 
\foreach \i in {3.8,7,8.6}\draw (\i,-0.2) node {$r$};
\draw (9,0.2) node {$r$}; 

\end{tikzpicture}
\begin{tikzpicture}
\draw (0.6,-0.5) node[below] {(h)}(2.2,-0.5) node[below] {(i)}(3.8,-0.5) node[below] {(l)}(5.4,-0.5) node[below] {(m)}(7,-0.5) node[below] {(n)};

\draw (-1,-0.5) node[below] {(g)};
\draw (-1.6,0)rectangle(-0.4,0.4)(-1.2,-0.4)rectangle(-0.8,0.8);
\foreach \i in {-1.4,-1,-0.6}\draw (\i,0.2) node {$m$};
\foreach \i in {-0.2,0.6}\draw (-1,\i) node {$r$};

\draw (0,0)rectangle(1.2,0.4)(0.4,-0.4)rectangle(0.8,0.8);
\foreach \i in {0.2,0.6}\draw (\i,0.2) node {$m$};
\draw (0.6,0.6) node {$m$};
\draw (0.6,-0.2) node {$r$};\draw (1,0.2) node {$s$};

\draw (1.6,0)rectangle(2.8,0.4)(2,-0.4)rectangle(2.4,0.8);
\foreach \i in {1.8,2.2,2.6}\draw (\i,0.2) node {$m$};
\draw (2.2,-0.2) node {$r$}; \draw (2.2,0.6) node {$s$}; 

\draw (3.2,0)rectangle(4.4,0.4)(3.6,-0.4)rectangle(4,0.8);
\foreach \i in {3.4,3.8}\draw (\i,0.2) node {$m$};
\draw (3.8,0.6) node {$m$};
\draw (3.8,-0.2) node {$r$};\draw (4.2,0.2) node {$1$};

\draw (4.8,0)rectangle(6,0.4)(5.2,-0.4)rectangle(5.6,0.8);
\foreach \i in {5,5.4,5.8}\draw (\i,0.2) node {$m$};
\draw (5.4,-0.2) node {$r$};\draw (5.4,0.6) node {$1$};

\draw (6.4,0)rectangle(7.6,0.4)(6.8,-0.4)rectangle(7.2,0.8);
\foreach \i in {6.6,7}\draw (\i,0.2) node {$1$};
\draw (7,0.6) node {$1$};
\draw (7,-0.2) node {$r$};\draw (7.4,0.2) node {$r$};

\end{tikzpicture}
\begin{tikzpicture}
\draw (-1,-0.5) node[below] {(p)} (0.6,-0.5) node[below] {(q)}(2.2,-0.5) node[below] {(r)}(3.8,-0.5) node[below] {(s)}(-2.6,-0.5) node[below] {(o)};

\draw (-3.2,0)rectangle(-2,0.4)(-2.8,-0.4)rectangle(-2.4,0.8);
\foreach \i in {-2.2,-3,-2.6}\draw (\i,0.2) node {$1$};
\draw (-2.6,-0.2) node {$r$};\draw (-2.6,0.6) node {$r$};

\draw (-1.6,0)rectangle(-0.4,0.4)(-1.2,-0.4)rectangle(-0.8,0.8);
\foreach \i in {-1.4,-1}\draw (\i,0.2) node {$1$};
\draw (-1,0.6) node {$1$};
\draw (-1,-0.2) node {$r$};\draw (-0.6,0.2) node {$s$};

\draw (0,0)rectangle(1.2,0.4)(0.4,-0.4)rectangle(0.8,0.8);
\foreach \i in {0.2,0.6,1}\draw (\i,0.2) node {$1$};
\draw (0.6,0.6) node {$s$};
\draw (0.6,-0.2) node {$r$};

\draw (1.6,0)rectangle(2.8,0.4)(2,-0.4)rectangle(2.4,0.8);
\draw (2.2,0.2) node {$m$}; \draw (2.2,0.6) node {$m$};
\draw (1.8,0.2) node {$t$};\draw (2.2,-0.2) node {$r$};\draw (2.6,0.2) node {$s$};

\draw (3.2,0)rectangle(4.4,0.4)(3.6,-0.4)rectangle(4,0.8);
\draw (3.8,0.2) node {$1$}; \draw (3.8,0.6) node {$1$};
\draw (3.4,0.2) node {$t$};\draw (3.8,-0.2) node {$r$};\draw (4.2,0.2) node {$s$};

\end{tikzpicture}
\caption{\label{figmattonellepos} Stable tiles centered in any $v\in V$ for a $q$-Potts configuration on $\Lambda$ for any $m,r,s,t\in S\backslash\{1\}$ different from each other. The tiles are depicted up to a rotation of $\alpha\frac\pi 2$, $\alpha\in\mathbb Z$. }
\end{figure}\FloatBarrier
In particular, if $\sigma(v)=m$, then 
\begin{align}\label{gapenergyformula}
H_{\emph{pos}}(\sigma^{v,r})-H_{\emph{pos}}(\sigma)=n_m(v)-n_r(v)+h\mathbbm{1}_{\{m=1\}}-h\mathbbm{1}_{\{r=1\}}.
\end{align}
\end{lemma}
\textit{Proof}. Let $\sigma\in\mathcal X$ and let $v\in V$. To find if a $v$-tile is stable for $\sigma$ we reduce ourselves to flip the spin on vertex $v$ from $\sigma(v)=m$ to a spin $r$ such that $n_r(v)>1$.
Indeed, otherwise the energy difference \eqref{energydifference3} is for sure strictly positive.
Let us divide the proof in several cases.\\
\textbf{Case 1.} Assume that $n_m(v)=0$ in $\sigma$. Then the corresponding $v$-tile is not stable for $\sigma$. Indeed, for any $m\in S$ and $r\notin\{1,m\}$, by flipping the spin on vertex $v$ from $m$ to $r$ we get
\begin{align}\label{firstalignpropmattonelle}
H_{\text{pos}}(\sigma^{v,r})-H_{\text{pos}}(\sigma)=-n_r(v)+h\mathbbm{1}_{\{m=1\}}.
\end{align}
Moreover, by flipping the spin on vertex $v$ from $m\neq 1$ to $1$ we have
\begin{align}
H_{\text{pos}}(\sigma^{v,1})-H_{\text{pos}}(\sigma)=-n_1(v)-h.
\end{align}
Hence, for any $m\in S$, if $v$ has spin $m$ and it has four nearest neighbors with spins different from $m$, then the tile centered in $v$ is not stable for $\sigma$ since the energy difference \eqref{energydifference3} is always strictly negative.\\
\textbf{Case 2.} Assume that $v\in V$ has three nearest neighbors with spin value different from $m$ in $\sigma$, i.e., $n_m(v)=1$. Then, in view of the energy difference \eqref{energydifference3}, for any $m\in S$ and $r\notin\{1,m\}$, by flipping the spin on vertex $v$ from $m$ to $r$ we have
\begin{align}
H_{\text{pos}}(\sigma^{v,r})-H_{\text{pos}}(\sigma)=1-n_r(v)+h\mathbbm{1}_{\{m=1\}}.
\end{align}
Furthermore, for any $m\neq 1$, if $r=1$ by flipping the spin on vertex $v$ from $m$ to $1$ we have
\begin{align}
H_{\text{pos}}(\sigma^{v,1})-H_{\text{pos}}(\sigma)=1-n_1(v)-h.
\end{align}
Hence, for any $m\in S$, $m\neq 1$, if $v$ has only one nearest neighbor with spin $m$, a tile centered in $v$ is stable  for $\sigma$ only if $v$ has nearest neighbors with spins different from each other and from $1$, see Figure \ref{figmattonellepos}(r). While, if $m=1$, if $v$ has only one nearest neighbor with spin $1$, a tile centered in $v$ is stable for $\sigma$ only if $v$ has nearest neighbors with spins different from each other, see Figure \ref{figmattonellepos}(s). \\
\textbf{Case 3.} Assume that $v\in V$ has two nearest neighbors with spin $m$ in $\sigma$, i.e., $n_m(v)=2$. Then, in view of the energy difference \eqref{energydifference3}, for any $m\in S$ and $r\notin\{1,m\}$, by flipping the spin on vertex $v$ from $m$ to $r$ we have
\begin{align}
H_{\text{pos}}(\sigma^{v,r})-H_{\text{pos}}(\sigma)=2-n_r(v)+h\mathbbm{1}_{\{m=1\}}.
\end{align}
Furthermore, by flipping the spin on vertex $v$ from $m\neq 1$ to $1$ we get
\begin{align}
H_{\text{pos}}(\sigma^{v,1})-H_{\text{pos}}(\sigma)=2-n_1(v)-h.
\end{align}
Hence, for any $m\in S$, if $v$ has two nearest neighbors with spin $m$ and two nearest neighbors with spin $1$, then the corresponding $v$-tile is not stable. In all the other cases, for any $m\in S$, if $v$ has two nearest neighbors with spin $m$, the corresponding $v$-tile is stable for $\sigma$, see Figure \ref{figmattonellepos}(f)--(q).\\
\textbf{Case 4.} Assume that $v\in V$ has three nearest neighbors with spin $m$ and one nearest neighbor $r$ in $\sigma$, i.e., $n_m(v)=3$ and $n_r(v)=1$. Then, for any $m\in S$ and $r\notin\{1,m\}$, by flipping the spin on vertex $v$ from $m$ to $r$ we have
\begin{align}
H_{\text{pos}}(\sigma^{v,r})-H_{\text{pos}}(\sigma)=2+h\mathbbm{1}_{\{m=1\}}.
\end{align}
Moreover, by flipping the spin on vertex $v$ from $m\neq 1$ to $1$ we get
\begin{align}
H_{\text{pos}}(\sigma^{v,1})-H_{\text{pos}}(\sigma)=2-h.
\end{align}
\textbf{Case 5.} Assume that $v\in V$ has four nearest neighbors with spin $m$, i.e., $n_m(v)=4$ in $\sigma$. Then, we have $n_r(v)=0$ and
\begin{align}
H_{\text{pos}}(\sigma^{v,r})-H_{\text{pos}}(\sigma)=4+h\mathbbm{1}_{\{m=1\}}.
\end{align}
Furthermore, by flipping the spin on vertex $v$ from $m\neq 1$ to $1$ we get
\begin{align}\label{lastalignpropmattonelle}
H_{\text{pos}}(\sigma^{v,1})-H_{\text{pos}}(\sigma)=4-h.
\end{align}
From Case 4 and Case 5, for any $m\in S$, we get that a $v$-tile is always stable for $\sigma$ if $v$ has at least three nearest neighbors with spin $m$, see Figure \ref{figmattonellepos}(a)--(e). In particular, in all the cases $1-5$ we verify that  \eqref{gapenergyformula} is satisfied by \eqref{firstalignpropmattonelle}--\eqref{lastalignpropmattonelle}.
$\qed$\\

We are now able to define the set of the local minima $\mathscr M_{\text{pos}}$ and the set of the stable plateaux $\bar{\mathscr M}_{\text{pos}}$ of the energy function \eqref{hamiltonianpos}. More precisely, the set of local minima $\mathscr M_\text{pos}$ is the set of stable points, i.e., $\mathscr M_\text{pos}:=\{\sigma\in\mathcal X:\ H_\text{pos}(\mathscr F(\partial\{\sigma\}))>H_\text{pos}(\sigma)\}$. While, a  plateau $D\subset\mathcal X$, namely a maximal connected set of equal energy states, is said to be \textit{stable} if $H_\text{pos}(\mathscr F(\partial D))>H_\text{pos}(D)$.  These will be given by the following sets:
\begin{align}
&\mathscr M^1_\text{pos}:=\{\bold 1,\bold2,\dots,\bold q\};\notag\\
&\mathscr M^2_\text{pos}:=\{\sigma\in\mathcal X:\ \text{$\sigma$ has strips of any spin $m\in S$ of thickness larger than or equal to two}\};\notag\\
&\mathscr M^3_{\text{pos}}:=\{\sigma\in\mathcal X:\ \text{$\sigma$ has $1$-rectangles that are not interacting with sides of length larger than}\notag\\ &\text{\hspace*{38pt} or equal to two, either in a sea of spins $m$ or in an $m$-strip for some $m\in S\backslash\{1\}$}\};\notag\\
&\mathscr M^4_\text{pos}:=\{\sigma\in\mathcal X:\ \text{for any $r\in S$, $\sigma$ is covered by $r$-rectangles. Every $r$-rectangle with $r\neq 1$}\notag\\ &\text{\hspace*{38pt} has minimum side of length larger than or equal to two, while every $1$-rectangle}\notag\\ &\text{\hspace*{38pt} has minimum side of length larger than or equal to one}\}\cup\{\sigma\in\mathcal X:\ \text{for any $m,r,$}\notag\\ &\text{\hspace*{38pt} $s\in S$ different from each other},\ \sigma\ \text{has an $r$-strip of thickness one and adjacent to}\notag\\ &\text{\hspace*{38pt} an $m$-strip and to an $s$-strip}\};\notag\\
&\bar{\mathscr M}^1_\text{pos}:=\{\sigma\in\mathcal X:\ \text{for any $m,r\in S\backslash\{1\}$, $\sigma$ has an $m$-rectangle with two consecutive sides}\notag\\ &\text{\hspace*{38pt} adjacent to two $r$-rectangles and the sides on the interfaces are of different}\notag\\ &\text{\hspace*{38pt} length}\}.\notag
\end{align}

\begin{proposition}[Sets of local minima and of stable plateaux]\label{proplocalminimapos}
If the external magnetic field is positive, then
\begin{align}
\mathscr M_{\emph{pos}}\cup\bar{\mathscr M}_{\emph{pos}}:=\mathscr M^1_{\emph{pos}}\cup\mathscr M^2_{\emph{pos}}\cup\mathscr M^3_{\emph{pos}}\cup\mathscr M^4_{\emph{pos}}\cup\bar{\mathscr M}^1_\emph{pos}.
\end{align}
\end{proposition}
\textit{Proof.}  A configuration $\sigma\in\mathcal X$ is said to be a local minimum, respectively a stable plateau, when for any $v\in V$ and $s\in S$ the energy difference \eqref{energydifference3} is either strictly positive, respectively null. When $\sigma$ has at least one unstable $v$-tile, for some $v\in V$, by flipping the spin on vertex $v$ it is possible to define a configuration with energy value strictly lower than $H_{\text{pos}}(\sigma)$. Thus, in this case $\sigma$ does not belong to $\mathscr M_{\text{pos}}\cup\bar{\mathscr M}_{\text{pos}}$.
Hence, below we give a geometric characterization of any $\sigma\in\mathscr M_{\text{pos}}\cup\bar{\mathscr M}_{\text{pos}}$ under the constraint that for any $v\in V$ the corresponding $v$-tile is stable for $\sigma$. In order to do this, we exploit Lemma \ref{stabletilespositive} and during the proof we often refer to Figure \ref{figmattonellepos}. In fact, a local minimum and a stable plateau are necessarily the union of one or more classes of stable tiles in Figure \ref{figmattonellepos}. Hence, we obtain all the local minima by enumerating all the possible ways in which the tiles in Figure \ref{figmattonellepos} may be combined. \\

\textbf{Step 1}. If $\sigma$ has only stable tiles as in Figure \ref{figmattonellepos}(a)--(b), then there are no interfaces, i.e., no disagreeing edges, and $\sigma\in\mathscr M^1_\text{pos}$. \\

\textbf{Step 2}. If $\sigma$ has only stable tiles as in Figure \ref{figmattonellepos}(a)--(e), then either there are not interfaces or there are only either horizontal or vertical interfaces of length $L$ and $K$, respectively. Thus, $\sigma\in\mathscr M^1_\text{pos}\cup\mathscr M^2_\text{pos}$.\\

\textbf{Step 3}. Assume that $\sigma$ has only stable tiles as in Figure \ref{figmattonellepos}(a)--(g), (n) and (o). For sake of semplicity we proceed by several steps.\\
\hspace*{12pt}\textbf{Step 3.1}. The $v$-tiles as in Figure \ref{figmattonellepos}(n) and (o) are such that any spin-update on vertex $v$ strictly increases the energy of at least $h$. Henceforth, for the case Figure \ref{figmattonellepos}(n) and for any $m\in S$, we get that $\sigma$ may have some $1$-rectangles both in a sea of spins $m$ and inside a strip of spins $m$. Note that the minimum side of these rectangles must be larger than or equal to two. Indeed, from Lemma \ref{stabletilespositive} we know that a $v$-tile where $v$ has spin $1$ and it has one nearest neighbor $1$ and three nearest neighbors with spin $m$ is not stable for $\sigma$. On the other hand, at this step a stable tile as in Figure \ref{figmattonellepos}(o) belongs to the configuration  $\sigma\in\mathscr M_{\text{pos}}\cup\bar{\mathscr M}_{\text{pos}}$ if and only if it belongs to a $1$-strip between two $r$-strips.\\
\hspace*{12pt}\textbf{Step 3.2}. Now we focus on the stable tile for $\sigma$ as in Figure \ref{figmattonellepos}(f) and (g). If $r=1$, the stable tile as in Figure \ref{figmattonellepos}(f) and (g) is not stable for $\sigma$ by item (1) of Lemma \ref{stabletilespositive}. It follows that $\sigma$ has not $m$-rectangles neither inside a sea of spins $1$ nor inside a strip of spins $1$. Hence, assume that $m,r\neq1$, $r\neq m$. In this case, we have to study the clusters as in Figure \ref{step33pos}.

\begin{figure}[h!]
\centering
\begin{tikzpicture}[scale=0.7,transform shape]
\fill[black!88!white] (0,1.5)rectangle(2,3);
\fill[black!20!white] (0,0)rectangle(2,1.5)(2,1.5)rectangle(4,3);
\draw (0,0)rectangle(2,3) (0,1.5)--(4,1.5)--(4,3)--(2,3);
\draw (1.88,1.6) node[white] {{$v$}};\draw (1.65,1.645) node[white] {{$\tilde v$}};
\draw[<->] (0,3.1)--(2,3.1); \draw (1,3.05) node[above] {\footnotesize{$\ell$}};
\draw (2,-0.1) node[below] {\large(a)};

\fill[black!88!white] (7.85,0)rectangle(8.15,3);
\fill[black!20!white] (6,0)rectangle(7.85,3)(8.15,0)rectangle(10,3);
\draw (6,0)rectangle(10,3)(7.85,0)--(7.85,3)(8.15,0)--(8.15,3);
\draw (8,1.6) node[white] {{$v$}};\draw (8,1.87) node[white] {{$\hat v$}};
\draw[<->] (7.85,3.1)--(8.15,3.1); \draw (8,3.05) node[above] {\footnotesize{$1$}};
\draw (8,-0.1) node[below] {\large(b)};
\end{tikzpicture}
\caption{\label{step33pos} Possible clusters of spins different from $1$ by considering the stable tile depicted in Figure \ref{figmattonellepos}(f) and (g) for $\sigma$ centered in the vertex $v$. We color black the $m$-rectangle and gray the $r$-rectangles.}
\end{figure}
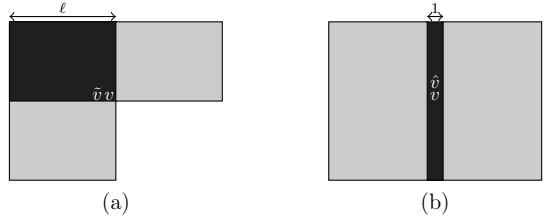\FloatBarrier

We claim that these types of clusters do not belong to $\sigma\in\mathscr M_{\text{pos}}\cup\bar{\mathscr M}_{\text{pos}}$ regardless of what is around them. Let us begin by considering the cluster depicted in Figure \ref{step33pos}(a). We construct a downhill path in which, starting from the vertex $v$, all the spins $m$ on those $\ell$ vertices next to a gray rectangle are flipped from $m$ to $r$, see Figure \ref{step33pos}(a).  In particular, during the first $\ell-1$ steps of this path, the spin which is flipping from $m$ to $r$ has two nearest neighbors with spin $m$ and two nearest neighbors with spin $r$. Indeed, after the spin update on vertex $v$ from $m$ to $r$, the spin $m$ on the left side of $v$, say $\tilde v$, has two nearest neighbors with spin $m$ and two nearest neighbors with spin $r$ and the path flips it to $r$ at zero energy cost. Hence, after this step the spin on the left side of $\tilde v$ has two nearest neighbors with spin $m$ and two nearest neighbors with spin $r$ and the path flips it from $m$ to $r$ without changing the energy. This construction is repeated until the last step, when the path flips the last spin $m$ to $r$. Indeed, in this case the spin $m$ that is flipping has two nearest neighbors with spin $r$, one nearest neighbor $m$ and the fourth nearest neighbor has spin different from $r$ and $m$. Thus, after this last flip the energy is reduced by $1$ and the claim is proved.\\
Let us now consider the cluster depicted in Figure \ref{step33pos}(b). In this case, we construct a downhill path given by two steps. First, it flips from $m$ to $r$ the spin $m$ on the vertex $v$ and the energy does not change. Then, it flips from $m$ to $r$ a spin $m$, say on vertex $\hat v$, which, after the previous flip, has three nearest neighbors with spin $r$ and one nearest neighbor $m$. Thus, by this flip the energy is reduced by two and the claim is verified.\\
In view of above, we conclude that any configuration $\sigma\in\mathscr M_{\text{pos}}\cup\bar{\mathscr M}_{\text{pos}}$ does not contain any $m$-rectangle in a sea of spin $r$ or with at least two consecutive sides adjacent to clusters of spin $r$ for any $m,r\in S$, see Figure \ref{step33pos}(a). Furthermore, for any $m,r\in S$, $r\neq m$, we get that in $\sigma$ there is not a spin $m$ that has two nearest neighbors with spin $m$ inside a rectangle with minimum side of length one and two nearest neighbors with spin $r$ belonging to two different $m$-rectangles, see Figure \ref{step33pos}(b). \\

\textbf{Step 4}. Let us now consider the case in which $\sigma$ have stable tiles as in Figure \ref{figmattonellepos}(a)--(q). First, we note that the stable $v$-tiles as in Figure \ref{figmattonellepos}(i), (m) and (q) are such that any spin-update on vertex $v$ strictly increases the energy by at least one. Thus, every $\sigma\in\mathscr M_{\text{pos}}\cup\bar{\mathscr M}_{\text{pos}}$ may have strips of thickness one as long as they are adjacent to two strips with different spins.
\begin{figure}[h!]
\centering
\begin{tikzpicture}[scale=0.7,transform shape]
\fill[black!12!white] (10.5,0)rectangle(12.2,0.7)(12.2,0.7)rectangle(13.3,2) (12.2,0.7)rectangle(11,1.6);
\fill[black!88!white] (12.2,0.7)rectangle(11,1.6);
\fill[black!50!white] (12.2,1.6) rectangle (11.3,2.1);
\draw (10.5,0)rectangle(12.2,0.7)(12.2,0.7)rectangle(13.3,2) (12.2,0.7)rectangle(11,1.6)(12.2,1.6) rectangle (11.3,2.1);
\draw[<->] (12.2,0.6)--(11,0.6); \draw(11.6,0.65) node[below] {\footnotesize{$\ell$}};
\draw  (12.08,1) node[white,below]{\footnotesize{$v$}};
\draw  (12.08,1.88) node[white,below]{\footnotesize{$w$}};
\draw(12,-0.1) node[below] {\Large(a)};

\fill[black!12!white] (14,0)rectangle(16,0.7)(16,0.7)rectangle(17,1.5)(14.5,0.7)rectangle(16,2);
\fill[black!88!white] (14.5,0.7)rectangle(16,2);
\fill[black!50!white] (14.5,0.7) rectangle (13.8,1.6);
\draw (14,0)rectangle(16,0.7)(16,0.7)rectangle(17,1.5)(14.5,0.7)rectangle(16,2) (14.5,0.7) rectangle (13.8,1.6);
\draw[<->] (14.5,0.6)--(16,0.6); \draw(15.25,0.65) node[below] {\footnotesize{$\ell$}};
\draw  (15.88,1) node[white,below]{\footnotesize{$v$}};
\draw  (14.38,1) node[white,below]{\footnotesize{$w$}};
\draw(15.5,-0.1) node[below] {\Large(b)};

\fill[black!12!white] (18.4,0)rectangle(19.5,0.7)(19.5,0.7)rectangle(20.5,1.5);
\fill[black!88!white] (18,0.7)rectangle(19.5,2);
\fill[black!50!white] (18.4,0.7) rectangle (17.4,0);
\draw (18.4,0)rectangle(19.5,0.7)(19.5,0.7)rectangle(20.5,1.5)(18,0.7)rectangle(19.5,2)(18.4,0.7) rectangle (17.4,0);
\draw  (19.38,1) node[white,below]{\footnotesize{$v$}};
\draw  (18.28,0.78) node[white,below]{\footnotesize{$w$}};
\draw(19,-0.1) node[below] {\Large(c)};
\end{tikzpicture}
\caption{\label{figvertexv} Illustration of an $m$-rectangle, that we color black, adjacent to two $r$-rectangles, that we color light gray. Furthermore, we color gray those $t$-rectangles with $t\in S\backslash\{r,m\}$.}
\end{figure}
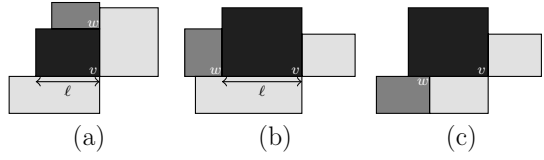\FloatBarrier
Now note that in view of the stable $v$-tiles depicted in Figure \ref{figmattonellepos}(h), (l) and (p), we have that for any $m,r\in S$, $m\neq r$, in $\sigma$ an $m$-rectangle may interact with an $r$-rectangle. More precisely, either an $m$-rectangle has a side adjacent to a side of an $r$-rectangle or there exists a vertex with spin $s\in S\backslash\{m,r\}$ that has two nearest neighbors with spin $s$, a nearest neighbor $m$ inside the $m$-rectangle and a nearest neighbor $r$ belonging to the $r$-rectangle.\\
Now note that in the previous step we studied the case in which an $m$-rectangle has two consecutive sides adjacent two $r$-rectangles, for any $m,r\in S\backslash\{1\}$, with the sides on the interfaces with the same length, see Figure \ref{step33pos}(a). By taking into account also the stable tiles depicted in Figure \ref{figmattonellepos}(h), (l) and (p), now $\sigma$ may have an $m$-rectangle which has two consecutive sides adjacent to two $r$-rectangles but the sides on the interfaces have not the same length, see Figure \ref{step33pos} where the black rectangle denotes the $m$-rectangle. 
It is easy to state that any configuration $\sigma\in\mathscr M_{\text{pos}}\cup\bar{\mathscr M}_{\text{pos}}$ has not the clusters depicted in Figure \ref{figvertexv}(a) and (b). Indeed, following the same strategy as in Step 3.3 above, we construct a path that flips from $m$ to $r$ all the spins $m$ adjacent to the side of length $\ell$ such that along this path the energy decreases. Hence, let us consider the case depicted in Figure \ref{figvertexv}(c). We prove that if $\sigma$ has this type of cluster surrounded by stable tiles, then $\sigma\in\bar{\mathscr M}_{\text{pos}}$. Indeed, in this case $\sigma$ does not communicate with a configuration with energy strictly lower but there exists a path which connects $\sigma$ and other configurations with the same energy. Indeed, we define a path which flips from $m$ to $r$ the spin $m$ adjacent to an $r$-rectangle. In particular, at any step the spin $m$ that is flipping has two nearest neighbors with spin $m$ and two nearest neighbors with spin $r$, see for instance the path depicted in Figure \ref{figzeronodiscesa} where the black rectangle denotes the $m$-rectangle. Since the sequence of configurations in Figure \ref{figzeronodiscesa} have the same energy value and they are connected with each other by means a path, they belong to a stable plateau. In particular, the energy along the path is constant. 

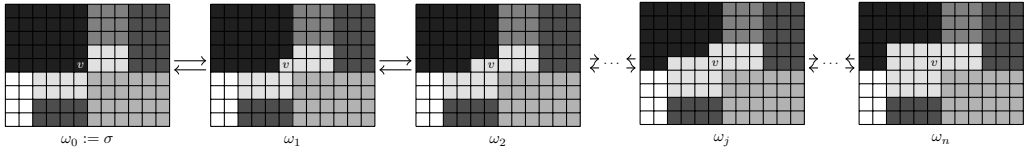
\begin{figure}[h!]
\centering
\begin{tikzpicture}[scale=0.6, transform shape]
\fill[black!12!white] (0,0) rectangle (3.6,2.7);
\fill[black!88!white] (0,2.7)rectangle(1.8,1.2);
\fill[black!0!white] (0,1.2)rectangle(0.6,0);
\fill[black!70!white] (0.6,0)rectangle(1.8,0.6)(3.6,2.7)rectangle(2.7,1.2);
\fill[black!30!white] (1.8,1.2)rectangle(3.6,0);
\fill[black!50!white] (1.8,2.7)rectangle(2.7,1.8);
\draw (1.65,1.35) node[white] {\footnotesize{$v$}};
\draw[step=0.3cm,color=black] (0,0) grid (3.6,2.7);
\draw(1.8,-0.1) node[below] {{$\omega_0:=\sigma$}};
\draw[->] (3.7,1.45)--(4.4,1.45); \draw [<-]  (3.7,1.25)--(4.4,1.25);
\fill[black!12!white] (4.5,0) rectangle (8.1,2.7);
\fill[black!88!white] (4.5,2.7)rectangle(6.3,1.2);
\fill[black!0!white] (4.5,1.2)rectangle(5.1,0);
\fill[black!70!white] (5.1,0)rectangle(6.3,0.6)(8.1,2.7)rectangle(7.2,1.2);
\fill[black!30!white] (6.3,1.2)rectangle(8.1,0);
\fill[black!50!white] (6.3,2.7)rectangle(7.2,1.8);
\fill[black!12!white] (6,1.2)rectangle(6.3,1.5);
\draw(6.3,-0.1) node[below] {{$\omega_1$}};
\draw (6.15,1.35) node {\footnotesize{$v$}};
\draw[step=0.3cm,color=black] (4.5,0) grid (8.1,2.7);
\draw[->] (8.2,1.45)--(8.9,1.45); \draw [<-]  (8.2,1.25)--(8.9,1.25);
\fill[black!12!white](9,0) rectangle (12.6,2.7);
\fill[black!88!white] (9,2.7)rectangle(10.8,1.2);
\fill[black!0!white] (9,1.2)rectangle(9.6,0);
\fill[black!70!white] (9.6,0)rectangle(10.8,0.6)(12.6,2.7)rectangle(11.7,1.2);
\fill[black!30!white] (10.8,1.2)rectangle(12.6,0);
\fill[black!50!white] (10.8,2.7)rectangle(11.7,1.8);
\fill[black!12!white] (10.8,1.2)rectangle(10.2,1.5);
\draw(10.8,-0.1) node[below] {{$\omega_2$}};
\draw (10.65,1.35) node {\footnotesize{$v$}};
\draw[step=0.3cm,color=black] (9,0) grid (12.6,2.7);
\end{tikzpicture}
\begin{tikzpicture}[scale=0.6, transform shape]
\draw[->] (12.7,1.45)--(12.9,1.45); \draw [<-]  (12.7,1.25)--(12.9,1.25);
\draw (13.2,1.35) node {\footnotesize$\dots$};
\draw[->] (13.7,1.45)--(13.5,1.45); \draw [<-]  (13.7,1.25)--(13.5,1.25);
\fill[black!12!white](13.8,0) rectangle (17.4,2.7);
\fill[black!88!white] (13.8,2.7)rectangle(15.6,1.2);
\fill[black!0!white] (13.8,1.2)rectangle(14.4,0);
\fill[black!70!white] (14.4,0)rectangle(15.6,0.6)(17.4,2.7)rectangle(16.5,1.2);
\fill[black!30!white] (15.6,1.2)rectangle(17.4,0);
\fill[black!50!white] (15.6,2.7)rectangle(16.5,1.8);
\fill[black!12!white] (15.6,1.2)rectangle(14.4,1.5) (15.6,1.5)rectangle(15.3,1.8);
\draw(15.6,-0.1) node[below] {{$\omega_j$}};
\draw (15.45,1.35) node {\footnotesize{$v$}};
\draw[step=0.3cm,color=black] (13.8,0) grid (17.4,2.7);
\draw[->] (17.5,1.45)--(17.7,1.45); \draw [<-]  (17.5,1.25)--(17.7,1.25);
\draw (18,1.35) node {\footnotesize$\dots$};
\draw[->] (18.5,1.45)--(18.3,1.45); \draw [<-]  (18.5,1.25)--(18.3,1.25);
\fill[black!12!white](18.6,0) rectangle (22.2,2.7);
\fill[black!88!white] (18.6,2.7)rectangle(20.4,1.2);
\fill[black!0!white] (18.6,1.2)rectangle(19.2,0);
\fill[black!70!white] (19.2,0)rectangle(20.4,0.6)(22.2,2.7)rectangle(21.3,1.2);
\fill[black!30!white] (20.4,1.2)rectangle(22.2,0);
\fill[black!50!white] (20.4,2.7)rectangle(21.3,1.8);
\fill[black!12!white] (20.4,1.2)rectangle(19.2,1.8);
\draw(20.4,-0.1) node[below] {{$\omega_n$}};
\draw (20.25,1.35) node {\footnotesize{$v$}};
\draw[step=0.3cm,color=black] (18.6,0) grid (22.2,2.7);
\end{tikzpicture}
\caption{\label{figzeronodiscesa} Example of a path $\omega:=(\omega_0,\dots,\omega_n)$ such that $H_\text{pos}(\omega_{i})=H_\text{pos}(\omega_j)$, for any $i,j=0,\dots,n$. Since all the configurations depicted have the same energy value and they are connected by means a path, they belong to a stable plateau.}
\end{figure}\FloatBarrier

\textbf{Step 5}. Let us now consider the case in which $\sigma$ may have stable tiles as in Figure \ref{figmattonellepos}(a)--(r). Let us focus on the stable tile depicted in Figure \ref{figmattonellepos}(r). Note that we have only to study the case in which the three nearest neighbors of $v$ with spins $r,s,t\in S\backslash\{m\}$ are different from each other. Indeed, all the other cases are not stable tiles for $\sigma$ by item (1) of Lemma \ref{stabletilespositive}. 
To aid the reader we refer to Figure \ref{figcolorvertx} for a pictorial illustration of the tile depicted in Figure \ref{figmattonellepos}(r), where we represent $r,m,s,t$ respectively by $\mycirclelightgray,\mycircleblack,\mycircleintgray,\mycirclegray$ and where we assume that $r,s,t\in S\backslash\{1,m\}$. Assume that this type of tile belongs to the configuration $\sigma$.
\begin{figure}[h!]
\centering
\begin{tikzpicture}[scale=1.1]
\fill[black!88!white] (-0.4,0.4)rectangle(0.4,0.8);
\fill[black!12!white] (0,0.8)rectangle(0.4,1.2);
\fill[black!50!white] (0,0)rectangle(0.4,0.4);
\fill[black!28!white] (0.4,0.4)rectangle(0.8,0.8);
\draw[step=0.4cm,color=black] (0,0) grid (0.4,1.2) (-0.4,0.4)grid(0.8,0.8);
\draw (0.2,0.6) node[white] {$v$} (-0.18,0.6) node[white] {\footnotesize{$v_1$}}(0.6,0.6) node {\footnotesize{$v_3$}} (0.22,1) node {\footnotesize{$v_2$}} (0.22,0.2) node {\footnotesize{$v_4$}} (-0.18,0.2) node{\footnotesize{$w_2$}} (-0.18,1) node{\footnotesize{$w_1$}};
\end{tikzpicture}
\caption{\label{figcolorvertx} Example of a $v$-tile equal to the one depicted in Figure \ref{figmattonellepos}(r). We do not color the vertices $w_1$ and $w_2$ since in the proof they assume different value in different steps.}
\end{figure}
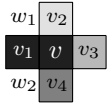\FloatBarrier
\hspace*{12pt}\textbf{Step 5.1}. If $n_m(v_1)=4$, then both $v_2$ and $v_4$ must have two nearest neighbors with their spin value and two nearest neighobors $m$, i.e., $n_m(v_i)=2$ and $n_{\sigma(v_i)}(v_i)=2$ for $i=2,4$. Otherwise, the $v_2$-tile and the $v_4$-tile would be not stable for $\sigma$ by item (1) of Lemma \ref{stabletilespositive} and we decrease the energy by a single spin update. It follows that the $v_2$-tile and $v_4$-tile coincide with those depicted in Figure \ref{figmattonellepos}(f), thus we are in the same situation of Step 3 and Step 4 and we use the same argument.\\
\hspace*{12pt}\textbf{Step 5.2}. If $n_m(v_1)=3$, similarly to the previous step, at least one between $v_2$ and $v_4$ have two nearest neighbors with its same spin value and two nearest neighbors of spin $m$. Hence, we are in an analogous situation as in Step 3 and in Step 4 and we conclude that such a tile belongs to a configuration which is either unstable for $\sigma$ or it belongs to a stable plateau.\\
\hspace*{12pt}\textbf{Step 5.3}. Let us now consider the case in which $v_1$ has only a nearest neighbor with spin $m$, i.e., the spin $m$ on vertex $v$.  We anticipate this case because it will be useful to study the case $n_m(v_1)=2$ in the next step. Let $r'$ be a spin value such that $n_{r'}(v_1)\ge1$. Along the path $\omega:=(\sigma,\sigma^{v,r},(\sigma^{v,r})^{v_1,r'})$ the energy decreases. Indeed, we have
\begin{align}
&H_{\text{pos}}(\sigma^{v,r})-H_{\text{pos}}(\sigma)=0,\notag \\
&H_{\text{pos}}((\sigma^{v,r})^{v_1,r'})-H_{\text{pos}}(\sigma^{v,r})=\begin{cases} -1-\mathbbm 1_{\{r=r'\}}-h\mathbbm 1_{\{r'=1\}}, &\text{if}\ n_{r'}(v_1)=1\ \text{in}\ \sigma;\\ -2-\mathbbm 1_{\{r=r'\}}-h\mathbbm 1_{\{r'=1\}},\ &\text{if}\  n_{r'}(v_1)=2\ \text{in}\ \sigma;\\ -3-\mathbbm 1_{\{r=r'\}}-h\mathbbm 1_{\{r'=1\}},\ &\text{if}\ n_{r'}(v_1)=3\ \text{in}\ \sigma.\end{cases}\notag
\end{align}
Note that the first equality follows by the fact that $\sigma^{v,r}$ is obtained by $\sigma$ by flipping the spin on vertex $v$ from $m$ to $r$ and $m, r\in S\backslash\{1\}$. On the other hand, the second equality follows by the fact that $(\sigma^{v,r})^{v_1,r'}$ is defined by $\sigma^{v,r}$ by flipping the spin on vertex $v_1$ from $m$ to $r'$ and $v_1$ has not any nearest neighbor with spin $m$ in $\sigma^{v,r}$. Thus, the energy decreases by the number of the nearest neighbor with spin $r'$ of $v_1$  in $\sigma^{v,r}$ and also by $h$ when $r'=1$, see \eqref{energydifference3}. It follows that in this case the tiles as depicted in Figure \ref{figmattonellepos}(r) do not belong to any configuration $\sigma\in\mathscr M_{\text{pos}}\cup\bar{\mathscr M}_{\text{pos}}$.\\
\hspace*{12pt}\textbf{Step 5.4}. The last case that we have to study is the one in which $v_1$ has two nearest neighbors with spin $m$. Obviously one of these is the vertex $v$. If the other spin $m$ lies either on the vertex $w_1$ or on the vertex $w_2$, then, it is in a configuration that falls in Step 5.1.  Hence, we conclude and in the sequel we assume that $v_1$ has the second nearest neighbor $m$ on the vertex which belongs to the same row where there are $v_1$ and $v$.  If $v_1$ has at least one nearest neighbor with a spin among $r,s,t$, say $r$, then along the path $\omega=(\sigma,\sigma^{v,r},(\sigma^{v,r})^{v_1,r})$ the energy decreases. Indeed, we have
\begin{align}
&H_{\text{pos}}(\sigma^{v,r})-H_{\text{pos}}(\sigma)=0,
&H_{\text{pos}}((\sigma^{v,r})^{v_1,r})-H_{\text{pos}}(\sigma^{v,r})\le-1.
\end{align}
In particular, the second inequality follows by \eqref{energydifference3}, by $r\neq 1$ and by the fact that $v_1$ has at least one nearest neighbor with spin $r$ in $\sigma$.
Thus in this case a tile as depicted in Figure \ref{figmattonellepos}(r) does not belong to any configuration $\sigma\in\mathscr M_{\text{pos}}\cup\bar{\mathscr M}_{\text{pos}}$. Hence assume that $v_1$ has two nearest neighbors with spin $m$ on vertices $v$ and $v_5$ and that for any $r'\notin\{r,t,z\}$, $n_{r'}(v_1)\in\{1,2\}$. If $n_{r'}(v_1)=2$, then we construct a path $\omega:=(\sigma,\sigma^{v,r},(\sigma^{v,r})^{v_1,r'})$ along which the energy decreases. Indeed, we have
\begin{align}
&H_{\text{pos}}(\sigma^{v,r})-H_{\text{pos}}(\sigma)=0,
&H_{\text{pos}}((\sigma^{v,r})^{v_1,r'})-H_{\text{pos}}(\sigma^{v,r})=-2-h\mathbbm 1_{\{r'=1\}},
\end{align}
where the first equality holds because $m$ and $r$ are different from $1$. On the other hand, the second equality follows by \eqref{energydifference3} the fact that $v_1$ has always two nearest neighbors with spin $r'$ in $\sigma^{v,r}$ since $r\neq r'$. Otherwise, $n_{r'}(v_1)=1$, see Figure \ref{step5finale}(a) where we represent $r'$ by $\mycircleerreprimo$. 

\begin{figure}[h!]
\centering
\begin{tikzpicture}[scale=1.1]
\fill[black!75!white] (-0.4,0.8)rectangle(0,1.2);
\fill[black!88!white] (-0.8,0.4)rectangle(0.4,0.8);
\fill[black!12!white] (0,0.8)rectangle(0.4,1.2);
\fill[black!50!white] (0,0)rectangle(0.4,0.4);
\fill[black!28!white] (0.4,0.4)rectangle(0.8,0.8);
\draw[step=0.4cm,color=black] (-0.4,0) grid (0.4,1.2) (-0.8,0.4)grid(0.8,0.8);
\draw (0.2,0.6) node[white] {$v$} (-0.18,0.6) node[white] {\footnotesize{$v_1$}}(0.6,0.6) node {\footnotesize{$v_3$}} (0.22,1) node {\footnotesize{$v_2$}} (0.22,0.2) node {\footnotesize{$v_4$}} (-0.58,0.6) node[white] {\footnotesize{$v_5$}};
\draw (0.1,-0.2) node[below]{(a)};
\end{tikzpicture}\ \ \ \ \ \ \ \ \ \
\begin{tikzpicture}[scale=1.1]
\fill[black!75!white] (-0.4,0.8)rectangle(0,1.2);
\fill[black!88!white] (-1.2,0.4)rectangle(0,0.8);
\fill[black!12!white] (0,0.8)rectangle(0.4,1.2)(0,0.4)rectangle(0.4,0.8);
\fill[black!50!white] (0,0)rectangle(0.4,0.4);
\fill[black!28!white] (0.4,0.4)rectangle(0.8,0.8);
\draw[step=0.4cm,color=black] (-0.4,0) grid (0.4,1.2) (-1.2,0.4)grid(0.8,0.8);
\draw (0.2,0.6) node {$v$} (-0.18,0.6) node[white] {\footnotesize{$v_1$}}(0.6,0.6) node {\footnotesize{$v_3$}} (0.22,1) node {\footnotesize{$v_2$}} (0.22,0.2) node {\footnotesize{$v_4$}} (-0.58,0.6) node[white] {\footnotesize{$v_5$}} (-1,0.6) node[white] {\footnotesize{$v_6$}};
\draw (0.1,-0.2) node[below]{(b)};
\end{tikzpicture}\ \ \ \ \ \ \ \ \ \
\begin{tikzpicture}
\fill[black!88!white] (-1.2,0.4)rectangle(-0.4,0.8);
\fill[black!12!white] (0.8,0.4)rectangle(1.2,1.2);
\fill[black!90!white] (0,0.4)rectangle(0.8,0.8);
\fill[black!35!white](0.4,0)rectangle(0.8,0.4)(-1.2,0)rectangle(-0.8,0.4);
\fill[black!70!white] (-1.6,0.4)rectangle(-1.2,0.8);
\fill[black!60!white] (0.4,0.8)rectangle(0.8,1.2);
\fill[black!28!white] (1.2,0.4)rectangle(1.6,0.8);
\fill[black!75!white] (0.4,0.8)rectangle(0.8,1.2);
\draw[step=0.4cm,color=black] (0.4,0) grid (1.2,1.2) (0,0.4)grid(1.6,0.8);
\draw[step=0.4cm,color=black](-1.2,0)grid(-0.8,1.2)(-1.6,0.4)rectangle(-0.4,0.8);
\draw (1,0.6) node {$v$} (0.6,0.6) node[white] {\footnotesize{$v_1$}}(1.4,0.6) node {\footnotesize{$v_3$}} (1.02,1) node {\footnotesize{$v_2$}} (1.02,0.2) node {\footnotesize{$v_4$}} (0.2,0.6) node[white] {\footnotesize{$v_5$}} (-1,0.6)node[white] {\footnotesize{$u$}}(-0.6,0.6) node[white] {\footnotesize{$v_n$}};
\draw (-0.17,1)  node {$\dots$} (-0.17,0.2) node {$\dots$};
\draw (-0.2,-0.2) node[below]{(c)};
\end{tikzpicture}
\caption{\label{step5finale} Illustration of the Step 5.4.}
\end{figure}
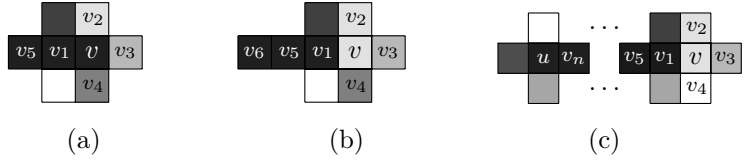\FloatBarrier
Without changing the energy, we flip from $m$ to a spin among $r,s,t$, say $r$, the spin on vertex $v$, see Figure \ref{step5finale}(b). We may repeat the procedure of above discussion by replacing $v$ with $v_1$ and $v_1$ with $v_5$. In any case, we conclude that a $v$-tile as in Figure \ref{figmattonellepos}(r) does not belong to any configuration $\sigma\in\mathscr M_{\text{pos}}\cup\bar{\mathscr M}_{\text{pos}}$. But at a certain point, we surely find a vertex $u$ such that the $u$-tile is analogous to the one centered in the vertex $v$ in Step 5.3, see Figure \ref{step5finale}(c). Note that by periodic boundary conditions this vertex may be the one on the right side of the vertex $v_3$, thus we conclude the proof of Case 5. Additionally, we note that in view of the above construction we have that the stable $v$-tiles as in Figure \ref{figmattonellepos}(i) belong only to a configuration $\sigma\in\mathscr M_{\text{pos}}\cup\bar{\mathscr M}_{\text{pos}}$ when they belong to a strip, i.e., when $v$ lies either on a row or on a column in which all the spins have the same value.  \\

 \textbf{Step 6}. Finally, we consider the case in which $\sigma$ has any possible stable tiles depicted in Figure \ref{figmattonellepos}(a)--(s). Note that the $v$-tile depicted in Figure \ref{figmattonellepos}(s) is such that any spin-update on vertex $v$ strictly increases the energy by at least $h$. 
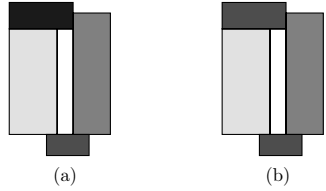
\begin{figure}[h!]
\centering
\begin{tikzpicture}[scale=0.7,transform shape]
\fill[black!12!white] (-0.9,0)rectangle(0,2);
\fill[black!50!white] (0.3,2.3)rectangle(1,0);
\fill[black!70!white] (-0.2,0)rectangle(0.6,-0.4);
\fill[black!90!white] (-0.9,2)rectangle(0.3,2.5);
\draw (0,0)rectangle(0.3,2)(-0.9,0)rectangle(0,2)(0.3,2.3)rectangle(1,0)(-0.2,0)rectangle(0.6,-0.4)(-0.9,2)rectangle(0.3,2.5);
\draw(0.15,-0.5) node[below] {(a)};

\fill[black!12!white] (3.1,0)rectangle(4,2);
\fill[black!50!white] (4.3,2.3)rectangle(5,0);
\fill[black!70!white] (3.8,0)rectangle(4.6,-0.4);
\fill[black!70!white] (3.1,2)rectangle(4.3,2.5);
\draw (4,0)rectangle(4.3,2)(3.1,0)rectangle(4,2)(4.3,2.3)rectangle(5,0)(3.8,0)rectangle(4.6,-0.4)(3.1,2)rectangle(4.3,2.5);
\draw(4.15,-0.5) node[below] {(b)};
\end{tikzpicture}
\caption{\label{figurecase6pos} Pictorial illustration of Step 6. The white rectangle denote a $1$-rectangle of minimum side of length $1$.}
\end{figure}\FloatBarrier
Hence, using these stable tiles together with those as in Figure \ref{figmattonellepos}(o) and (q), we get that a configuration $\sigma\in\mathscr M_{\text{pos}}\cup\bar{\mathscr M}_{\text{pos}}$ may have a cluster as depicted in Figure \ref{figurecase6pos}, i.e., $\sigma$ may have a $1$-rectangle with a side of length one and the other one of length larger than or equal to two. This rectangle staisfies the following conditions: there are no two consecutive sides adjacent to rectangles of the same spins and the sides of length different from one have to be adjacent to rectangles with spins different from each other. \\

Hence, we conclude that if $\sigma$ is characterized by the stable tiles depicted in Figure \ref{figmattonellepos}(a)--(s), then $\sigma\in\mathscr M^1_{\text{pos}}\cup\mathscr M^2_{\text{pos}}\cup\mathscr M^3_{\text{pos}}\cup\mathscr M^4_{\text{pos}}\cup\bar{\mathscr M}^1_\text{pos}$. $\qed$\\

We are now able to prove Proposition \ref{ricorrenzapos}.

\textit{Proof of Proposition \ref{ricorrenzapos}.}  Let $\widetilde{\mathscr M}_{\text{pos}}:=(\mathscr M_{\text{pos}}\backslash\{\bold 1,\dots,\bold q\})\cup\bar{\mathscr M_{\text{pos}}}.$ In order to conclude the proof, it is enough to focus on the configurations belonging to $\widetilde{\mathscr M}_{\text{pos}}$. Given $\eta\in\widetilde{\mathscr M}_{\text{pos}}$, we prove that $V_\eta^\text{pos}$ is smaller than or equal to $V^*:=2 <\Gamma_{\text{pos}}(\bold m,\mathcal X^s_{\text{pos}})$ for any $\bold m\in\mathcal X^m_\text{pos}$. \\

Let us first give an outline of the proof. First, we estimate the stability level for those configurations in $\mathscr M_\text{pos}^2\subset\widetilde{\mathscr M}_{\text{pos}}$ that are characterized by two or more strips of different spins, see Figure \ref{figurerecurrencepositive}(a)--(b). Second, we compute the stability level for those configurations in $\mathscr M_\text{pos}^3\subset\widetilde{\mathscr M}_{\text{pos}}$ that are characterized by a sea of $m$, for any $m\in S\backslash\{1\}$, with some rectangles of spins $1$ which do not interact among each other, see Figure \ref{figurerecurrencepositive}(c). Third, we consider the case of those configurations in $\mathscr M_\text{pos}^4\cup\bar{\mathscr M_\text{pos}}\subset\widetilde{\mathscr M}_{\text{pos}}$ that are covered by interacting rectangles, see Figure \ref{figurerecurrencepositive}(d).
Note that any configuration in $\widetilde{\mathscr M}_{\text{pos}}$ belongs to one of the three cases above. Indeed, $\widetilde{\mathscr M}_{\text{pos}}=\mathscr M^2_{\text{pos}}\cup\mathscr M^3_{\text{pos}}\cup\mathscr M^4_{\text{pos}}\cup\bar{\mathscr M}_\text{pos}$ and, any $\eta\in\widetilde{\mathscr M}_{\text{pos}}$ is such that it has at least two strips, or at least an isolated $1$-rectangle inside either a sea of spins $m$ or a strip of spins $m$, or a couple of interacting rectangles. See Figure \ref{figurerecurrencepositive}(e).
\begin{figure}[h!]
\centering
\begin{tikzpicture}[scale=0.7,transform shape]
\fill[black!40!white] (0.9,0)rectangle(1.5,2.7);
\fill[black!60!white] (2.4,0)rectangle(3.6,2.7);
\fill[black!15!white] (0,0)rectangle(0.9,2.7);
\draw[step=0.3cm,color=black] (0,0) grid (3.6,2.7);
\draw (1.8,-0.1) node[below] {\large(a)};

\fill[black!40!white] (4.5,0)rectangle(8.1,0.6)(4.5,2.4)rectangle(8.1,2.7);
\fill[black!70!white] (4.5,0.6)rectangle(8.1,1.5);
\draw[step=0.3cm,color=black] (4.5,0) grid (8.1,2.7);
\draw (6.3,-0.1) node[below] {\large(b)};

\fill[black!50!white] (9,0) rectangle (12.6,2.7);
\fill[white] (9.3,0.3)rectangle(10.2,1.8)(10.8,0.3)rectangle(11.7,1.2)(11.1,1.8)rectangle(12.3,2.4);
\draw[step=0.3cm,color=black] (9,0) grid (12.6,2.7);
\draw (10.8,-0.1) node[below] {(c)};
\end{tikzpicture}\\
\begin{tikzpicture}[scale=0.7, transform shape]
\fill[black!20!white] (0,0) rectangle (3.6,2.7);
\fill[black!30!white] (0.9,1.2)rectangle(1.8,1.8);
\fill[black!10!white] (0,1.8)rectangle(1.8,2.7);
\fill[black!50!white] (0.6,0)rectangle(1.8,1.2);
\fill[black!70!white] (1.8,1.5)rectangle(2.7,0.9)(0,0)rectangle(0.6,1.8);
\fill[white] (2.7,1.5)rectangle(3.6,0.9);
\draw[step=0.3cm,color=black] (0,0) grid (3.6,2.7);
\draw (1.8,-0.1) node[below] {\large(d)};

\fill[black!40!white] (4.5,0)rectangle(8.1,0.6)(4.5,2.4)rectangle(8.1,2.7);
\fill[black!70!white] (4.5,0.3)rectangle(8.1,1.8)(4.5,0)rectangle(8.1,0.6)(4.5,2.4)rectangle(8.1,2.7);
\fill[white] (5.1,0.3)rectangle(6.3,1.2)(6.9,0.6)rectangle(7.5,1.2);
\fill[black!20!white] (4.5,1.8)rectangle(5.7,2.4);
\fill[black!88!white] (8.1,1.8)rectangle(5.7,2.4);
\draw[step=0.3cm,color=black] (4.5,0) grid (8.1,2.7);
\draw (6.3,-0.1) node[below] {\large(e)};
\end{tikzpicture}
\caption{\label{figurerecurrencepositive} Example of configurations belonging to $\widetilde{\mathscr M}_{\text{pos}}$. We color white the vertices with spin $1$ and with the other colors the vertices with spins different from $1$.}
\end{figure}
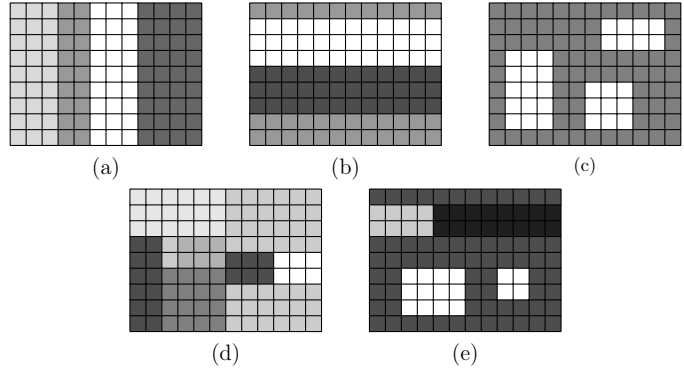\FloatBarrier

\textbf{Case 1.} Let us begin by assuming that $\eta$ only contains either horizontal or vertical strips. For concreteness, consider the case represented in Figure \ref{figurerecurrencepositive}(a). The case represented in Figure \ref{figurerecurrencepositive}(b) may be studied similarly. 
Assume that $\eta$ has an $m$-strip $a\times K$ adjacent to an $r$-strip $b\times K$ with $a,b\in\mathbb Z$, $a,b\ge 1$. Assume that $m,r\in S$, $m\neq 1$. Let $\bar\eta$ be the configuration obtained from $\eta$ by flipping  all the spins $m$ belonging to the $m$-strip from $m$ to $r$. We define a path $\omega:\eta\to\bar\eta$ as the concatenation of $a$ paths $\omega^{(1)},\dots,\omega^{(a)}$.  In particular, for any $i=1,\dots,a-1$, \[\omega^{(i)}:=(\omega^{(i)}_0=\eta_{i-1},\omega^{(i)}_1,\dots,\omega^{(i)}_{K}=\eta_{i}),\] where $\eta_{0}\equiv\eta$ and $\eta_i$ is the configuration in which the initial $m$-strip is reduced to a strip $(a-i)\times K$ and the $r$-strip to a strip $(b+i)\times K$. For any $i=1,\dots,a-1$, the path $\omega^{(i)}$ flips consecutively from $m$ to $r$ those spins $m$ belonging to the column next to the $r$-strip. In particular, note that for any $i=1,\dots,a-1$, $\omega^{(i)}_1:=(\omega^{(i)}_0)^{(v,r)}$ such that $\omega^{(i)}_0(v)=m$ and, in $\omega^{(i)}_0$, $v$ is the center of a tile as in Figure \ref{figmattonellepos}(c). Thus, using \eqref{energydifference3}, we have
\begin{align}
&H_{\text{pos}}(\omega^{(i)}_1)-H_{\text{pos}}(\omega^{(i)}_0)=2-h\mathbbm 1_{\{r=1\}}.\label{recurrenceposrecone} 
\end{align}
If $j=2,\dots,K-1$, then $\omega^{(i)}_j:=(\omega^{(i)}_{j-1})^{(v,r)}$ such that $\omega^{(i)}_{j-1}(v)=m$ and, in $\omega^{(i)}_{j-1}$, $v$ is the center of a tile as in Figure \ref{figmattonellepos}(f). Thus, using \eqref{energydifference3}, we have
\begin{align}
&H_{\text{pos}}(\omega^{(i)}_j)-H_{\text{pos}}(\omega^{(i)}_{j-1})=-h\mathbbm 1_{\{r=1\}}, 
\end{align}
for $j=2,\dots,K-1$. Finally, $\omega^{(i)}_K:=(\omega^{(i)}_{K-1})^{(v,r)}$ such that $\omega^{(i)}_{K-1}(v)=m$ and, in $\omega^{(i)}_{K-1}$, $v$ has three nearest neighbors $r$ and only one nearest neighbor $m$. Thus, using \eqref{energydifference3}, we have
\begin{align}
&H_{\text{pos}}(\omega^{(i)}_K)-H_{\text{pos}}(\omega^{(i)}_{K-1})=-2-h\mathbbm 1_{\{r=1\}}.
\end{align}  Then, for any $i=1,\dots,a-1$, 
the maximum energy value along $\omega^{(i)}$ is reached at the first step.  Finally, we construct a path $\omega^{(a)}:=(\omega^{(a)}_0=\eta_{a-1},\dots,\omega^{(a)}_{K}=\bar\eta)$ that flips consecutively from $m$ to $r$ the spins $m$ of the remaining column of the initial $m$-strip in $\eta$. Note that for the energy difference there are two possible values depending on whether the strips next to the initial $m$ strip $a\times K$ have the same value $r$ or one has value $r$ and the other $s\neq r,m$. 
Hence, if $v_i$ is the vertex whose spin is flipping in $\omega^{(a)}_{i-1}$ at the step $i$, if $i=1$, then $v_1$ has two nearest neighbors with spin $m$ and at most two nearest neighbors with spin $r$, i.e.,
\begin{align}
&H_{\text{pos}}(\omega^{(a)}_1)-H_{\text{pos}}(\eta_{a-1})=\begin{cases}1-h\mathbbm 1_{\{r=1\}},\ &\text{if $n_r(v_1)=1$,
}\\ -h\mathbbm 1_{\{r=1\}},\ &\text{if  $n_r(v_1)=2$.
 }  \end{cases}
 \end{align}
 Otherwise, if $i=2,\dots,K$, then $v_i$ has one nearest neighbor with spin $m$ and at most three nearest neighbors with spin $r$. Thus, using \eqref{energydifference3}, we have 
 \begin{align}
&H_{\text{pos}}(\omega^{(a)}_i)-H_{\text{pos}}(\omega^{(a)}_{i-1})=\begin{cases}-1-h\mathbbm 1_{\{r=1\}},\ &\text{if $n_r(v_i)=2$,
}\\ -2-h\mathbbm 1_{\{r=1\}},\ &\text{if $n_r(v_i)=3$,
}  \end{cases} \label{recurrenceposreconebis}
\end{align}
In view of \eqref{recurrenceposrecone}--\eqref{recurrenceposreconebis}, we get 
 $H(\eta)>H(\bar\eta)$. Furthermore, since the maximum energy value is reached at the first step, we get $V_\eta^\text{pos}\le2= V^*$.\\

\textbf{Case 2.} For any $m\in S$, $m\neq 1$, let us now consider the local minimum $\eta$ characterized by a sea of spins $m$ with some not interacting $1$-rectangles, see for instance Figure \ref{figurerecurrencepositive}(c). We distinguish two cases:
\begin{itemize}
\item[(i)] $\eta$ has at least a rectangle $R_{\ell_1\times\ell_2}$ of spins $1$ with its minimum side of length $\ell:=\min\{\ell_1,\ell_2\}$ larger than or equal to $\ell^*$,
\item[(ii)] $\eta$ has only rectangles $R_{\ell_1\times\ell_2}$ of spins $1$ with a side of length $\ell$ smaller than $\ell^*$.
\end{itemize}
In the first case, we define a path $\omega=(\omega_0,\dots,\omega_{\ell-1})$, where $\omega_0=\eta$ and $\omega_{\ell-1}=\tilde\eta$, that flips consecutively from $m$ to $1$ those spins $m$ adjacent to a side of length $\ell\ge\ell^*$. Given $v$ the vertex in which $\omega_{i-1}$ and $\omega_i$ differ, if $i=1$, $v$ has three nearest neighbors with spin $m$ and one nearest neighbor with spin $1$ in $\omega_0=\eta$. Otherwise, if $i=2,\dots,\ell-1$, $v$ has two nearest neighbors with spin $m$ and two nearest neighbors with spin $1$ in $\omega_{i-1}$. Thus, according to \eqref{energydifference3}, we obtain  
\begin{align}
&H_{\text{pos}}(\omega_1)-H_{\text{pos}}(\eta)=(3-1)-h=2-h,\label{align1}\\
&H_{\text{pos}}(\omega_i)-H_{\text{pos}}(\omega_{i-1})=(2-2)-h=-h, \ \text{for}\ i=2,\dots,\ell-1 \label{align2}.
\end{align}
Using \eqref{align1} and \eqref{align2}, it follows that $H_\text{pos}(\tilde\eta)-H_\text{pos}(\eta)=2-h\ell$. If $\ell>\ell^*$, then  using the definition of $\ell^*$ in \eqref{ellestar} we have $2-h\ell<2-h\ell^*<0$. Thus, $H_\text{pos}(\tilde\eta)<H_\text{pos}(\eta)$. Since the maximum energy is reached at the first step, from \eqref{align1} we have $V_\eta^\text{pos}=2-h< V^*$. 
Otherwise, $\eta$ has only rectangles $R_{\ell^*\times\ell^*}$ of spins $1$, thus $\ell=\ell^*$, the final configuration $\tilde\eta$ has a $1$-rectangle $\tilde R:=R_{\ell^*\times(\ell^*+1)}$. Either $\tilde R$ does not interact with the other rectangles of $\tilde\eta$, or $\tilde R$ interacts with another $1$-rectangle $\hat R$. In the former case, we
conclude by repeating the above construction along the side of length $\ell^*+1>\ell^*$. In the latter case, there exists in $\tilde\eta$ a vertex $w$ that interacts with both $\tilde R$ and $\hat R$. This vertex $w$ has two nearest neighbors with spin $m$ and two nearest neighbors with spin $1$ inside the $1$-rectangles $\tilde R$ and $\hat R$. Hence, set $\hat\eta:={\tilde\eta}^{(w,1)}$ and using \eqref{energydifference3}, we get
\begin{align}
H_{\text{pos}}(\hat\eta)-H_{\text{pos}}(\tilde\eta)=(2-2)-h=-h<0.\label{align3recurrencepos}
\end{align}
Hence, in view of \eqref{align1}--\eqref{align3recurrencepos}, along the path $(\eta,\omega_1,\dots,\omega_{\ell-2},\tilde\eta,\hat\eta)$  the maximum energy is reached at the first step. Thus, using \eqref{align1} we conclude that $V_\eta^\text{pos}=2-h<V^*$.\\

Now, let us focus on the case (ii). We define a path $\omega=(\omega_0,\dots,\omega_{\ell-1})$ that flips consecutively, from $1$ to $m$, those spins next to a side of length $\ell<\ell^*$. Given $v$ the vertex in which $\omega_{i-1}$ and $\omega_i$ differ, if $i=2,\dots,\ell-2$, in $\omega_{i-1}$ the vertex $v$ has two nearest neighbors with spin $m$ and  two nearest neighbors with spin $1$. Otherwise, if $i=\ell-1$, in $\omega_{\ell-2}$ the vertex $v$ has one neighbor with spin $1$ and three neighbors with spin $m$. Hence, according to \eqref{energydifference3}, we have
\begin{align}
&H_{\text{pos}}(\omega_i)-H_{\text{pos}}(\omega_{i-1})= (2-2)+h=h, \ \text{for}\ i=1,\dots,\ell-2;\\
&H_{\text{pos}}(\omega_{\ell-1})-H_{\text{pos}}(\omega_{\ell-2})=(1-3)+h=-(2-h).
\end{align}
Thus, the maximum energy along $\omega$ is reached at the step $\ell-1$ and, $V_\eta^\text{pos}=h(\ell-1)\le V^*$ since $h(\ell-1)<h(\ell^*-1)<2-h$ where the last inequality follows by \eqref{ellestar}.\\

\textbf{Case 3.} Finally, we focus on a local minimum $\eta$ covered by adjacent rectangles, see for instance Figure \ref{figurerecurrencepositive}(d). Let $m,r\in S\backslash\{1\}$, $m\neq r$. Indeed, any $m$-rectangle has to interact with at least an $r$-rectangle with $r\neq 1$ since otherwise it could be surrounded by spins $1$ that is in contradiction with $\eta\in\widetilde{\mathscr M}_\text{pos}$ by Proposition \ref{proplocalminimapos}.\\ Assume that $\eta$ has an $m$-rectangle $\bar R:=R_{a\times b}$ and an $r$-rectangle $\tilde R:=R_{c\times d}$ such that the $m$-rectangle $\bar R$ has a side of length $a$ adjacent to a side of the $r$-rectangle $\tilde R$ of length $c\ge a$. The case $c<a$ may be studied by interchanging the role of the spins $m$ and $r$. Given $\bar\eta$ the configuration obtained from $\eta$ by flipping from $m$ to $r$ all the spins $m$ belonging to $\bar R$, we construct a path $\omega:\eta\to\bar\eta$ as the concatenation of $b$ paths $\omega^{(1)},\dots,\omega^{(b)}$. In particular, for any $i=1,\dots,b-1$, \[\omega^{(i)}:=(\omega^{(i)}_0=\eta_{i-1},\omega^{(i)}_1,\dots,\omega^{(i)}_{a}=\eta_{i}),\] where $\eta_{0}\equiv\eta$ and $\eta_i$ is the configuration in which the initial $r$-rectangle $\tilde R$ is reduced to a rectangle $c\times d$ with a protuberance $a\times i$ and
 the initial $m$-rectangle $\bar R$ is reduced to a rectangle $a\times(b-i)$. For any $i=1,\dots,b-1$, the path $\omega^{(i)}$ flips consecutively from $m$ to $r$ those spins $m$ adjacent to the side of length $a$ of the $m$-rectangle $a\times(b-i)$. Given $v$ the vertex in which $\omega_{j-1}^{(i)}$ and $\omega_j^{(i)}$ differ, if $j=1$, then in $\omega^{(i)}_{0}$ the vertex  $v$ has two nearest neighbors with spin $m$, one nearest neighbor with spin $r$ and one nearest neighbor with spin different from $m,r$. Thus, using \eqref{energydifference3}, we get
\begin{align}
&H_{\text{pos}}(\omega^{(i)}_1)-H_{\text{pos}}(\omega^{(i)}_0)=1. \label{alignrecposlast}
\end{align} 
 If $j=2,\dots,a-1$, since in $\omega^{(i)}_{j-1}$ the vertex $v$ has two nearest neighbors $m$ and two nearest neighbors $r$, then
 \begin{align}
&H_{\text{pos}}(\omega^{(i)}_j)-H_{\text{pos}}(\omega^{(i)}_{j-1})=0.
\end{align}
 Finally, in $\omega^{(i)}_{a-1}$ the vertex  $v$ has one nearest neighbor with spin $m$, two nearest neighbors with spin $r$ and one nearest neighbor with spin different from $m,r$. Thus, according to \eqref{energydifference3} we get
\begin{align}
&H_{\text{pos}}(\omega^{(i)}_a)-H_{\text{pos}}(\omega^{(i)}_{a-1})=-1.
\end{align}
Then, for any $i=1,\dots,b-1$, 
 the maximum energy value along $\omega^{(i)}$ is reached at the first step. Finally, we define a path $\omega^{(b)}:=(\omega^{(b)}_0=\eta_{b-1},\dots,\omega^{(b)}_{\ell}=\bar\eta)$ that flips consecutively from $m$ to $r$ the spins $m$ belonging to the remaining $m$-rectangle $a\times1$. In particular, if $v_i$ is the vertex whose spin is flipping in $\omega_{i-1}^{(b)}$ at step $i$, if $i=1$, then $v_1$ has two nearest neighbors with spin $m$ and at most two nearest neighbors with spin $r$, i.e.,
\begin{align}
&H_{\text{pos}}(\omega^{(b)}_1)-H_{\text{pos}}(\eta_{b-1})=\begin{cases}0,\ &\text{if $n_r(v_1)=1$,
}\\ -1,\ &\text{if  $n_r(v_1)=2$.
}  \end{cases}
\end{align}
 Otherwise, if $i=2,\dots,a$, then $v_i$ has one nearest neighbor with spin $m$ and at most three nearest neighbors with spin $r$. Thus, using \eqref{energydifference3}, we have 
\begin{align}
&H_{\text{pos}}(\omega^{(b)}_i)-H_{\text{pos}}(\omega^{(b)}_{i-1})=\begin{cases}-1,\ &\text{if  $n_r(v_i)=2$,
}\\ -2,\ &\text{if  $n_r(v_i)=3$,
} \end{cases}\label{alignrecposlastbis}
\end{align}
Thanks to \eqref{alignrecposlast}--\eqref{alignrecposlastbis} 
and since by Proposition \ref{proplocalminimapos} we have $a\ge2$, we get $H(\eta)>H(\bar\eta)$. Moreover, since the maximum along $\omega$ is reached at the first step, by \eqref{alignrecposlast} we get $V_\eta^\text{pos}=1<V^*$.
 $\qed$
 
\subsection{Energy landscape and asymptotic behavior: proof of the main results}\label{proofenergylandscape}
We are now able to prove Corollary \ref{corollarygammatildepos} and Theorem \ref{teotimetargetG}.

\textit{Proof of Corollary \ref{corollarygammatildepos}.} By \cite[Lemma 3.6]{nardi2016hitting} we have that $\widetilde{\Gamma}_\text{pos}(\mathcal X\backslash\mathcal X^s_\text{pos})$ is the maximum energy that the process started in $\eta\in\mathcal X\backslash\mathcal X^s_\text{pos}$ has to overcome in order to arrive in $\mathcal X^s_\text{pos}=\{\bold 1\}$, i.e., 
\begin{align}
\widetilde{\Gamma}_\text{pos}(\mathcal X\backslash\mathcal X^s_\text{pos})=\max_{\eta\in\mathcal X\backslash\mathcal X^s_\text{pos}} \Gamma_\text{pos}(\eta,\mathcal X^s_\text{pos}). 
\end{align}
Hence, let us proceed to estimate $\Gamma_\text{pos}(\eta,\mathcal X^s_\text{pos})$ for any $\eta\in\mathcal X\backslash\mathcal X^s_\text{pos}$. Let $\bold m\in\mathcal X^m_\text{pos}$. Note that for any $\eta\in\mathcal X\backslash(\mathcal X^s_\text{pos}\cup\mathcal X^m_\text{pos})$ there are not initial cycles $\mathcal C^\eta_{\mathcal X^s_\text{pos}}(\Gamma_\text{pos}(\eta,\mathcal X^s_\text{pos}))$ deeper that $\mathcal C^\bold m_{\mathcal X^s_\text{pos}}(\Gamma_\text{pos}^m)$. While for any $\bold z\in\mathcal X^m_\text{pos}\backslash\{\bold m\}$, the initial cycles $\mathcal C^\bold z_{\mathcal X^s_\text{pos}}(\Gamma_\text{pos}^m)$ are as deep as $\mathcal C^\bold m_{\mathcal X^s_\text{pos}}(\Gamma_\text{pos}^m)$. By this fact, that holds since we are in the metastability scenario as in \cite[Subsection 3.5, Example 1]{nardi2016hitting}, we get that for any $\bold m\in\mathcal X^m_\text{pos}$
\begin{align}
\Gamma_\text{pos}(\eta,\mathcal X^s_\text{pos})=\Phi_\text{pos}(\eta,\mathcal X^s_\text{pos})-H_\text{pos}(\eta)\le\Phi_\text{pos}(\bold m,\mathcal X^s_\text{pos})-H_\text{pos}(\bold m)=\Gamma_\text{pos}(\bold m,\mathcal X^s_\text{pos}),
\end{align}
where the last equality follows by Theorem \ref{teometastablepos}. Thus, we conclude that \eqref{aligngammatildepos} is verified since, using \eqref{estimatestablevelmetapos}, we have $\Gamma_\text{pos}(\bold m,\mathcal X^s_\text{pos})=\Gamma^m_\text{pos}$. $\qed$\\

\textit{Proof of Theorem \ref{teotimetargetG}.} The theorem follows by \cite[Theorem 2.3]{fernandez2015asymptotically}.  
In order to apply this result it is enough to show that the pair $(\bold m,G)$ verifies the assumption 
\begin{align}\label{recurrencetimeanderror}
\sup_{\eta\in\mathcal X}\left(\tau^\eta_{\{\bold m,G\}}>R\right)\le\delta,
\end{align} 
with $R<\mathbb E(\tau^\bold m_G)$ and $\delta$ sufficiently small. To prove \eqref{recurrencetimeanderror}, let us distinguish two cases. \\
\textbf{Case 1.} Let $\eta\in\mathcal C^\bold m_{\mathcal X^s_\text{pos}}(\Gamma^m_\text{pos})$. By Equation (2.20) of \cite[Theorem 2.17]{manzo2004essential} applied to the cycle $\mathcal C^\bold m_{\mathcal X^s_\text{pos}}(\Gamma^m_\text{pos})$ we get that almost surely the process visits $\bold m$ before exiting from the cycle $\mathcal C^\bold m_{\mathcal X^s_\text{pos}}(\Gamma^m_\text{pos})$. More precisely, we have that there exists $k_1>0$ such that for any $\beta$ sufficiently large
\begin{align}\label{dismnosgm}
\mathbb P(\tau^\eta_\bold m>\tau^\eta_{\partial\mathcal C^\bold m_{\mathcal X^s_\text{pos}}(\Gamma^m_\text{pos})})\le e^{-k_1\beta}.
\end{align}
Since $\tau^\eta_G>\tau^\eta_{\partial\mathcal C^\bold m_{\mathcal X^s_\text{pos}}(\Gamma^m_\text{pos})}$, it follows that the process almost surely visits $\bold m$ before hitting $G$. Furthermore, since $\{\bold 1,\dots,\bold q\}\backslash\{\bold m\}\subset G$, we obtain that almost surely the process starting from $\eta$ visits $\bold m$ before hitting $\{\bold 1,\dots,\bold q\}\backslash\{\bold m\}$, i.e., $\tau^\eta_\bold m=\tau^\eta_{\{\bold 1,\dots,\bold q\}}$. Using the recurrence property given in Proposition \ref{teorecproppos}, we conclude that \eqref{recurrencetimeanderror} is satisfied by choosing $R=e^{\beta(2+\epsilon)}$ with $2+\epsilon<\Gamma^m_\text{pos}$ and $\delta=e^{-e^{k_2\beta}}$ with $k_2>0$. \\
\textbf{Case 2.} Let $\eta\in\mathcal X\backslash\mathcal C^\bold m_{\mathcal X^s_\text{pos}}(\Gamma^m_\text{pos})$. In this case \eqref{recurrencetimeanderror} is trivially verified for any $R$ and $\delta$ sufficiently small since $\eta$ belongs to the target. 
$\qed$

\section{Minimal gates and tube of typical trajectories}\label{secgatetube}
In this section we investigate on the minimal gates and the tube of typical paths for the transition from any $\mathbf m\in\mathcal X^m_\text{pos}$ to $\mathcal X^s_\text{pos}=\{\mathbf 1\}$. We further identify the union of all minimal gates also for the transition from a metastable state to the other metastable states.
\subsection{Identification of critical configurations for the transition from a metastable to the stable state}\label{mingatespos}
The goal of this subsection is to investigate the set of critical configurations for the transition from any $\bold m\in\mathcal X^m_\text{pos}$ to $\mathcal X^s_\text{pos}=\{\bold 1\}$. The idea of the proof of the following lemmas and proposition generalizes the proof of similiar results given in \cite[Section 6]{cirillo2013relaxation} for the Blume Capel model.\\
First we need to give some further definitions. For any $m\in S\backslash\{1\}$ we define $\mathscr D_{\text{pos}}^m\subset\mathcal X$ as the set of those configurations with $|\Lambda|-(\ell^*(\ell^*-1)+1)$ spins equal to $m$
\begin{align}\label{setDP}
\mathscr D_{\text{pos}}^m:=\{\sigma\in\mathcal X: N_m(\sigma)=|\Lambda|-(\ell^*(\ell^*-1)+1)\}.
\end{align}
 Furthermore, we define 
\begin{align}\label{setD+P}
\mathscr D_{\text{pos}}^{m,+}:=\{\sigma\in\mathcal X: N_m(\sigma)>|\Lambda|-(\ell^*(\ell^*-1)+1)\},
\end{align} note that $\bold m\in\mathscr D_{\text{pos}}^{m,+}$, and
\begin{align}\label{setD-P}
\mathscr D_{\text{pos}}^{m,-}:=\{\sigma\in\mathcal X: N_m(\sigma)<|\Lambda|-(\ell^*(\ell^*-1)+1)\}.
\end{align}
For any $\sigma\in\mathscr D_{\text{pos}}^m$, we remark that $\sigma$ has $\ell^*(\ell^*-1)+1$ spins different from $m$ and they may have not the same spin value and may belong to one or more clusters, see Figure \ref{figureDposexample}.
\begin{figure}[h!]
\centering
\begin{tikzpicture}[scale=0.7,transform shape]
\fill[black!85!white] (-0.6,0) rectangle (2.7,2.4);
\fill[black!50!white] (0,0)rectangle(0.9,0.6);
\fill[black!20!white] (0.3,1.2)rectangle(0.9,2.1)(0.3,0.6)rectangle(1.2,0.9)(2.1,1.8)rectangle(2.4,2.1);
\fill[white] (1.8,0.3)rectangle(2.4,0.9)(2.1,1.2)rectangle(2.4,1.5);
\draw[step=0.3cm,color=black] (-0.6,0) grid (2.7,2.4);
\draw (1.05,-0.1) node[below] {\large(a)};

\fill[black!85!white] (3.6,0) rectangle (6.9,2.4);
\fill[black!30!white](4.2,0.6)rectangle(5.4,1.8)(4.2,0.3)rectangle(4.5,0.6) (6,1.5)rectangle(6.3,2.1);
\fill[black!85!white] (4.8,1.2)rectangle(5.1,1.5);
\fill[black!30!white] (5.7,0.9)rectangle(6.6,1.2);
\draw[step=0.3cm,color=black] (3.6,0) grid (6.9,2.4);
\draw (5.25,-0.1) node[below] {\large(b)};

\fill[black!85!white] (7.8,0) rectangle (11.1,2.4);
\fill[black!20!white](8.4,0.6)rectangle(9,1.8)rectangle(9.9,0.9)(9.9,1.8)rectangle(10.5,1.2);
\draw[step=0.3cm,color=black] (7.8,0) grid (11.1,2.4);

\draw (9.45,-0.1) node[below] {\large(c)};
\end{tikzpicture}
\caption{\label{figureDposexample} Illustration of three examples of $\sigma\in{\mathscr D}_{\text{pos}}^m$ when $\ell^*=5$. We color black the vertices with spin $m$. In (a) the $\ell^*(\ell^*-1)+1=21$ spins different from $m$ have not all the same spin value and they belong to more clusters. In (b) these spins different from $m$ have the same spin value and they belong to three different clusters. In (c) the spins different from $m$ have the same spin value and they belong to a single cluster.}
\end{figure}
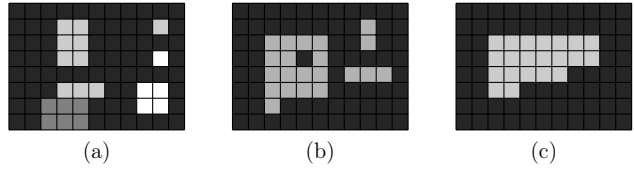\FloatBarrier
A \textit{two dimensional polyomino} on $\mathbb{Z}^2$ is a finite union of unit squares. The area of a polyomino is the number of its unit squares, while its perimeter is the cardinality of its boundary, namely, the number of unit edges of the dual lattice which intersect only one of the unit squares of the polyomino itself. Thus, the perimeter is the number of interfaces on $\mathbb Z^2$ between the sites inside the polyomino and those outside. The polyominoes with minimal perimeter among those with the same area are said to be \textit{minimal polyominoes}.
\begin{lemma}\label{bottomDP}
If the external magnetic field is positive, then for any $m\in\{2,\dots,q\}$ the minimum of the energy in $\mathscr D_{\emph{pos}}^m$ is achieved by those configurations in which the $\ell^*(\ell^*-1)+1$ spins different from $m$ are $1$ and they belong to a unique cluster of perimeter $4\ell^*$. 
More precisely,
\begin{align}\label{bottomDmpos}
\mathscr F(\mathscr D_{\emph{pos}}^m)=\{\sigma\in\mathscr D_{\emph{pos}}^m:\ & \sigma\ \text{has all spins $m$ except those in a unique cluster $C^1(\sigma)$ of spins}\notag \\ &\text{ $1$ of perimeter}\ 4\ell^*\}.
\end{align}
Moreover, 
\begin{align}\label{HbottomP}
H_{\text{\emph{pos}}}(\mathscr F(\mathscr D_{\emph{pos}}^m))=H_{\emph{pos}}(\bold m)+\Gamma_{\emph{pos}}(\bold m,\mathcal X^s_{\emph{pos}})=\Phi_{\emph{pos}}(\bold m,\mathcal X^s_{\emph{pos}}).
\end{align} 
\end{lemma}
\textit{Proof.} Let $m\in\{2,\dots,q\}$. From the definition of the Hamiltonian $H_\text{pos}$, \eqref{hamiltonianpos}, we get that the presence of disagreeing edges increases the energy, thus in order to identify the bottom of $\mathscr D_{\text{pos}}^m$ we have to consider those configurations $\sigma\in\mathscr D_{\text{pos}}^m$ in which the $\ell^*(\ell^*-1)+1$ spins different from $m$ belong to a single cluster. Moreover, given the number of the disagreeing edges, the presence of each spin $1$ decreases the energy by $h$ compared of the presence with other spins. Hence, the single cluster is full of spins $1$, say $C^1(\sigma)$, and it is inside a homogenous sea of spins $m$. 
Arguing like in the second part of the proof of Proposition \ref{lowerboundpos}, we have that $4\ell^*$ is the minimal perimeter of a polyomino of area $\ell^*(\ell^*-1)+1$. Thus, for any $\sigma\in\mathscr F(\mathscr D_{\text{pos}}^m)$,  $C^1(\sigma)$ must have perimeter $4\ell^*$. Hence, all the characteristics given in \eqref{bottomDmpos} are verified. Let us now prove \eqref{HbottomP}. By  \eqref{gatexm1} we get that $\mathcal W_{\text{pos}}(\bold m,\mathcal X^s_{\text{pos}})\subset\mathscr F(\mathscr D_{\text{pos}}^m)$, that is $H_\text{pos}(\mathscr F(\mathscr D_{\text{pos}}^m))=H_\text{pos}(\mathcal W_{\text{pos}}(\bold m,\mathcal X^s_{\text{pos}}))$. Thus, \eqref{HbottomP} is satisfied since for any $\eta\in\mathcal W_{\text{pos}}(\bold m,\mathcal X^s_{\text{pos}})$,
\begin{align}
H_{\text{pos}}(\eta)-H_{\text{pos}}(\bold m)=4\ell^*-h(\ell^*(\ell^*-1)+1)=\Gamma_{\text{pos}}(\bold m,\mathcal X^s_{\text{pos}}).
\end{align}
$\qed$\\

In the next corollary we show that every optimal path from $\bold m\in\mathcal X^m_{\text{pos}}$ to $\mathcal X^s_{\text{pos}}=\{\bold 1\}$ visits at least once $\mathscr F(\mathscr D_{\text{pos}}^m)$, i.e., we prove that $\mathscr F(\mathscr D_{\text{pos}}^m)$ is a gate for the transition from $\bold m$ to $\mathcal X^s_\text{pos}$.
\begin{cor}\label{Hmaxlemma3}
Let $\bold m\in\mathcal X^m_{\emph{pos}}$ and let $\omega\in\Omega_{\bold m,\mathcal X^s_\emph{pos}}^{opt}$. If the external magnetic field is positive, then $\omega\cap\mathscr F(\mathscr D_{\emph{pos}}^m)\neq\varnothing$. Hence, $\mathscr F(\mathscr D_{\emph{pos}}^m)$ is a gate for the transition $\bold m\to\mathcal X^s_\emph{pos}$.
\end{cor}
\textit{Proof.} Every path from $\bold m\in\mathcal X^m_{\text{pos}}$ to the stable configuration $\bold 1$ has to pass through the set $\mathcal V_k^m:=\{\sigma\in\mathcal X:\ N_m(\sigma)=k\}$ for any $k=|V|,\dots,0$. In particular, given $k^*:=\ell^*(\ell^*-1)+1$, any $\omega=(\omega_0,\dots,\omega_n) \in\Omega_{\bold m,\mathcal X^s_\text{pos}}^{opt}$ visits at least once the set $\mathcal V^m_{|\Lambda|-k^*}\equiv\mathscr D_\text{pos}^m$.  Hence, there exists $i\in\{0,\dots n\}$ such that $\omega_i\in\mathscr D_{\text{pos}}^m$. Since from \eqref{HbottomP} we have that the energy value of any configuration belonging to $\mathscr F(\mathscr D_{\text{pos}}^m)$ is equal to the min-max reached by any optimal path from $\bold m$ to $\mathcal X^s_{\text{pos}}$, we conclude that $\omega_i\in\mathscr F(\mathscr D_{\text{pos}}^m)$.  $\qed$\\ 

In the last result of this subsection, we prove that, for any $\bold m\in\mathcal X^m_{\text{pos}}$, every optimal path $\omega\in\Omega_{\bold m,\mathcal X^s_{\text{pos}}}^{opt}$ is such that $\omega\cap\mathcal W_{\text{pos}}(\bold m,\mathcal X^s_{\text{pos}})\neq\varnothing$. Hence, we show that  $\mathcal W_{\text{pos}}(\bold m,\mathcal X^s_{\text{pos}})$ is a gate for the transition $\bold m\to\mathcal X^s_{\text{pos}}$.

\begin{proposition}\label{propgatepos}
If the external magnetic field is positive, then for any $\bold m\in\mathcal X^m_{\emph{pos}}$ each path $\omega\in\Omega_{\bold m,\mathcal X^s_{\emph{pos}}}^{opt}$ visits $\mathcal W_{\emph{pos}}(\bold m,\mathcal X^s_{\emph{pos}})$. Hence, $\mathcal W_{\emph{pos}}(\bold m,\mathcal X^s_{\emph{pos}})$ is a gate for the transition $\bold m\to\mathcal X^s_\emph{pos}$.
\end{proposition} 
\textit{Proof.} For any $m\in S$, $m\neq 1$, let $\tilde{\mathscr D}^m_{\text{pos}}$ and $\hat{\mathscr D}^m_{\text{pos}}$ be the subsets of $\mathscr F(\mathscr D_{\text{pos}}^m)$ defined as follows. $\tilde{\mathscr D}^m_{\text{pos}}$ is  the set of those configurations of $\mathscr F(\mathscr D_{\text{pos}}^m)$ in which the boundary of the polyomino $C^1(\sigma)$ intersects each side of the boundary of its smallest surrounding rectangle $R(C^1(\sigma))$ on a set of the dual lattice $\mathbb{Z}^2+(1/2,1/2)$ made by at least two consecutive unit segments, see Figure \ref{figureexamplebis}(a). On the other hand, $\hat{\mathscr D}^m_{\text{pos}}$ is the set of those configurations of $\mathscr F(\mathscr D_{\text{pos}}^m)$ in which the boundary of the polyomino $C^1(\sigma)$ intersects at least one side of the boundary of $R(C^1(\sigma))$ in a single unit segment, see Figure \ref{figureexamplebis}(b) and (c).
\begin{figure}[h!]
\centering
\begin{tikzpicture}[scale=0.7,transform shape]
\fill[black!85!white] (0,0)rectangle(3.3,2.4);
\fill[white] (0.6,0.6)rectangle(1.2,1.2)(0.9,0.3)rectangle(1.5,1.8)(1.5,0.3)rectangle(1.8,2.1)(1.8,1.2)rectangle(2.1,2.1);
\draw[step=0.3cm,color=black] (0,0) grid (3.3,2.4);
\draw[thick,dotted,white] (0.6,0.3)rectangle(2.1,2.1);
\draw (1.65,-0.1) node[below] {\large(a)};

\fill[black!85!white] (3.6,0) rectangle (6.9,2.4);
\fill[white](4.5,0.3)rectangle(5.7,1.8)(5.1,1.8)rectangle(5.4,2.1);
\draw[step=0.3cm,color=black] (3.6,0) grid (6.9,2.4);
\draw[thick,dotted,white] (4.5,0.3)rectangle(5.7,2.1);
\draw (5.25,-0.1) node[below] {\large(b)};

\fill[black!85!white] (7.2,0) rectangle (10.5,2.4);
\fill[white](8.4,0.6)rectangle(9.6,2.1)(9.6,1.2)rectangle(9.9,1.5);
\draw[step=0.3cm,color=black] (7.2,0) grid (10.5,2.4);
\draw[thick,dotted,white] (8.4,0.6)rectangle(9.9,2.1);
\draw (8.85,-0.1) node[below] {\large(c)};

\fill[black!85!white] (10.8,0) rectangle (14.1,2.4);
\fill[white](12,0.6)rectangle(12.9,2.1)(11.7,0.6)rectangle(12,1.5)(12.9,1.2)rectangle(13.2,1.5)(12,0.3)rectangle(12.6,0.6);
\draw[step=0.3cm,color=black] (10.8,0) grid (14.1,2.4);
\draw[thick,dotted,white] (11.7,0.3)rectangle(13.2,2.1);
\draw (12.45,-0.1) node[below] {\large(d)};
\end{tikzpicture}
\caption{\label{figureexamplebis} Examples of $\sigma\in\tilde{\mathscr D}^m_{\text{pos}}$ (a) and of $\sigma\in\hat{\mathscr D}^m_{\text{pos}}$ (b) and (c)  when $\ell^*=5$. We associate the color black to the spin $m$, the color white to the spin $1$. The dotted rectangle represents $R(C^1(\sigma))$. Figure (d) is an example of configuration that does not belong to $\hat{\mathscr D}^m_{\text{pos}}$.}
\end{figure}
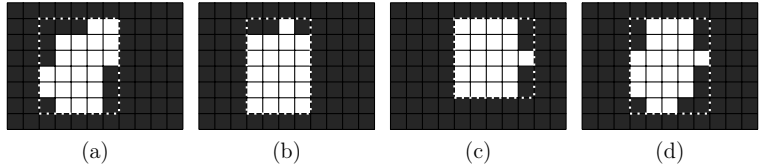\FloatBarrier
In particular note that $\mathscr F(\mathscr D_{\text{pos}}^m)=\tilde{\mathscr D}^m_{\text{pos}}\cup\hat{\mathscr D}^m_{\text{pos}}$. The proof proceeds in five steps.\\

\textbf{Step 1}. Our first aim is to prove that 
\begin{align}\label{Dhatmunion}
\hat{\mathscr D}^m_{\text{pos}}=\mathcal W_{\text{pos}}(\bold m,\mathcal X^s_{\text{pos}})\cup\mathcal W_{\text{pos}}'(\bold m,\mathcal X^s_{\text{pos}}).
\end{align}  
From \eqref{gatexm1} we have  $\mathcal W_{\text{pos}}(\bold m,\mathcal X^s_{\text{pos}})\cup\mathcal W_{\text{pos}}'(\bold m,\mathcal X^s_{\text{pos}})\subseteq\hat{\mathscr D}^m_{\text{pos}}$. Thus we reduce our proof to show that $\sigma\in\hat{\mathscr D}^m_{\text{pos}}$ implies $\sigma\in\mathcal W_{\text{pos}}(\bold m,\mathcal X^s_{\text{pos}})\cup\mathcal W_{\text{pos}}'(\bold m,\mathcal X^s_{\text{pos}})$. Note that this implication is not straightforward, since given $\sigma\in\hat{\mathscr D}^m_{\text{pos}}$, the boundary of the polyomino $C^1(\sigma)$ could intersect the other three sides of the boundary of its smallest surrounding rectangle $R(C^1(\sigma))$ in a proper subsets of the sides itself, see Figure \ref{figureexamplebis}(d) for an illustration of this hypothetical case. Hence, consider $\sigma\in\hat{\mathscr D}^m_{\text{pos}}$ and let $R(C^1(\sigma))=R_{(\ell^*+a)\times(\ell^*+b)}$ with $a,b \in\mathbb{Z}$. 
In view of the proof of Lemma \ref{bottomDP} we have that $C^1(\sigma)$ is a minimal polyomino and by \cite[Lemma 6.16]{cirillo2013relaxation} it is also convex and monotone, i.e., its perimeter of value $4\ell^*$ is equal to the one of $R(C^1(\sigma))$. Hence, the following equality holds
\begin{align}\label{equalityst}
4\ell^*=4\ell^*+2(a+b).
\end{align}
In particular, \eqref{equalityst} is satisfied only by $a=-b$. Now, let $\tilde R$ be the smallest rectangle surrounding the polyomino, say $\tilde C^1(\sigma)$, obtained by removing the unit protuberance from $C^1(\sigma)$. If $C^1(\sigma)$ has the unit protuberance adjacent to a side of length $\ell^*+a$, then $\tilde R$ is a rectangle $(\ell^*+a)\times(\ell^*-a-1)$. Note that $\tilde R$ must have an area larger than or equal to the number of spins $1$ of the polyomino $\tilde C^1(\sigma)$, that is $\ell^*(\ell^*-1)$. Thus, we have
\begin{align}\label{areartildepos1}
\text{Area}(\tilde R)=(\ell^*+a)(\ell^*-a-1)=\ell^*(\ell^*-1)-a^2-a\ge \ell^*(\ell^*-1) \iff -a^2-a\ge 0.
\end{align}
Since $a\in\mathbb{Z}$, $-a^2-a\ge 0$ is satisfied only if either $a=0$ or $a=-1$. 
Otherwise, if $C^1(\sigma)$ has the unit protuberance adjacent to a side of length $\ell^*-a$, then $\tilde R$ is a rectangle $(\ell^*+a-1)\times(\ell^*-a)$ and
\begin{align}\label{areartildepos2}
\text{Area}(\tilde R)=(\ell^*+a-1)(\ell^*-a)=\ell^*(\ell^*-1)-a^2+a\ge \ell^*(\ell^*-1) \iff -a^2+a\ge 0.
\end{align}
Since $a\in\mathbb{Z}$, $-a^2+a\ge 0$ is satisfied only if either $a=0$ or $a=1$. In both cases we get that $\tilde R$ is a rectangle of side lengths $\ell^*$ and $\ell^*-1$. Thus, if the protuberance is attached to one of the longest sides of $\tilde R$, then $\sigma\in\mathcal W_{\text{pos}}(\bold m,\mathcal X^s_{\text{pos}})$, otherwise $\sigma\in\mathcal W'_{\text{pos}}(\bold m,\mathcal X^s_{\text{pos}})$. In any case we conclude that \eqref{Dhatmunion} is satisfied.\\

\textbf{Step 2}. For any $\bold m\in\mathcal X^m_{\text{pos}}$ and for any path $\omega=(\omega_0,\dots,\omega_n)\in\Omega_{\bold m,\mathcal X^s_{\text{pos}}}^{opt}$, let 
\begin{align}\label{defgmomega}
g_m(\omega):=\{i\in\mathbb{N}: \omega_i\in\mathscr F(\mathscr D_{\text{pos}}^m),\  N_1(\omega_{i-1})=\ell^*(\ell^*-1),\ N_m(\omega_{i-1})=|\Lambda|-\ell^*(\ell^*-1)\}.
\end{align}
We claim that $g_m(\omega)\neq\varnothing$. Let $\omega=(\omega_0,\dots,\omega_n)\in\Omega_{\bold m,\mathcal X^s_{\text{pos}}}^{opt}$ and let
$j^*\le n$ be the smallest integer such that after $j^*$ the path leaves $\mathscr D^{m,+}_\text{pos}$, i.e., $(\omega_{j^*},\dots,\omega_n)\cap\mathscr D^{m,+}_\text{pos}=\varnothing$.
 Since $\omega_{j^*-1}$ is the last configuration in $\mathscr D^{m,+}_\text{pos}$, it follows  that $\omega_{j^*}\in\mathscr D^m_\text{pos}$ and, by the proof of Corollary \ref{Hmaxlemma3}, we have that $\omega_{j^*}\in\mathscr F(\mathscr D^m_\text{pos})$. Moreover, since $\omega_{j^*-1}$ is the last configuration in $\mathscr D^{m,+}_\text{pos}$, we have that $N_m(\omega_{j^*-1})=|\Lambda|-\ell^*(\ell^*-1)$ and $\omega_{j^*}$ is obtained by $\omega_{j^*-1}$ by flipping a spin $m$ from $m$ to $s\neq m$. Note that $N_m(\omega_{j^*-1})=|\Lambda|-\ell^*(\ell^*-1)$ implies $N_s(\omega_{j^*-1})\le\ell^*(\ell^*-1)$ for any $s\in S\backslash\{m\}$. By Lemma \ref{bottomDP}, $\omega_{j^*}\in\mathscr F(\mathscr D^m_\text{pos})$ implies  $N_1(\omega_{j^*})=\ell^*(\ell^*-1)+1$, thus $N_1(\omega_{j^*-1})<\ell^*(\ell^*-1)$ is not feasible since $\omega_{j^*}$ and $\omega_{j^*-1}$ differ by a single spin update which increases the number of spins $1$ of at most one. Then, $j^*\in g_m(\omega)$ and the claim is proved.\\

\textbf{Step 3}.  We claim that for any path $\omega\in\Omega_{\bold m,\mathcal X^s_{\text{pos}}}^{opt}$ one has $\omega_i\in\hat{\mathscr D}^m_{\text{pos}}$  for any $i\in g_m(\omega)$.
We argue by contradiction. Assume that there exists $i\in g_m(\omega)$ such that  $\omega_i\notin\hat{\mathscr D}^m_{\text{pos}}$ and $\omega_i\in\tilde{\mathscr D}^m_{\text{pos}}$.  Since $\omega_{i-1}$ is obtained from $\omega_i$ by flipping a spin $1$ to $m$ and since any configuration belonging to $\tilde{\mathscr D}^m_{\text{pos}}$ has all the spins $1$ with at least two nearest neighbors with spin $1$, using \eqref{energydifference3} we have
\begin{align}\label{disequalityHi}
H_{\text{pos}}(\omega_{i-1})-H_{\text{pos}}(\omega_i)\ge(2-2)+h=h>0.
\end{align} 
In particular, from \eqref{disequalityHi} we get a contradiction. Indeed, 
\begin{align}\label{absurdgate31}
\Phi_\omega^\text{pos}\ge H_{\text{pos}}(\omega_{i-1})>H_{\text{pos}}(\omega_i)=H_{\text{pos}}(\bold m)+\Gamma_{\text{pos}}(\bold m,\mathcal X^s_{\text{pos}})=\Phi_{\text{pos}}(\bold m,\mathcal X^s_{\text{pos}}),
\end{align}
where the first equality follows by \eqref{HbottomP}. Thus by \eqref{absurdgate31} $\omega$ is not an optimal path, which is a contradiction, the claim is proved and we conclude the proof of Step 3.\\

\textbf{Step 4}. Now we claim that for any $\bold m\in\mathcal X^m_{\text{pos}}$ and for any path $\omega\in\Omega_{\bold m,\mathcal X^s_{\text{pos}}}^{opt}$,
\begin{align}
\omega_i\in\mathscr F(\mathscr D_{\text{pos}}^m)\implies \omega_{i-1}, \omega_{i+1}\notin\mathscr D_{\text{pos}}^m.
\end{align}
Using Corollary \ref{Hmaxlemma3}, for any $\bold m\in\mathcal X^m_\text{pos}$ and any path $\omega\in\Omega_{\bold m,\mathcal X^s_{\text{pos}}}^{opt}$ there exists an integer $i$ such that $\omega_i\in\mathscr F(\mathscr D_{\text{pos}}^m)$.  Assume by contradiction that $\omega_{i+1}\in\mathscr D_{\text{pos}}^m$. In particular,  since $\omega_i$ and $\omega_{i+1}$ have the same number of spins $m$, note that $\omega_{i+1}$ is obtained by flipping a spin $1$ from $1$ to $t\neq 1$. Since $\omega_i(v)\neq t$ for every $v\in V$,  the above flip increases the energy, i.e., $H_{\text{pos}}(\omega_{i+1})>H_{\text{pos}}(\omega_i)$. Hence, using this inequality and \eqref{HbottomP}, we have 
\begin{align}
\Phi_\omega^\text{pos}\ge H_{\text{pos}}(\omega_{i+1})>H_{\text{pos}}(\omega_i)=H_{\text{pos}}(\bold m)+\Gamma_{\text{pos}}(\bold m,\mathcal X^s_{\text{pos}})=\Phi_{\text{pos}}(\bold m,\mathcal X^s_{\text{pos}}),
\end{align}
which implies the contradiction because $\omega$ is not optimal. Thus $\omega_{i+1}\notin\mathscr D_{\text{pos}}^m$ and similarly we show that also $\omega_{i-1}\notin\mathscr D_{\text{pos}}^m$.\\

\textbf{Step 5}. In this last step of the proof we claim that for any $\bold m\in\mathcal X^m_{\text{pos}}$ and for any path $\omega\in\Omega_{\bold m,\mathcal X^s_{\text{pos}}}^{opt}$ there exists a positive integer $i$ such that $\omega_i\in\mathcal W_{\text{pos}}(\bold m,\mathcal X^s_{\text{pos}})$. Arguing by contradiction, assume that there exists $\omega\in\Omega_{\bold m,\mathcal X^s_{\text{pos}}}^{opt}$ such that $\omega\cap\mathcal W_{\text{pos}}(\bold m,\mathcal X^s_{\text{pos}})=\varnothing$. Thanks to Corollary \ref{Hmaxlemma3}, we know that $\omega$ visits $\mathscr F(\mathscr D^m_{\text{pos}})$ and thanks to Step 4 we have that the configurations along $\omega$ belonging to $\mathscr F(\mathscr D^m_{\text{pos}})$ are not consecutive. More precisely, they are linked by a sub-path that belongs either to $\mathscr D_{\text{pos}}^{m,+}$ or $\mathscr D_{\text{pos}}^{m,-}$. If $n$ is the length of $\omega$, then let $j\le n$ be the smallest integer such that $\omega_j\in\mathscr F(\mathscr D^m_{\text{pos}})$ and such that $(\omega_j,\dots,\omega_n)\cap\mathscr D_{\text{pos}}^{m,+}=\varnothing$, thus, $j\in g_m(\omega)$ since $j$ plays the same role of $j^*$ in the proof of Step 2. Using \eqref{Dhatmunion}, Step 3 and the assumption $\omega\cap\mathcal W_{\text{pos}}(\bold m,\mathcal X^s_{\text{pos}})=\varnothing$, it follows that $\omega_j\in\mathcal W_{\text{pos}}'(\bold m,\mathcal X^s_{\text{pos}})$. Moreover, starting from $\omega_j\in\mathscr F(\mathscr D^m_{\text{pos}})$ the energy along the path decreases only by either
\begin{itemize}
\item[(i)] flipping the spin in the unit protuberance from $1$ to $m$, or
\item[(ii)] flipping a spin, with two nearest neighbors with spin $1$, from $m$ to $1$.
\end{itemize}
Since by the definition of $j$ we have that $\omega_{j-1}$ is the last that visits $\mathscr D_{\text{pos}}^{m,+}$, $\omega_{j+1}\notin\mathscr D_{\text{pos}}^{m,+}$, (i) is not feasible. Considering (ii), we have $H_{\text{pos}}(\omega_{j+1})=H_{\text{pos}}(\bold m)+\Gamma_{\text{pos}}(\bold m,\mathcal X^s_{\text{pos}})-h$. Starting from $\omega_{j+1}$ we consider only moves which imply either a decrease of energy or an increase by at most $h$. Since $C^1(\omega_{j+1})$ is a polyomino $\ell^*\times(\ell^*-1)$ with a bar made of two adjacent unit squares on a shortest side, the only feasible moves are
\begin{itemize}
\item[(iii)] flipping a spin, with two nearest neighbors with spin $m$, from $m$ to $1$,
\item[(iv)] flipping  a spin, with two nearest neighbors with spin $1$, from $1$ to $m$.
\end{itemize}
By means of the moves (iii) and (iv), the process reaches a configuration $\sigma$ in which all the spins are equal to $m$ except those, that are $1$, in a connected polyomino $C^1(\sigma)$ that is convex and such that $R(C^1(\sigma))=R_{(\ell^*+1)\times(\ell^*-1)}$.
We cannot repeat the move (iv) otherwise we get a configuration that does not belong to $\mathscr D^m_\text{pos}$. While applying one time (iv) and iteratively (iii), until we fill the rectangle $R_{(\ell^*+1)\times(\ell^*-1)}$ with spins $1$, we get a set of configurations in which the one with the smallest energy is $\sigma$ such that $C^1(\sigma)\equiv R(C^1(\sigma))$. Moreover, from any configuration in this set, a possible move
is reached by flipping from $m$ to $1$ a spin $m$ with three nearest neighbors with spin $m$ that implies to enlarge the circumscribed rectangle. This spin-flip increases the energy by $2-h$. Thus, we obtain
\begin{align}\label{absurdfinallygate3}
\Phi_\omega^\text{pos}&\ge 4\ell^*-h(\ell^*+1)(\ell^*-1)+2-h+H_{\text{pos}}(\bold m)\notag\\
&=4\ell^*-h(\ell^*)^2+2+H_{\text{pos}}(\bold m)\notag\\
&>\Gamma_{\text{pos}}(\bold m,\mathcal X^s_{\text{pos}})+H_{\text{pos}}(\bold m)=\Phi_\text{pos}(\bold m,\mathcal X^s_{\text{pos}}),
\end{align}
which is a contradiction by the definition of an optimal path. Note that the last inequality follows by $2>h(\ell^*-1)$ since $0<h<1$, see Assumption \ref{remarkconditionpos}. It follows that it is not possible to have $\omega\cap\mathcal W_{\text{pos}}(\bold m,\mathcal X^s_{\text{pos}})=\varnothing$ for any $\omega\in\Omega_{\bold m,\mathcal X^s_{\text{pos}}}^{opt}$, namely $\mathcal W_{\text{pos}}(\bold m,\mathcal X^s_{\text{pos}})$ is a gate for this type of transition. $\qed$

\subsection{Minimal gates for the transition from a metastable state to the other metastable states }\label{subsectionnew}
This subsection is devoted to the study of the transition from a metastable state to the set of the other metastable configurations. In Propositions \ref{propgateposmetameta} and \ref{unionWmingatemetameta} we identify geometrically two gates for this type of transition and in Theorem \ref{mingatetranmetameta} we show that the union of these sets gives the union of all the minimal gates for the same transition. Furthermore, in this subsection we also give some more details for the transition from any metastable state to the stable configuration $\mathbf 1$. More precisely, in Proposition \ref{thenocrossmeta} we prove that for any $\mathbf m\in\mathcal X^m_\text{pos}$ almost surely any optimal path $\omega\in\Omega_{\mathbf m,\mathcal X^s_\text{pos}}^{opt}$ does not visit any metastable state different from the initial one during the transition.
Let us begin by proving the following useful lemma.
\begin{lemma}\label{lemmastepsella}
For any $m\in\{2,\dots,q\}$, let $\eta\in\bar B_{\ell^*-1,\ell^*}^1(m,1)$ and let $\bar\eta\in\mathcal X$ a configuration which communicates with $\eta$ by one step of the dynamics. If the external magnetic field is positive, then either $H_{\text{\emph{pos}}}(\eta)<H_{\text{\emph{pos}}}(\bar\eta)$ or $H_{\text{\emph{pos}}}(\eta)>H_{\text{\emph{pos}}}(\bar\eta)$.
\end{lemma}
\textit{Proof.} Since $\eta$ and $\bar\eta$ differ by a single-spin update, let us define $\bar\eta:=\eta^{v,t}$ for some $v\in V$ and $t\in S$, $t\neq\eta(v)$. Note that $\eta\in\bar B_{\ell^*-1,\ell^*}^1(m,1)$ implies that $\eta$ is characterized by all spins $m$ except those that are $1$ in a quasi-square $(\ell^*-1)\times\ell^*$ with a unit protuberance on one of the longest sides. In particular, for any $v\in V$, either $\eta(v)=m$ or $\eta(v)=1$. If $\eta(v)=m$, then for any $t\in S\backslash\{m\}$, depending on the distance between the vertex $v$ and the $1$-cluster, we have 
\begin{align}
&H_{\text{pos}}(\bar\eta)-H_{\text{pos}}(\eta)=
\begin{cases}
4-h\mathbbm 1_{\{t=1\}},\ &\text{if $n_m(v)=4$ };\\
3-\mathbbm 1_{\{t=1\}}-h\mathbbm 1_{\{t=1\}},\ &\text{if $n_m(v)=3$, $n_1(v)=1$};\\
2-2\mathbbm 1_{\{t=1\}}-h\mathbbm 1_{\{t=1\}},\ &\text{if $n_m(v)=2$, $n_1(v)=2$}.
\end{cases}
\end{align}
Otherwise, if $\eta(v)=1$, for any $t\in S\backslash\{1\}$, depending on the distance between the vertex $v$ and the boundary of the $1$-cluster,  we get
\begin{align}
&H_{\text{pos}}(\bar\eta)-H_{\text{pos}}(\eta)=
\begin{cases}
4+h,\ &\text{if $n_1(v)=4$};\\
3-\mathbbm 1_{\{t=m\}}+h,\ &\text{if $n_m(v)=1$, $n_1(v)=3$};\\
2-2\mathbbm 1_{\{t=m\}}+h,\ &\text{if $n_m(v)=2$, $n_1(v)=2$};\\
1-3\mathbbm 1_{\{t=m\}}+h,\ &\text{if $n_m(v)=3$, $n_1(v)=1$}.
\end{cases}
\end{align}
We conclude that $H_{\text{{pos}}}(\eta)\neq H_{\text{{pos}}}(\bar\eta)$. 
$\qed$

In the next proposition we prove that the communication energy between metastable states is equal to the one between a metastable state and the stable state $\mathcal X^s_\text{pos}=\{\mathbf 1\}$.
\begin{proposition}\label{theoremcomparisonpos}
Consider the $q$-state Potts model on a $K\times L$ grid $\Lambda$, with periodic boundary conditions and with positive external magnetic field. For any $\bold m\in\mathcal X^m_{\emph{pos}}$, 
\begin{align}\label{phimetametanonres}
\Phi_{\emph{pos}}(\bold m,\mathcal X^m_\emph{pos}\backslash\{\bold m\})=4\ell^*-h(\ell^*(\ell^*-1)+1)+H_\emph{pos}(\bold m)=\Phi_{\emph{pos}}(\bold m,\mathcal X^s_\emph{pos}).
\end{align}
\end{proposition}
\textit{Proof}. Let us divide the proof in two steps. First we compute an upper bound of $\Phi_\text{pos}(\bold m,\mathcal X^m_\text{pos}\backslash\{\bold m\})$, second a lower bound.\\
\textit{Upper bound}. For any $\bold z\in\mathcal X^m_\text{pos}\backslash\{\bold m\}$, we use Definition \ref{refpathmiopos} to construct two reference paths $\tilde\omega^{(1)}:\bold m\to\bold 1$ and $\tilde\omega^{(2)}:\bold z\to\bold 1$. Thus, we define the reference path  $\omega^*:\bold m\to\bold z$ as the concatenation of the reference path $\tilde\omega^{(1)}$ and the time reversal $(\tilde\omega^{(2)})^T:\bold 1\to\bold z$. Thus, $\omega^*=(\tilde\omega^{(1)},(\tilde\omega^{(2)})^T)$. By this definition of $\omega^*$, we have $\max_{\xi\in\omega^*} H_\text{pos}(\xi)=\max\{\max_{\xi\in\tilde\omega^{(1)}}H_\text{pos}(\xi),\max_{\xi\in(\tilde\omega^{(2)})^T}H_\text{pos}(\xi)\}$. In the proof of Lemma \ref{lemmaduepositive}, using equations \eqref{alignrecallrefpathmax1}--\eqref{arg maxposvero} and \eqref{lebelHconfrontomax}, we get that
\begin{align}
\max_{\xi\in\omega^*} H_\text{pos}(\xi)=4\ell^*-h(\ell^*(\ell^*-1)+1)+H_\text{pos}(\bold m).
\end{align}
Thus, applying the definition of the communication energy, we conclude that
\begin{align}\label{upperphinoresprop}
\Phi_\text{pos}(\bold m,\bold z)=\min_{\omega:\bold m\to\bold z}\max_{\xi\in\omega} H_\text{pos}(\xi)\le 4\ell^*-h(\ell^*(\ell^*-1)+1)+H_\text{pos}(\bold m).
\end{align}
\textit{Lower bound}. During the transition from $\bold m$ to any $\bold z\in\mathcal X^m_\text{pos}\backslash\{\bold m\}$ the process has to intersect at least once the set $\mathcal V^m_k:=\{\sigma\in\mathcal X:\ N_m(\sigma)=k\}$ for any $k=KL,\dots,0$. In particular, given $k^*:=\ell^*(\ell^*-1)+1$, the process has to visit at least once the set $\mathcal V^m_{|\Lambda|-k^*}$. Since $\mathcal V^m_{|\Lambda|-k^*}\equiv\mathscr D_{\text{pos}}^m$, from Lemma \ref{bottomDP}, we have  
\begin{align}\label{lowerboundphinoresalign}
H_\text{pos}(\mathscr F(\mathcal V^m_{|\Lambda|-k^*}))=4\ell^*-h(\ell^*(\ell^*-1)+1)+H_{\text{pos}}(\bold m).
\end{align} 
It follows that the process visits at least once the set $\mathcal V^m_{|\Lambda|-k^*}$ in a configuration with energy larger than or equal to the r.h.s. of \eqref{lowerboundphinoresalign}, i.e., we obtain the following lower bound for the communication height between metastable states
\begin{align}\label{lowerphinoresprop}
\Phi_\text{pos}(\bold m,\bold z)\ge\ 4\ell^*-h(\ell^*(\ell^*-1)+1)+H_\text{pos}(\bold m).
\end{align}
Thanks to \eqref{upperphinoresprop} and \eqref{lowerphinoresprop}, we conclude that \eqref{phimetametanonres} is satisfied. $\qed$

Exploiting the equality $\Phi_\text{pos}(\mathbf m,\mathcal X^s_\text{pos})=\Phi_\text{pos}(\mathbf m,\mathcal X^m_\text{pos}\backslash\{\mathbf m\})$ for any $\mathbf m\in\mathcal X^m_\text{pos}$, we are now able to state the following corollary and proposition.
\begin{cor}\label{Hmaxlemmametameta}
Let $\bold m\in\mathcal X^m_{\emph{pos}}$ and let $\omega\in\Omega_{\bold m,\mathcal X^m_\emph{pos}\backslash\{\bold m\}}^{opt}$. If the external magnetic field is positive, then $\omega\cap\mathscr F(\mathscr D_{\emph{pos}}^m)\neq\varnothing$. Hence, $\mathscr F(\mathscr D_{\emph{pos}}^m)$ is a gate for the transition $\bold m\to\mathcal X^m_\emph{pos}\backslash\{\bold m\}$.
\end{cor}
\textit{Proof.} Thanks to \eqref{phimetametanonres} the proof is analogous to the one of Corollary \ref{Hmaxlemma3}. We refer to Appendix \ref{appendixproofcormetameta} for the explicit proof . $\qed$

\begin{proposition}\label{propgateposmetameta}
If the external magnetic field is positive, for any $\bold m\in\mathcal X^m_{\emph{pos}}$, 
 $\mathcal W_{\emph{pos}}(\bold m,\mathcal X^s_{\emph{pos}})$ is a gate for the transition $\bold m\to\mathcal X^m_\emph{pos}\backslash\{\bold m\}$.
\end{proposition} 
\textit{Proof.} Thanks to \eqref{phimetametanonres} the proof is analogous to the one of Proposition \ref{propgatepos}. See Appendix \ref{appendicepropgatemetameta} for the detailed proof. $\qed$

\noindent Given $\mathbf m\in\mathcal X^m_\text{pos}$, the reader may be surprised that $\mathcal W_\text{pos}(\bold m,\mathcal X^s_\text{pos})$ is a minimal gate for both the transitions $\mathbf m\to\mathcal X^s_\text{pos}$ and $\mathbf m\to\mathcal X^m_\text{pos}\backslash\{\mathbf m\}$. Intuitively, the set $\Omega_{\bold m,\mathcal X^m_\text{pos}\backslash\{\mathbf m\}}$ is partitioned in two non-empty subsets, i.e., the set containing those paths $\omega\in\Omega_{\mathbf m,\mathcal X^m_\text{pos}\backslash\{\mathbf m\}}$ such that $\omega\cap\mathcal C^{\mathbf 1}_{\mathcal X^m_\text{pos}}(\Gamma_\text{pos}^m(\mathbf 1,\mathcal X^m_\text{pos}))\neq\varnothing$ and the set cointaining those paths that do not enter this cycle. Corollary \ref{corposmetameta} points out that the $\Omega_{\mathbf m,\mathcal X^m_\text{pos}\backslash\{\mathbf m\}}^{opt}$ is a subset of the first set, i.e.,
\begin{align}\label{inclusionoptmetameta}
\Omega_{\mathbf m,\mathcal X^m_\text{pos}\backslash\{\mathbf m\}}^{opt}\subseteq\{\omega\in\Omega_{\mathbf m,\mathcal X^m_\text{pos}\backslash\{\mathbf m\}}: \omega\cap\mathcal C^\mathbf 1_{\mathcal X^m_\text{pos}}(\Gamma_\text{pos}^m(\mathbf 1,\mathcal X^m_\text{pos}))\neq\varnothing\}.
\end{align}
More precisely, in Proposition \ref{thenocrossmeta} we show that almost surely the process started in $\mathbf m\in\mathcal X^m_\text{pos}$ does not visit any other metastable states before hitting the stable configuration $\mathcal X^s_\text{pos}=\{\mathbf 1\}$. In order to prove this result, first we need to introduce the following habitat and to show that almost surely during the transition from a metastable to the stable state the process does not exit from it.
For any $\mathbf m\in\mathcal X^m_\text{pos}$, let
\begin{align}\label{habitatpos}
\mathcal A_\text{pos}:=\{\sigma\in\mathcal X:\ H_{\text{pos}}(\sigma)<H_{\text{pos}}(\bold m)+\Gamma_{\text{pos}}(\bold m,\mathcal X^s_{\text{pos}})+\hat\delta/2\},
\end{align}
where $\hat\delta$ is the minimum energy gap between an optimal and a non-optimal path from $\mathbf m$ to $\mathcal X^s_\text{pos}$. Note that $\mathcal A_\text{pos}$ is a cycle and that the choice to give some results on the dynamics from a metastable to the stable state inside $\mathcal A_\text{pos}$ is justified by the following result.
\begin{proposition}\label{remarkoptpath}
Let $\mathcal A_\emph{pos}$ be the habitat defined in \eqref{habitatpos}. 
Then, $\mathscr F(\mathcal A_\emph{pos})=\mathcal X^s_\emph{pos}$ and $V(\mathcal A_\emph{pos})=\Gamma^m_\text{pos}$.
Moreover, for any $\omega\in\Omega_{\bold m,\mathcal X^s_{\text{pos}}}^{opt}$, during the transition from any $\bold m\in\mathcal X^m_\emph{pos}$ to $\mathcal X^s_\emph{pos}$ the process does not exit almost surely from $\mathcal A_\emph{pos}$, i.e., 
\begin{align}\label{insidehabitat}
\lim_{\beta\to\infty} \mathbb P_\beta(\tau^\bold m_{\mathcal X^s_\emph{pos}}<\tau^\bold m_{\partial\mathcal A_\emph{pos}})=1.
\end{align}
\end{proposition} 
\textit{Proof}. By Proposition \ref{stablesetposprop} we have $\mathcal X^s_\text{pos}=\{\bold 1\}$ and by the definition \eqref{habitatpos} we have $\bold 1\in\mathcal A_\text{pos}$. Hence, $\mathscr F(\mathcal A_\text{pos})=\mathcal X^s_\text{pos}=\{\bold 1\}$. Furthermore, by \eqref{habitatpos} we also get that $\mathcal X^m_\text{pos}\subset\mathcal A_\text{pos}$. Hence, using Theorem \ref{teometastablepos} we have that $V(\mathcal A_\text{pos})=\Gamma^m_\text{pos}$. Finally, \eqref{insidehabitat} is verified thanks to Equation (2.20) of \cite[Theorem 2.17]{manzo2004essential} applied to the cycle $\mathcal A_\text{pos}$. $\qed$\\
We are now able to prove the following result.
\begin{proposition}[Study of the transition from any $\bold m\in\mathcal X^m_\text{pos}$ to $\mathcal X^s_\text{pos}$]\label{thenocrossmeta}
If the external magnetic field is positive, then for any $\bold m\in\mathcal X^m_\emph{pos}$ we have that every optimal path from $\bold m$ to $\mathcal X^s_\emph{pos}$ almost surely does not intersect other metastable states. More precisely,
\begin{align}\label{alignpropbis}
\lim_{\beta\to\infty} \mathbb P(\tau^\bold m_{\mathcal X^m_\emph{pos}\backslash\{\bold m\}}>\tau^\bold m_{\mathcal X^s_\emph{pos}})=1.
\end{align}
\end{proposition}
\textit{Proof.} Let $\omega=(\omega_0=\bold m,\dots,\omega_n)\in\Omega_{\bold m,\mathcal X^s_\text{pos}}^{opt}$ and, for some $j<n$, let $\omega_j\in\mathcal W_{\text{pos}}(\bold m,\mathcal X^s_{\text{pos}})$. By Corollary \ref{corpos} and by Corollary \ref{corposmetameta} we get that almost surely the process started in $\bold m\in\mathcal X^m_\text{pos}$ visits $\mathcal W_{\text{pos}}(\bold m,\mathcal X^s_{\text{pos}})$ before hitting $\mathcal X^s_\text{pos}\cup\mathcal X^m_\text{pos}\backslash\{\bold m\}$. Hence, almost surely we have 
\begin{align}\label{pathalmostsurelyint}
(\omega_0,\dots,\omega_j)\cap(\mathcal X^s_\text{pos}\cup\mathcal X^m_\text{pos}\backslash\{\bold m\})=\varnothing.
\end{align}
Thus, our claim is to show that starting from $\omega_j$, the process arrives in $\mathcal X^s_\text{pos}$ before visiting $\mathcal X^m_\text{pos}\backslash\{\bold m\}$. By Lemma \ref{lemmastepsella}, we have that $\omega_{j+1}$ does not have the same energy value of $\omega_j$. Thus, starting from $\omega_j$, the path passes to a configuration with energy strictly lower or strictly higher than $H_\text{pos}(\omega_j)$. More precisely, for some $v\in V$ and some $t\in S$, let $\omega_{j+1}:=\omega_j^{v,t}$. We have to consider the following possibilities:
\begin{itemize}
\item[(a)] $v$ is the vertex in the unit protuberance in the $1$-cluster in $\omega_j$ and $t=m$;
\item[(b)] $v$ is a vertex with spin $m$ with two nearest-neighbors with spin $m$ and two nearest-neighbors with spin $1$ in $\omega_j$ and $t=1$;
\item[(c)] $v$ has spin $1$ (respectively $m$), $t=m$ (respectively $t=1$) and $v$ is not a vertex that follows in case (a) (respectively case (b));
\item[(d)] $t\in S\backslash\{1,m\}$.
\end{itemize}
If $v$ and $t$ are as in case (a), then $\omega_{j+1}\in\mathcal C^\bold m_{\mathcal X^s_\text{pos}}(\Gamma^m_\text{pos})$ and, starting from this configuration, almost surely the process comes back to $\bold m$. Indeed, by Equation (2.20) of \cite[Theorem 2.17]{manzo2004essential} applied the the cycle $\mathcal C^\bold m_{\mathcal X^s_\text{pos}}(\Gamma^m_\text{pos})$ we have that there exists $k_1>0$ such that for every $\beta$ sufficiently large
\begin{align}
\mathbb P(\tau^{\omega_{j+1}}_\bold m>\tau^{\omega_{j+1}}_{\partial\mathcal C^\bold m_{\mathcal X^s_\text{pos}}(\Gamma^m_\text{pos})})\le e^{-k_1\beta}.
\end{align}
Thus, we repeat the same arguments to reduce the proof again to the cases (a)--(d) above.
If $v$ and $t$ are as in case (b), then $\omega_{j+1}\in\mathcal C^\bold 1_{\mathcal X^m_\text{pos}}(\Gamma_\text{pos}(\bold 1,\mathcal X^m_\text{pos}))$ and almost surely the process visits $\mathcal X^s_\text{pos}=\{\bold 1\}$ before exiting from this cycle. Indeed, by Equation (2.20) of \cite[Theorem 2.17]{manzo2004essential} applied to the cycle $\mathcal C^\bold 1_{\mathcal X^m_\text{pos}}(\Gamma_\text{pos}(\bold 1,\mathcal X^m_\text{pos}))$ we have that there exists $k_2>0$ such that for every $\beta$ sufficiently large
\begin{align}
\mathbb P(\tau^{\omega_{j+1}}_{\mathcal X^s_\text{pos}}>\tau^{\omega_{j+1}}_{\partial\mathcal C^\bold 1_{\mathcal X^m_\text{pos}}(\Gamma_\text{pos}(\bold 1,\mathcal X^m_\text{pos}))})\le e^{-k_2\beta}.
\end{align}
Since almost surely  \eqref{pathalmostsurelyint} holds, we conclude that \eqref{alignpropbis} is verified. \\
Finally, we consider $v$ and $t$ as in case (c) and (d) and our claim is to prove that almost surely $\omega_{j+1}$ as in these two cases does not belong to any optimal path from $\bold m$ to $\mathcal X^s_\text{pos}$. Indeed, $H_\text{pos}(\omega_{j+1})>H_\text{pos}(\omega_j)$ and since the minimum increase of energy is $h$, it follows that 
\begin{align}\label{incresaeminimumen}
H_\text{pos}(\omega_{j+1})\ge H_\text{pos}(\omega_j)+h=\Phi_\text{pos}(\bold m,\mathcal X^s_\text{pos})+h,
\end{align}
where the last equality follows by $\omega_j\in\mathcal W_{\text{pos}}(\bold m,\mathcal X^s_{\text{pos}})$. Hence, by \eqref{incresaeminimumen} and by the definition of the habitat $\mathcal A_\text{pos}$, we get that $\omega_{j+1}\notin\mathcal A_\text{pos}$. However, by Proposition \ref{remarkoptpath} we have that almost surely the process started in $\bold m$ does not exit from $\mathcal A_\text{pos}$ before hitting its bottom, and thus the cases (c) and (d) do not belong to any optimal path from $\bold m$ to $\mathcal X^s_\text{pos}$. $\qed$
Exploiting Proposition \ref{thenocrossmeta}, we are now able to identify another gate for the transition from a metastable state to the set of the other metastable states.
\begin{proposition}\label{unionWmingatemetameta}
$\bigcup_{\mathbf z\in\mathcal X^m_\emph{pos}\backslash\{\mathbf m\}} \mathcal W_\emph{pos}(\mathbf z,\mathcal X^s_\emph{pos})$ is a gate for the transition $\mathbf m\to\mathcal X^m_\emph{pos}\backslash\{\mathbf m\}$.
\end{proposition}
\textit{Proof.} Using Proposition \ref{thenocrossmeta} we get that the process started in $\mathbf m$ almost surely visits $\mathcal X^s_\text{pos}=\{\mathbf1\}$ earlier than $\mathcal X^m_\text{pos}\backslash\{\mathbf m\}$. It follows that almost surely $\mathbf 1\in\omega$  for any $\omega\in\Omega_{\mathbf m,\mathcal X^m_\text{pos}\backslash\{\mathbf m\}}^{opt}$, and 
since it is necessary to complete the transition to $\mathcal X^m_\text{pos}\backslash\{\bold m\}$, $\omega$ almost surely has a subpath which goes from $\mathbf 1$ to some $\mathcal X^m_\text{pos}\backslash\{\mathbf m\}$.  Thus, exploiting the reversibility of the dynamics, Proposition \ref{propgatepos} and since any optimal path $\omega\in\Omega_{\mathbf m,\mathcal X^m_\text{pos}\backslash\{\mathbf m\}}^{opt}$ hits $\mathcal X^m_\text{pos}\backslash\{\mathbf m\}$ in any metastable state different from $\mathbf  m$ with the same probability, i.e., $\frac{1}{q-2}$, we get that $\bigcup_{\mathbf z\in\mathcal X^m_\text{pos}\backslash\{\mathbf m\}} \mathcal W_\text{pos}(\mathbf z,\mathcal X^s_\text{pos})$ is a gate for the transition $\mathbf 1\to\mathcal X^m_\text{pos}\backslash\{\mathbf m\}$. $\qed$

\subsection{Minimal gates: proof of the main results}\label{proofgates}
We are now able to prove Theorems \ref{teogateposset} and \ref{mingatetranmetameta}.

\textit{Proof of Theorem \ref{teogateposset}}. For any $\bold m\in\mathcal X^m_{\text{pos}}$, by Proposition \ref{propgatepos} we get that $\mathcal W_{\text{pos}}(\bold m,\mathcal X^s_{\text{pos}})$ is a gate for the transition from $\bold m$ to $\mathcal X^s_{\text{pos}}=\{\bold 1\}$. In order to prove that 
$\mathcal W_{\text{pos}}(\bold m,\mathcal X^s_{\text{pos}})$ is a minimal gate, we exploit \cite[Theorem 5.1]{manzo2004essential} and we show that any $\eta\in\mathcal W_{\text{pos}}(\bold m,\mathcal X^s_{\text{pos}})$ is an essential saddle. To this end, in view of the definition of an essential saddle given in Subsection \ref{modinddef}, for any $\eta\in\mathcal W_{\text{pos}}(\bold m,\mathcal X^s_{\text{pos}})$ we define an optimal path from $\bold m$ to $\mathcal X^s_{\text{pos}}$ that passes through $\eta$ and such that it reaches its maximum energy only in this configuration. In particular, the optimal path is defined by modifying the reference path $\tilde\omega$ of Definition \ref{refpathmiopos} in a such a way that $\tilde\omega_{\ell^*(\ell^*-1)+1}=\eta$ in which $C^1(\eta)$ is a quasi-square $\ell^*\times(\ell^*-1)$ with a unit protuberance. This is possible by choosing the intial vertex $(i,j)$ such that during the construction the cluster $C^1(\tilde\omega_{\ell^*(\ell^*-1)})$ coincides with the quasi-square in $\eta$ and in the next step the unit protuberance is added in the site as in $\eta$. It follows that $\tilde\omega\cap\mathcal W_{\text{pos}}(\bold m,\mathcal X^s_{\text{pos}})=\{\eta\}$ and by the proof of Lemma \ref{lemmaduepositive} we get $\text{arg max}_{\tilde\omega} H_\text{pos}=\{\eta\}$. To conclude, we prove \eqref{setmingatespos}, i.e., that $\mathcal W_{\text{pos}}(\bold m,\mathcal X^s_{\text{pos}})$ is the unique minimal gate for the transition $\bold m\to\mathcal X^s_\text{pos}$. Note that the above reference paths $\tilde\omega$ reach the energy $\Phi_\text{pos}(\bold m,\mathcal X^s_\text{pos})$ only in $\mathcal W_{\text{pos}}(\bold m,\mathcal X^s_{\text{pos}})$. Thus, we get  that for any $\eta_1\in\mathcal W_{\text{pos}}(\bold m,\mathcal X^s_{\text{pos}})$, the set $\mathcal W_{\text{pos}}(\bold m,\mathcal X^s_{\text{pos}})\backslash\{\eta_1\}$ is not a gate for the transition 
$\bold m\to\mathcal X^s_{\text{pos}}$ since, in view of the above construction, we have that there exists an optimal path $\tilde\omega$ such that $\tilde\omega\cap\mathcal W_{\text{pos}}(\bold m,\mathcal X^s_{\text{pos}})\backslash\{\eta_1\}=\varnothing$. Note that the uniqueness of the minimal gate follows by the condition $\frac 2 h\notin\mathbb N$, see Assumption \ref{remarkconditionpos}(ii).
 $\qed$ 
 \begin{remark}\label{remuness}
A saddle $\eta\in\mathcal S(\sigma,\sigma')$ is unessential if for any $\omega\in\Omega_{\sigma,\sigma'}^{opt}$ such that $\omega\cap\eta\neq\varnothing$ the following conditions are both satisfied:
\begin{itemize}
\item[(i)] $\{\text{argmax}_\omega H\}\backslash \{\eta\}\neq\varnothing$,
\item[(ii)] there exists $\omega'\in\Omega_{\sigma,\sigma'}^{opt}$ such that $\{\text{argmax}_{\omega'} H\}\subseteq\{\text{argmax}_\omega H\}\backslash \{\eta\}$.
\end{itemize}
\end{remark}

\noindent\textit{Proof of Theorem \ref{mingatetranmetameta}} By Proposition \ref{propgateposmetameta} we have that the set given in (a) is a gate for the transition $\mathbf m\to\mathcal X^m_\text{pos}\backslash\{\mathbf m\}$. Hence, our aim is to prove that $\mathcal W_{\text{pos}}(\bold m,\mathcal X^s_\text{pos})$ is a minimal gate for the same transition. In order to show that this set satisfies the definition of minimal gate given in Subsection \ref{modinddef}, we show that for any $\eta\in\mathcal W_{\text{pos}}(\bold m,\mathcal X^s_\text{pos})$ there exists an optimal path $\omega'\in\Omega_{\mathbf m,\mathcal X^m_\text{pos}\backslash\{\mathbf m\}}^{opt}$ such that $\omega'\cap(\mathcal W_{\text{pos}}(\bold m,\mathcal X^s_\text{pos})\backslash\{\eta\})=\varnothing.$
We construct this optimal path $\omega'$ as the reference path $\omega^*$ defined in the proof of the upper bound of Proposition \ref{theoremcomparisonpos} in such a way that at the step $k^*-1$ the rectangular $\ell^*\times(\ell^*-1)$ $s$-cluster is as in $\eta$ without the protuberance. For $k^*\le k\le k^*+\ell^*-1$, we proceed as follows. At step $k^*$ the unit protuberance is added in the same position as in $\eta$, and in the following steps the same side is filled flipping consecutively to $s$ spins $1$ that have two nearest neighbors with spin $s$. Thus, $\omega'\cap\mathcal W_{\text{pos}}(\bold m,\mathcal X^s_\text{pos})=\{\eta\}$ and the condition of minimality is satisfied. 
\noindent By Proposition \ref{unionWmingatemetameta} the set depicted in (b) is a gate for the transition $\mathbf m\to\mathcal X^m_\text{pos}\backslash\{\mathbf m\}$. Thus, our aim is to prove that $\bigcup_{\mathbf z\in\mathcal X^m_\text{pos}\backslash\{\mathbf m\}} \mathcal W_\text{pos}(\mathbf z,\mathcal X^s_\text{pos})$ is a minimal gate for the same transition. Similarly to the previous case we show that for any $\eta\in\bigcup_{\mathbf z\in\mathcal X^m_\text{pos}\backslash\{\mathbf m\}} \mathcal W_\text{pos}(\mathbf z,\mathcal X^s_\text{pos})$ there exists an optimal path $\omega'\in\Omega_{\mathbf m,\mathcal X^m_\text{pos}\backslash\{\mathbf m\}}^{opt}$ such that $\omega'\cap(\bigcup_{\mathbf z\in\mathcal X^m_\text{pos}\backslash\{\mathbf m\}} \mathcal W_\text{pos}(\mathbf z,\mathcal X^s_\text{pos})\backslash\{\eta\})=\varnothing$. We define this optimal path $\omega'$ as the reference path $\omega^*$ constructed in the proof of the upper bound of Proposition \ref{theoremcomparisonpos} in such a way that at the step $k^*-1$ the rectangular $\ell^*\times(\ell^*-1)$ $s$-cluster is as in $\eta$ without the protuberance. For $k^*\le k\le k^*+\ell^*-1$, we proceed as follows. At step $k^*$ the unit protuberance is added in the same position as in $\eta$, and in the following steps the same side is filled flipping consecutively to $s$ spins $1$ that have two nearest neighbors with spin $s$. Thus, $\omega'\cap\bigcup_{\mathbf z\in\mathcal X^m_\text{pos}\backslash\{\mathbf m\}} \mathcal W_{\text{pos}}(\bold z,\mathcal X^s_\text{pos})=\{\eta\}$ and the condition of minimality is verified.
\noindent Thus $\bigcup_{\mathbf z\in\mathcal X^m_\text{pos}} \mathcal W_\text{pos}(\mathbf z,\mathcal X^s_\text{pos})\subseteq\mathcal G_\text{pos}(\mathbf m,\mathcal X^m_\text{pos}\backslash\{\mathbf m\})$, and we conclude exploiting \cite[Theorem 5.1]{manzo2004essential} and showing that any 
\begin{align}\label{unessentialmetametapos}
\eta\in\mathcal S_\text{pos}(\mathbf m,\mathcal X^m_\text{pos}\backslash\{\mathbf m\})\backslash\bigcup_{\mathbf z\in\mathcal X^m_\text{pos}} \mathcal W_\text{pos}(\mathbf z,\mathcal X^s_\text{pos})
\end{align}
 is an unessential saddle for the transition $\mathbf m\to\mathcal X^m_\text{pos}\backslash\{\mathbf m\}$. 
To this end we prove that any $\eta$ as in \eqref{unessentialmetametapos} satisfies conditions Remark \ref{remuness} (i) and (ii). Indeed, let $\omega\in\Omega_{\mathbf m,\mathcal X^m_\text{pos}\backslash\{\mathbf m\}}^{opt}$ such that $\omega\cap\{\eta\}\neq\varnothing$. Note that condition (i) in Remark \ref{remuness} is satisfied since $\omega$ intersects at least once both $\mathcal W_\text{pos}(\mathbf m,\mathcal X^s_\text{pos})$ and $\bigcup_{\mathbf z\in\mathcal X^m_\text{pos}\backslash\{\mathbf m\}} \mathcal W_\text{pos}(\mathbf z,\mathcal X^s_\text{pos})$. Next we define an optimal path $\omega'\in\Omega_{\mathbf m,\mathcal X^m_\text{pos}\backslash\{\mathbf m\}}^{opt}$ in order to prove that also condition (ii) in Remark \ref{remuness} is satisfied. From Propositions \ref{propgateposmetameta} and \ref{unionWmingatemetameta}, there exist $\eta_1^*\in\omega\cap\mathcal W_\text{pos}(\mathbf m,\mathcal X^s_\text{pos})$ and $\eta_2^*\in\omega\cap\bigcup_{\mathbf z\in\mathcal X^m_\text{pos}\backslash\{\mathbf m\}} \mathcal W_\text{pos}(\mathbf z,\mathcal X^s_\text{pos})$. Thus, we construct $\omega'$ as the reference path defined in the proof of the upper bound of Proposition \ref{theoremcomparisonpos} in such a way that $\omega'\cap\mathcal W_\text{pos}(\mathbf m,\mathcal X^s_\text{pos})=\{\eta_1^*\}$, $\omega'\cap\bigcup_{\mathbf z\in\mathcal X^m_\text{pos}\backslash\{\mathbf m\}} \mathcal W_\text{pos}(\mathbf z,\mathcal X^s_\text{pos})=\{\eta_2^*\}$ and $\{\text{argmax}_{\omega'} H\}=\{\eta_1^*,\eta_2^*\}$.
$\qed$

\subsection{Tube of typical paths: proof of the main results}\label{prooftube}
In order to give the proof of Theorem \ref{teotubesetpos}, first we prove the following lemmas.

\begin{lemma}\label{lemmapbsottocritical}
 For any $m\in S\backslash\{1\}$, consider the local minima $\eta\in\bar R_{\ell,\ell-1}(m,1)$ with $\ell\le\ell^*$ and $\zeta\in\bar R_{\ell,\ell}(m,1)$ with $\ell\le\ell^*-1$. Let $\mathcal C(\eta)$ and $\mathcal C(\zeta)$ be the non-trivial cycles whose bottom are $\eta$ and $\zeta$, respectively. Thus,  
\begin{align}
\mathcal B(\mathcal C(\eta))&=\bar B^1_{\ell-1,\ell-1}(m,1);\\
\label{lemmatubeone}
\mathcal B(\mathcal C(\zeta))&=\bar B^1_{\ell-1,\ell}(m,1).
\end{align}
\end{lemma}
\textit{Proof.} For any $m\in S\backslash\{1\}$, let $\eta_1\in\bar R_{\ell,\ell-1}(m,1)$ with $\ell\le\ell^*$. By Proposition \ref{proplocalminimapos}, $\eta_1\in\mathscr M^3_\text{pos}$ is a local minimum for the Hamiltonian $H_\text{pos}$. Using \eqref{principalboundary}, our aim is to prove the following
\begin{align}\label{claimonetubeuno}
\bar B^1_{\ell-1,\ell-1}(m,1)=\mathscr F(\partial\mathcal C(\eta_1)).
\end{align}
In $\eta_1$, for any $v\in V$ the corresponding $v$-tile (see before Lemma \ref{stabletilespositive} for the definition) is one among those depicted in Figure \ref{figmattonellepos}(a), (b), (d), (e) and (n) with $r=m$. Starting from $\eta_1$, the spin-flip to $m$ (resp.~$1$) the spin $1$ (resp.~$m$) on a vertex whose tile is one among those depicted in Figure \ref{figmattonellepos}(b), (e) (resp. (a), (d)), the process visits a configuration $\sigma_1$ such that
\begin{align}\label{maxincreases}
H_\text{pos}(\sigma_1)-H_\text{pos}(\eta_1)\ge2-h.
\end{align}
Thus, the smallest energy increase is given by $h$ by flipping to $m$ a spin $1$ on a vertex $v_1$ centered in a tile as in Figure \ref{figmattonellepos}(n) with $r=m$. Let $\eta_2:=\eta_1^{v_1,m}\in\bar B_{\ell-1,\ell-1}^{\ell-2}(m,1)$.
 In $\eta_2$, for any $v\in V$ the corresponding $v$-tile is one among those depicted in Figure \ref{figmattonellepos}(a), (b), (d), (e), (n) and (l) with $r=1$.
 Since $H_\text{pos}(\eta_2)=H_\text{pos}(\eta_1)+h$, the spin-flips on a vertex whose tile is one among those depicted in Figure \ref{figmattonellepos}(a), (b), (d), (e)
 lead to $H_\text{pos}(\sigma_2)-H_\text{pos}(\eta_1)\ge2$. Thus, as in the previous case, the smallest energy increase is given by flipping to $m$ a spin $1$ on a vertex $v_2$ centered in a tile as Figure \ref{figmattonellepos}(n). Note that starting from $\eta_2$ the only spin-flip which decreases the energy leads to the bottom of $\mathcal C(\eta_1)$, namely in $\eta_1$.
Iterating the strategy, the same arguments hold as long as the uphill path towards $\mathscr F(\partial\mathcal C(\eta_1))$ visits $\eta_{\ell-1}\in\bar B_{\ell-1,\ell-1}^1(m,1)$. Indeed, in this type of configuration for any $v\in V$ the corresponding $v$-tile is one among those depicted in Figure \ref{figmattonellepos}(a), (b), (d), (e), (n)  and unstable tile (s) with $t=r=s=m$, and it is possible to decrease the energy by passing to a configuration that does not belong to $\mathcal C(\eta_1)$. More precisely, there exists a vertex $w$ such that its tile is as the one in Figure \ref{figmattonellepos}(s) with $t=r=s=m$. By flipping to $m$ the the spin $1$ on $w$ the energy decreases by $2-h$, and the process enters a new cycle visiting its bottom, i.e., a local minimum belonging to $\bar R_{\ell-1,\ell-1}(m,1)\subset\mathscr M_\text{pos}^3$, see Proposition \ref{proplocalminimapos}. 

\noindent Let us now note that
\begin{align}\label{etaelleuno}
H_\text{pos}(\eta_{\ell-1})-H_\text{pos}(\eta_1)=h(\ell-2).
\end{align}
Since $\ell\le\ell^*$, comparing \eqref{maxincreases} with \eqref{etaelleuno}, we get that $\eta_{\ell-1}\in\mathscr F(\partial\mathcal C(\eta_1))$, and \eqref{claimonetubeuno} is verified.

\noindent Let us now consider for any $m\in S\backslash\{1\}$ the local minimum $\zeta_1\in\bar R_{\ell,\ell}(m,1)\subset\mathscr M^3_\text{pos}$ with $\ell\le\ell^*-1$.
Arguing similarly to the previous case, we verify \eqref{lemmatubeone} by proving that
\begin{align}\label{claimonetube}
\bar B^1_{\ell-1,\ell}(m,1)=\mathscr F(\partial\mathcal C(\zeta_1)).
\end{align}
$\qed$
\begin{lemma}\label{lemmapbsupercritical}
 For any $m\in S\backslash\{1\}$, consider the local minimum $\eta\in\bar R_{\ell_1,\ell_2}(m,1)$ with $\min\{\ell_1,\ell_2\}\ge\ell^*$. Let $\mathcal C(\eta)$ be the non-trivial cycle whose bottom is $\eta$. Thus,  
\begin{align}
\mathcal B(\mathcal C(\eta))&=\bar B^1_{\ell_1,\ell_2}(m,1)\cup\bar B^1_{\ell_2,\ell_1}(m,1).
\end{align}
\end{lemma}
\textit{Proof.} For any $m\in S\backslash\{1\}$, let $\eta_1\in\bar R_{\ell_1,\ell_2}(m,1)$ with $\ell^*\le\ell_1\le\ell_2$. By Proposition \ref{proplocalminimapos}, $\eta_1\in\mathscr M^3_\text{pos}$ is a local minimum for the Hamiltonian $H_\text{pos}$. Using \eqref{principalboundary}, our aim is to prove the following
\begin{align}\label{claimonetube}
\bar B^1_{\ell_1,\ell_2}(m,1)\cup\bar B^1_{\ell_2,\ell_1}(m,1)=\mathscr F(\partial\mathcal C(\eta_1)).
\end{align}
In $\eta_1$, for any $v\in V$ the corresponding $v$-tile is one among those depicted in Figure \ref{figmattonellepos}(a), (b), (d), (e) and (n) with $r=m$. Let $v_1\in V$ such that the $v_1$-tile is as the one depicted in Figure \ref{figmattonellepos}(d), and let $\eta_2:=\eta_1^{v_1,1}$. Note that if $v_1$ is adjacent to a side of length $\ell_2$, then $\eta_2\in\bar B^1_{\ell_1,\ell_2}(m,1)$, otherwise $\eta_2\in\bar B^1_{\ell_2,\ell_1}(m,1)$. Without loss of generality, let us assume that $\eta_2\in\bar B^1_{\ell_1,\ell_2}(m,1)$. By simple algebraic calculation we obtain that
\begin{align}
H_\text{pos}(\eta_2)-H_\text{pos}(\eta_1)=2-h.
\end{align}
In $\eta_2$ for any $v\in V$ the corresponding $v$-tile is one among those depicted in Figure \ref{figmattonellepos}(a), (b), (d), (e), (n) and (l) with $r=1$. By flipping to $1$ a spin $m$ on a vertex $w$ whose tile is as the one depicted in Figure \ref{figmattonellepos}(l) with $r=1$ the energy decreases by $h$ and the process enters a cycle different from the previous one that is either the cycle $\bar{\mathcal C}$ whose bottom is a local minimum belonging to $\bar R_{\ell_1+1,\ell_2}(m,1)$, or a trivial cycle for which iterating this procedure the process enters $\bar{\mathcal C}$. Thus, $\bar B^1_{\ell_1,\ell_2}(m,1)\subseteq\partial\mathcal C(\eta_1)$. Similarly we prove that $\bar B^1_{\ell_2,\ell_1}(m,1)\subseteq\partial\mathcal C(\eta_1)$. 

Let us now note that starting from $\eta_1$ the smallest energy increase is $h$, and it is given by flipping to $m$ a spin $1$ on a vertex whose tile is as the one depicted in Figure \ref{figmattonellepos}(n) with $r=m$. Let us consider the uphill path $\omega$ started in $\eta_1$ and constructed by flipping to $m$ all the spins $1$ along a side of the rectangular $\ell_1\times\ell_2$ $1$-cluster, say one of length $\ell_1$. Using the discussion given in the proof of Lemma \ref{lemmapbsottocritical} and the construction of $\omega$, we get that the process intersects $\partial\mathcal C(\eta)$ in a configuration $\sigma$ belonging to $\bar B^1_{\ell_2-1,\ell_1}(m,1)$. By simple algebraic computations, we obtain the following
\begin{align}\label{etaelleunobis}
H_\text{pos}(\sigma)-H_\text{pos}(\eta_1)=h(\ell_2-1).
\end{align}
Since $\ell_2\ge\ell^*$, it follows that $H_\text{pos}(\sigma)>H_\text{pos}(\eta_2)$. Since by flipping to $m$ (resp.~$1$) the vertex centered in a tile as depicted in Figure \ref{figmattonellepos}(b), (e) (resp.~(a)), the energy increase is largest than or equal to $2+h$, it follows that \eqref{claimonetube} is satisfied. $\qed$\\

We are now able to prove Theorem \ref{teotubesetpos}

\textit{Proof of Theorem \ref{teotubesetpos}}. Following the same approach as \cite[Section 6.7]{olivieri2005large}, we geometrically characterize the tube of typical trajectories for the transition using the so-called ``standard cascades''. See \cite[Figure 6.3]{olivieri2005large} for an example of these objects. We describe the standard cascades in terms of the paths that are started in $\mathbf m$ and are vtj-connected to $\mathcal X^s_\text{pos}$. See \eqref{defvtjpaths} for the formal definition and see \cite[Lemma 3.12]{nardi2016hitting} for an equivalent characterization of these paths. We remark that  any typical path from $\mathbf m$ to $\mathcal X^s_\text{pos}$ is also an optimal path for the same transition.
\noindent In order to describe these typical paths we proceed similarly to \cite[Section 7.4]{olivieri2005large}, where the authors apply the model-independent results given in Section 6.7 to identify the tube of typical paths in the context of the Ising model. Thus, we define a vtj-connected cycle-path that is the concatenation of both trivial and non-trivial cycles that satisfy \eqref{cyclepathvtj}.  In Theorem \ref{teogateposset} we give the geometric characterization of all the minimal gates for the transition $\mathbf m\to\mathcal X^s_\text{pos}$. Let $\eta_1$ be a configuration belonging to one of these minimal gates. We begin by studying the first descent from $\eta_1$ both to $\mathbf m$ and to $\mathcal X^s_\text{pos}$. Then, we complete the description of $\mathfrak T_{\mathcal X^s_\text{pos}}(\mathbf m)$ by joining the time reversal of the first descent from $\eta_1$ to $\mathbf m$ with the first descent from $\eta_1$ to $\mathcal X^s_\text{pos}$. 

\noindent Let us begin by studying the standard cascades from $\eta_1$ to $\mathbf m$. Since a spin-flip from $1$ to $t\notin\{1,m\}$ implies an increase of the energy value equal to the increase of the number of the disagreeing edges, we consider only the splin-flips from $1$ to $m$ on those vertices belonging to the $1$-cluster.
Thus, starting from $\eta_1$ and given $v_1$ a vertex such that $\eta_1(v_1)=1$, since $H_\text{pos}(\eta_1)=\Phi_\text{pos}(\mathbf m,\mathcal X^s_\text{pos})$, we get
\begin{align}
H_\text{pos}(\eta_1^{v_1,m})=\Phi_\text{pos}(\mathbf m,\mathcal X^s_\text{pos})+n_1(v_1)-n_m(v_1)+h.
\end{align}
It follows that the only possibility in which the assumed optimality of the path is not contradicted is the one where $n_1(v_1)=1$ and $n_m(v_1)=3$. Thus, along the first descent from $\eta_1$ to $\mathbf m$ the process visits $\eta_2$ in which all the vertices have spin $m$ except those, which are $1$, in a rectangular cluster $\ell^*\times(\ell^*-1)$, i.e., $\eta_2\in\bar R_{\ell^*-1,\ell^*}(m,1)$. By Proposition \ref{proplocalminimapos}, $\eta_2\in\mathscr M^3_\text{pos}$ is a local minimum, thus according to \eqref{cyclepathvtj} we have to describe its non-trivial cycle and its principal boundary. Starting from $\eta_2$, the next configuration along a typical path is defined by flipping to $m$ a spin $1$ on a vertex $v_2$ on one of the four corners of the rectangular $1$-cluster. Indeed, since $H_\text{pos}(\eta_2)=\Phi_\text{pos}(\mathbf m,\mathcal X^s_\text{pos})-2+h$, we have
\begin{align}
H_\text{pos}(\eta_2^{v_2,1})=H_\text{pos}(\eta_2)+n_1(v_2)-n_m(v_2)+h=\Phi_\text{pos}(\mathbf m,\mathcal X^s_\text{pos})-2+2h+n_1(v_2)-n_m(v_2),
\end{align}
and the only possibility in which the assumed optimality of the path is not contradicted is $n_1(v_2)=2$ and $n_m(v_2)=2$. Then, a typical path towards $\mathbf m$ proceeds by eroding the $\ell^*-2$ unit squares with spin $1$ belonging to a side of length $\ell^*-1$ that are corners of the $1$-cluster and that belong to the same side of $v_2$. Each of the first $\ell^*-3$ spin-flips increases the energy by $h$, i.e., the smallest energy increase for any single step of the dynamics, and these uphill steps are necessary in order to exit from the cycle whose bottom is the local minimum $\eta_2$. After these $\ell^*-3$ steps, the process hits the bottom of the boundary of this cycle in a configuration $\eta_{\ell^*}\in\bar B_{\ell^*-1,\ell^*-1}^1(m,1)$, see Lemma \ref{lemmapbsottocritical}. The last spin-update, that flips from $1$ to $m$ the spin $1$ on the unit protuberance of the $1$-cluster, decreases the energy by $2-h$. Thus, the typical path arrives in a local minimum $\eta_{\ell^*+1}\in\bar R_{\ell^*-1,\ell^*-1}(m,1)$, i.e., it enters a new cycle whose bottom is a configuration in which all the vertices have spin $1$, except those, which are $1$, in a square $(\ell^*-1)\times(\ell^*-1)$ $1$-cluster. Summarizing the construction above, we have the following sequence of vtj-connected cycles
\begin{align}
\{\eta_1\},\mathcal C^{\eta_2}_{\mathbf m}(h(\ell^*-2)),\{\eta_{\ell^*}\},\mathcal C^{\eta_{\ell^*+1}}_{\mathbf m}(h(\ell^*-2)).
\end{align}
Iterating this argument, we obtain that the first descent from $\eta_1\in\mathcal W_\text{pos}(\mathbf m,\mathcal X^s_\text{pos})$ to $\mathbf m$ is characterized by the concatenation of those vtj-connected cycle-subpaths between the cycles whose bottom is the local minima in which all the vertices have spin equal to $m$, except those, which are $1$, in either a quasi-square $(\ell-1)\times\ell$ or a square $(\ell-1)\times(\ell-1)$ for any $\ell=\ell^*,\dots,1$, and whose depth is given by $h(\ell -2)$. More precisely, from a quasi-square to a square, a typical path proceeds by flipping to $m$ those spins $1$ on one of the shortest sides of the $1$-cluster. On the other hand, from a square to a quasi-square, it proceeds by flipping to $m$ those spins $1$ belonging to one of the four sides of the square. 
Thus, a standard cascade from $\eta_1$ to $\mathbf m$ is characterized by the sequence of those configurations that belong to
\begin{align}\label{firstdescenteeta1}
\bigcup_{\ell=1}^{\ell^*}\biggr[\bigcup_{l=1}^{\ell-1}\bar B^l_{\ell-1,\ell}(m,1)\cup\bar R_{\ell-1,\ell}(m,1)\cup\bigcup_{l=1}^{\ell-2}\bar B^l_{\ell-1,\ell-1}(m,1)\cup\bar R_{\ell-1,\ell-1}(m,1)\biggl]. 
\end{align}

\noindent Let us now consider the first descent from $\eta_1\in\bar B^1_{\ell^*-1,\ell^*}(m,1)$ to $\mathcal X^s_\text{pos}=\{\mathbf 1\}$. In order to not contradict the definition of an optimal path, we have only to consider those steps which flip to $1$ a spin $m$. Indeed, adding a spin different from $m$ and $1$ leads to a configuration with energy value strictly larger than $\Phi_\text{pos}(\mathbf m,\mathcal X^s_\text{pos})$. Thus, let $w_1$ be a vertex such that $\eta_1(w_1)=m$. Flipping the spin $m$ on the vertex $w_1$, we get
\begin{align}
H_\text{pos}(\eta_1^{w_1,1})=\Phi_\text{pos}(\mathbf m,\mathcal X^s_\text{pos})+n_m(w_1)-n_1(w_1)-h,
\end{align}
and the only feasible choice is $n_m(w_1)=2$ and $n_1(w_1)=2$ in $\eta_1$. Thus, $\eta_1^{w_1,1}\in\bar B^2_{\ell^*-1,\ell^*}(m,1)$, namely the bar is now of length two. Arguing similarly, we get that along the descent to $\mathbf 1$ a typical path proceeds by flipping from $m$ to $1$ the spins $m$ with two nearest-neighbors with spin $1$ and two nearest-neighbors with spin $m$ belonging to the incomplete side of the $1$-cluster. More precisely, it proceeds downhill visiting $\bar\eta_i\in\bar B_{\ell^*-1,\ell^*}^{i}(m,1)$ for any $i=2,\dots,\ell^*-1$ and $\bar\eta_{\ell^*}\in\bar R_{\ell^*,\ell^*}(m,1)$, that is a local minimum by Proposition \ref{proplocalminimapos}. In order to exit from the cycle whose bottom is $\bar\eta_{\ell^*}$, the process crosses the bottom of its boundary by creating a unit protuberance of spin $1$ adjacent to one of the four edges of the $1$-square, i.e., visits $\{\bar\eta_{\ell^*+1}\}$ where $\bar\eta_{\ell^*+1}\in\bar B_{\ell^*,\ell^*}^1(m,1)$. Indeed, starting from $\bar\eta_{\ell^*}\in\bar R_{\ell^*,\ell^*}(m,1)$ the energy minimum increase is obtained by flipping a spin $m$ with three nearest-neighbors with spin $1$ and one nearest-neighbor with spin $m$. Starting from $\{\bar\eta_{\ell^*+1}\}$, a typical path towards $\mathbf 1$ proceeds by enlarging the protuberance to a bar of length two to $\ell^*-1$, thus it visits $\bar\eta_{\ell^*+i}\in\bar B_{\ell^*,\ell^*}^i(m,1)$ for any $i=2,\dots,\ell^*-1$. Each of these steps decreases the energy by $h$, and after them the descent arrives in the bottom of the cycle, i.e., in the local minimum $\bar\eta_{2\ell^*}\in\bar R_{\ell^*,\ell^*+1}(m,1)$. Then, the process exits from this cycle through the bottom of its boundary, i.e., by adding a unit protuberance of spin $1$ on \textit{any} one of the four edges of the rectangular $\ell^*\times(\ell^*+1)$ $1$-cluster in $\bar\eta_{2\ell^*}$, see Lemma \ref{lemmapbsupercritical}. Thus, it visits the trivial cycle $\{\bar\eta_{2\ell^*+1}\}$, where $\bar\eta_{2\ell^*+1}\in\bar B^1_{\ell^*,\ell^*+1}(m,1)\cup\bar B^1_{\ell^*+1,\ell^*}(m,1)$. Note that  the resulting standard cascade is different from the one towards $\mathbf m$. Thus, summarizing the construction above, we have defined the following sequence of vtj-connected cycles
\begin{align}
\{\eta_1\},\mathcal C^{\bar\eta_{\ell^*}}_{\mathbf 1}(h(\ell^*-1)),\{\bar\eta_{\ell^*+1}\},\mathcal C^{\bar\eta_{2\ell^*}}_{\mathbf 1}(h(\ell^*-1)),\{\bar\eta_{2\ell^*+1}\}. 
\end{align}
Note that if $\bar\eta_{2\ell^*}\in\bar B^1_{\ell^*,\ell^*+1}(m,1)$, then the process enters the cycle whose bottom is a configuration belonging to $\bar R_{\ell^*+1,\ell^*+1}(m,1)$. On the other hand, if $\bar\eta_{2\ell^*}\in\bar B^1_{\ell^*+1,\ell^*}(m,1)$, then the standard cascade enters the cycle whose bottom is a configuration belonging to $\bar R_{\ell^*,\ell^*+2}(m,1)$. In the first case the cycle has depth $h\ell^*$, in the second case the cycle has depth $h(\ell^*-1)$. Iterating this argument, 
we get that the first descent from $\eta_1$ to $\mathcal X^s_\text{pos}=\{\mathbf 1\}$ is characterized by vtj-connected cycle-subpaths from $\bar R_{\ell_1,\ell_2}(m,1)$ to $\bar R_{\ell_1,\ell_2+1}(m,1)$ 
defined as the sequence of those configurations belonging to $\bar B^l_{\ell_1,\ell_2}(m,1)$ for any $l=1,\dots,\ell_2-1$. Enlarging the $1$-cluster, at a certain point, the process arrives in a configuration in which this cluster is either a vertical or a horizontal strip, i.e., it intersects one of the two sets defined in \eqref{setverticalstrip}--\eqref{sethorizontalstrip}. If the descent arrives in $\mathscr S_\text{pos}^v(m,1)$, then it proceeds by enlarging the vertical strip column by column. Otherwise, if it arrives in $\mathscr S_\text{pos}^h(m,1)$, then it enlarges the horizontal strip row by row. In both cases, starting from a configuration with an $1$-strip, i.e., a local minimum in $\mathscr M^2_\text{pos}$ by Proposition \ref{proplocalminimapos}, the path exits from its cycle by adding a unit protuberance with a spin $1$ adjacent to one of the two vertical (resp. horizontal) edges and increasing the energy by $2-h$. Starting from the trivial cycle given by this configuration with an $1$-strip with a unit protuberance, the standard cascade enters a new cycle and it proceeds downhill by filling the column (resp.~row) with spins $1$. More precisely, the standard cascade visits $K-1$ (resp.~$L-1$) configurations such that each of them is defined by the previous one flipping from $m$ to $1$ a spin $m$ with two nearest-neighbors with spin $m$ and two nearest-neighbors with spin $1$. Each of these spin-updates decreases the energy by $h$. The process arrives in this way to the bottom of the cycle, i.e., in a configuration in which the thickness of the $1$-strip has been enlarged by a column (resp.~row). Starting from this state with the new $1$-strip, we repeat the same arguments above until the standard cascade arrives in the trivial cycle of a configuration $\sigma$ with an $1$-strip of thickness $L-2$ (resp.~$K-2$) and with a unit protuberance. Starting from $\{\sigma\}$, the process enters the cycle whose bottom is $\mathbf 1$ and it proceeds downhill either by flipping from $m$ to $1$ those spins $m$ with two nearest-neighbors with spin $m$ and two nearest-neighbors with spin $1$, or by flipping to $1$ all the spins $m$ with three nearest-neighbors with spin $1$ and one nearest-neighbor with spin $m$. The last step flips from $m$ to $1$ the last spin $m$ with four nearest-neighbors with spin $1$. Note that if the vtj-connected cycle path $(\mathcal C_1,\dots,\mathcal C_n)$ is such that $(\mathcal C_1,\dots,\mathcal C_n)\cap\mathscr S_\text{pos}^v(m,1)\neq\varnothing$ (resp. $(\mathcal C_1,\dots,\mathcal C_n)\cap\mathscr S_\text{pos}^h(m,1)\neq\varnothing$), then $(\mathcal C_1,\dots,\mathcal C_n)\cap\mathscr S_\text{pos}^h(m,1)=\varnothing$ (resp. $(\mathcal C_1,\dots,\mathcal C_n)\cap\mathscr S_\text{pos}^v(m,1)=\varnothing$). 
Thus, the first descent from $\eta_1$ to $\mathcal X^s_\text{pos}$ is characterized by the sequence of those configurations that belong to 
\begin{align}\label{firstdescenteeta2}
&\bigcup_{\ell_1=\ell^*}^{K-1}\bigcup_{\ell_2=\ell^*}^{K-1} \bar R_{\ell_1,\ell_2}(m,1)\cup\bigcup_{\ell_1=\ell^*}^{K-1}\bigcup_{\ell_2=\ell^*}^{K-1}\bigcup_{l=1}^{\ell_2-1}\bar B_{\ell_1,\ell_2}^l(m,1)\cup\bigcup_{\ell_1=\ell^*}^{L-1}\bigcup_{\ell_2=\ell^*}^{L-1} \bar R_{\ell_1,\ell_2}(m,1)\notag\\
&\cup\bigcup_{\ell_1=\ell^*}^{L-1}\bigcup_{\ell_2=\ell^*}^{L-1}\bigcup_{l=1}^{\ell_2-1}\bar B_{\ell_1,\ell_2}^l(m,1)\cup\mathscr S_\text{pos}^v(m,1)\cup\mathscr S_\text{pos}^h(m,1).
\end{align}
Finally, the standard cascade from $\mathbf m$ to $\mathcal X^s_\text{pos}$ is given by \eqref{firstdescenteeta1}--\eqref{firstdescenteeta2}. Finally, \eqref{ristimetubenzbneg} follows by \cite[Lemma 3.13]{nardi2016hitting}.
$\qed$ 

\appendix
\section{Appendix}
\subsection{Additional material for Subsection \ref{stablevpos}}\label{appendixproofpos}
\subsubsection{Definition \ref{refpathmiopos}}\label{appendixproofposdef}
For any $\bold m\in\mathcal X^m_{\text{pos}}$, we define a \textit{reference path} $\tilde\omega:\bold m\to\bold 1$, $\tilde\omega=(\omega_0^*,\dots,\omega_{KL}^*)$ as the concatenation of the two paths ${\tilde\omega}^{(1)}:=(\bold 1=\tilde\omega_0,\dots,\tilde\omega_{(K-1)^2})$ and ${\tilde\omega}^{(2)}:=(\tilde\omega_{(K-1)^2},\dots,$ $\bold m=\tilde\omega_{KL})$. The path ${\tilde\omega}^{(1)}$ is defined as follows. We set $\tilde\omega_0:=\bold m$. Then, we define $\tilde\omega_1:={\tilde\omega_0}^{(i,j),1}$, where $(i,j)$ denotes the vertex which belongs to the row $r_i$ and to the column $c_j$ of $\Lambda$, for some $i=0,\dots,K-1$ and $j=0,\dots,L-1$. Sequentially, we flip clockwise from $m$ to $1$ all the vertices that sourround the vertex $(i,j)$ in order to depict a $3\times 3$ square of spins $1$. We iterate this construction until we get $\tilde\omega_{(K-1)^2}\in\bar R_{K-1,K-1}(m,1)$. See \cite[Figure 6(a)]{bet2021metastabilityneg} for an illustration of this procedure.
 This time the white squares denote those vertices with spin $m$, the black ones denote the vertices with spin $1$.  Note that by considering the periodic boundary conditions the definition of $\tilde\omega$ is general for any $i$ and $j$.\\
The path ${\tilde\omega}^{(2)}$ is defined as follows. Without loss of generality, assume that $\tilde\omega_{(K-1)^2}\in\bar R_{K-1,K-1}(m,1)$ has the cluster of spin $1$ in the first $c_0,\dots,c_{K-2}$ columns, see \cite[Figure 6(b)]{bet2021metastabilityneg}. Starting from this last configuration $\tilde\omega_{(K-1)^2}$ of $\tilde\omega^{(1)}$, we define $\tilde\omega_{(K-1)^2+1},\dots,\tilde\omega_{(K-1)^2+K-1}$ as a sequence of configurations in which the cluster of spin $1$ grows gradually by flipping the spins $m$ on the vertices $(K-1,j)$, for $j=0,\dots,K-2$. Thus, $\tilde\omega_{(K-1)^2+K-1}\in\bar R_{K-1,K}(m,1)$, as depicted in \cite[Figure 6(c)]{bet2021metastabilityneg}. 
Finally, we define $\tilde\omega_{(K-1)^2+K},\dots,\tilde\omega_{KL}$ as a sequence of configurations in which the cluster of spin $s$ grows gradually column by column. More precisely, starting from $\tilde\omega_{(K-1)^2+K-1}\in\bar R_{K-1,K}(m,1)$, ${\tilde\omega}^{(2)}$ passes through configurations in which the spins $m$ on columns $c_K,\dots,c_{L-1}$ become $1$. The procedure ends with $\tilde\omega_{KL}=\bold m$. 
\subsubsection{Proof of Lemma \ref{lemmaposuno}}\label{appendixproofposlemma61}
Consider the reference path of Definition \ref{refpathmiopos} and note that for any $i=0,\dots,KL$, $N_1(\tilde\omega_i)=i$. The reference path may be constructed in such a way that $\tilde\omega_{\ell^*(\ell^*-1)}:=\sigma$. Let $\gamma:=(\tilde\omega_{\ell^*(\ell^*-1)}=\sigma,\tilde\omega_{\ell^*(\ell^*-1)-1},\dots,\tilde\omega_0=\bold m)$ be the time reversal of the subpath $(\tilde\omega_0,\dots,\tilde\omega_{\ell^*(\ell^*-1)})$ of $\tilde\omega$. We claim that $\max_{\xi\in\gamma} H_\text{pos}(\xi)<4\ell^*-h(\ell^*(\ell^*-1)+1)+H_\text{pos}(\bold m).$
Indeed, $\tilde\omega_{\ell^*(\ell^*-1)}=\sigma,\tilde\omega_{\ell^*(\ell^*-1)-1},\dots,\tilde\omega_1$ is a sequence of configurations in which all the spins are equal to $m$ except those, which are $1$, in either a quasi-square $\ell\times(\ell-1)$ or a square $(\ell-1)\times(\ell-1)$ possibly with one of the longest sides not completely filled. For any $\ell=\ell^*,\dots,2$, the path $\gamma$ moves from $\bar R_{\ell,\ell-1}(m,1)$ to $\bar R_{\ell-1,\ell-1}(m,1)$ by flipping to $m$ the $\ell-1$ spins $1$ on one of the shortest sides of the $1$-cluster. In particular, $\tilde\omega_{\ell(\ell-1)-1}$ is obtained by $\tilde\omega_{\ell(\ell-1)}\in\bar R_{\ell,\ell-1}(m,1)$ by flipping the spin on a corner of the quasi-square from $1$ to $m$ and this increases the energy by $h$. The next $\ell-3$ steps are defined by flipping the spins on the incomplete shortest side from $1$ to $m$,  thus each step increases the energy by $h$. Finally, $\tilde\omega_{(\ell-1)^2}\in\bar R_{\ell-1,\ell-1}(m,1)$ is defined by flipping the last spin $1$ to $m$ and this decreases the energy by $2-h$. For any $\ell=\ell^*,\dots,2$, $h(\ell-2)<2-h$. Indeed, from \eqref{ellestar} and from Assumption \ref{remarkconditionpos}, we have $2-h>h(\ell^*-2)\ge h(\ell-2)$. Hence, $\max_{\xi\in\gamma} H_\text{pos}(\xi)=H_\text{pos}(\sigma)=4\ell^*-h(\ell^*(\ell^*-1)+1)-(2-h)+H_\text{pos}(\bold m)$ and the claim is verified. 
$\qed$
\subsubsection{Proof of Lemma \ref{lemmaduepositive}}\label{appendixproofposlemma62}
Let $\bold m\in\mathcal X^m_\text{pos}$ and let $\sigma\in\bar B_{\ell^*-1,\ell^*}^2(m,1)$.  Consider the reference path of Definition \ref{refpathmiopos} and assume that it is constructed in such a way that $\tilde\omega_{\ell^*(\ell^*-1)+2}:=\sigma$. Let $\gamma:=(\tilde\omega_{\ell^*(\ell^*-1)+2}=\sigma,\tilde\omega_{\ell^*(\ell^*-1)+3},\dots,\tilde\omega_{KL-1},\bold 1)$. Our aim is to prove that $\max_{\xi\in\gamma} H_\text{pos}(\xi)<4\ell^*-h(\ell^*(\ell^*-1)+1)+H_\text{pos}(\bold m)$. In particular, we prove this claim by showing that $\max_{\xi\in\tilde\omega}H_\text{pos}(\xi)=4\ell^*-h(\ell^*(\ell^*-1)+1)+H_\text{pos}(\bold m)$ and that $\gamma$ does not visit the unique configuration in which this maximum is reached.\\
Consider $\ell\le K-2$. We recall that ${\tilde\omega}^{(1)}$ is defined as a sequence of configurations in which all the spins are equal to $m$ except those, which are $1$, in either a square $\ell\times\ell$ or a quasi-square $\ell\times(\ell-1)$ possibly with one of the longest sides not completely filled.
For some $\ell\le K-2$, let $\tilde\omega_{\ell(\ell-1)}\in\bar R_{\ell-1,\ell}(m,1)$ and $\tilde\omega_{\ell^2}\in\bar R_{\ell,\ell}(m,1)$, then 
\begin{align}\label{alignrecallrefpathmax1}
\max_{\sigma\in\{\tilde\omega_{\ell(\ell-1)},\tilde\omega_{\ell(\ell-1)+1},\dots,\tilde\omega_{\ell^2}\}} H_{\text{pos}}(\sigma)=H_{\text{pos}}(\tilde\omega_{\ell(\ell-1)+1})=4\ell-h\ell^2+h\ell-h+H_{\text{pos}}(\bold m).
\end{align}
Otherwise, if $\tilde\omega_{\ell^2}\in\bar R_{\ell,\ell}(m,1)$ and $\tilde\omega_{\ell(\ell+1)}\in\bar R_{\ell,\ell+1}(m,1)$, then 
\begin{align}
\max_{\sigma\in\{\tilde\omega_{\ell^2},\tilde\omega_{\ell^2+1},\dots,\tilde\omega_{\ell(\ell+1)}\}} H_{\text{pos}}(\sigma)=H_{\text{pos}}(\tilde\omega_{\ell^2+1})=4\ell-h\ell^2+2-h+H_{\text{pos}}(\bold m).
\end{align}
Let $k^*:=\ell^*(\ell^*-1)+1$. By recalling the condition $\frac{2}{h}\notin\mathbb N$ of Assumption \ref{remarkconditionpos} and by studying the maxima of $H_{\text{pos}}$ as a function of $\ell$, we have
\begin{align}\label{arg maxposvero}
\text{arg max}_{{\tilde\omega}^{(1)}} H_{\text{pos}}=\{\tilde\omega_{k^*}\}.
\end{align}
Note that if $\frac 2 h$ belonged to $\mathbb N$, then $\tilde\omega_{k^*}$ and $\tilde\omega_{(\ell^*)^2+1}$ would have the same energy value.\\
Let us now study the maximum energy value reached along ${\tilde\omega}^{(2)}$. This path is constructed as a sequence of configurations whose clusters of spins $1$ wrap around $\Lambda$. Moreover, the maximum of the energy is reached at the first configuration of ${\tilde\omega}^{(2)}$, see Figure \ref{Hmaxbis} for a qualitative representation of the energy of the configurations in ${\tilde\omega}^{(2)}$. Indeed, 
\begin{align}
&H_{\text{pos}}(\tilde\omega_{(K-1)^2+j})-H_{\text{pos}}(\tilde\omega_{(K-1)^2+j-1})=-2-h,\ j=2,\dots,K-1,\notag\\
&H_{\text{pos}}(\tilde\omega_{(K-1)^2+K})-H_{\text{pos}}(\tilde\omega_{(K-1)^2+K-1})=2-h,\notag\\
&H_{\text{pos}}(\tilde\omega_{(K-1)^2+j})-H_{\text{pos}}(\tilde\omega_{(K-1)^2+j-1})=-h,\ j=K+1,\dots,2K-1,\notag\\
&H_{\text{pos}}(\tilde\omega_{(K-1)^2+2K})-H_{\text{pos}}(\tilde\omega_{(K-1)^2+2K-1})=2-h.\notag
\end{align}
Note that 
\begin{align}
H_{\text{pos}}(\tilde\omega_{(K-1)^2+1})-H_{\text{pos}}(\tilde\omega_{(K-1)^2+K})&=4K-4-h(K-1)^2-(2K-h((K-1)^2+K))\notag\\
&=2K-4+h(K-1)>0,
\end{align}
where the last inequality follows by $K\ge 3\ell^*$. Moreover, 
\begin{align}
H_{\text{pos}}(\tilde\omega_{(K-1)^2+K})-&H_{\text{pos}}(\tilde\omega_{(K-1)^2+2K})\notag\\
&=2K+2-h((K-1)^2+2K-(2K+2-h((K-1)^2+K)=K>0.
\end{align}
\begin{figure}[h!]
\centering
\begin{tikzpicture}[scale=0.7, transform shape]
\fill[color=black] (1,10) circle (1.5pt) node[above] {\large$H_{\text{pos}}(\tilde\omega_{(K-1)^2+1})$}; 
\fill[color=black] (2,8) circle (1.5pt) node[left] {\large$H_{\text{pos}}(\tilde\omega_{(K-1)^2+K-1})$}; 
\fill[color=black] (2.2,8.3) circle (1.5pt) node[above,right] {\large$H_{\text{pos}}(\tilde\omega_{(K-1)^2+K})$}; 
\fill[color=black] (3.2,6.3) circle (1.5pt)  node[left] {\large$H_{\text{pos}}(\tilde\omega_{(K-1)^2+2K-1})$}; 
\fill[color=black] (3.4,6.6) circle (1.5pt) node[above,right] {\large$H_{\text{pos}}(\tilde\omega_{(K-1)^2+2K})$}; 
\draw[-] (1,10) -- (2,8) -- (2.2,8.3)--(3.2,6.3)--(3.4,6.6)--(3.6,6.3);
\draw[dotted,thick] (3.6,6.3)--(3.8,6) (4.05,5.6)--(4.25,5.3);
\draw [-] (4.25,5.3)--(4.75,4.5);
\fill[color=black] (4.75,4.5) circle (1.5pt) node[right] {\large$H_{\text{pos}}(\bold 1)$};
\end{tikzpicture}
\caption{\label{Hmaxbis} Qualitative illustration of the energy of the configurations belonging to ${\tilde\omega}^{(2)}$.}
\end{figure}
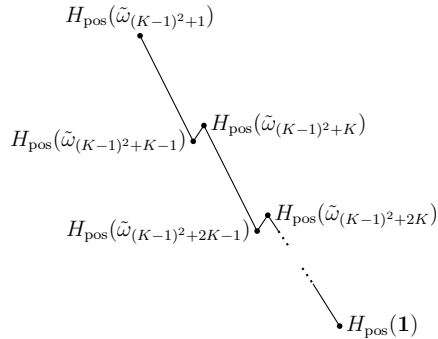\FloatBarrier

By iterating the analysis of the energy gap between two consecutive configurations along ${\tilde\omega}^{(2)}$, we conclude that 
\begin{align}\label{arg maxpos}
\text{arg max}_{{\tilde\omega}^{(2)}} H_{\text{pos}}=\{\tilde\omega_{(K-1)^2+1}\}.
\end{align}
In particular, 
\begin{align}\label{lebelHconfrontomax}
H_{\text{pos}}(\tilde\omega_{(K-1)^2+1})<H_{\text{pos}}(\tilde\omega_{k^*}).
\end{align}
This inequality is proved in \cite[Appendix A.1]{bet2021metastabilityneg}. Hence, $\text{arg max}_{\tilde\omega} H_{\text{pos}}=\{\tilde\omega_{k^*}\}$. Since $\gamma$ is constructed as a the subpath of $\tilde\omega$ which goes from $\tilde\omega_{\ell^*(\ell^*-1)+2}=\sigma$ to $\bold 1$, $\gamma$ does not visit the configuration $\tilde\omega_{k^*}$. Hence, the claim is verified.
$\qed$\\
\subsubsection{Proof of Proposition \ref{lowerboundpos}}\label{appendixproofposprop62}
For any $k=1,\dots,|V|$, let $\mathcal V^1_k:=\{\sigma\in\mathcal X: N_1(\sigma)=k\}$. Every path $\omega$ from any $\bold m\in\{\bold 2,\dots,\bold q\}$ to the stable configuration $\bold 1$ has to intersect the set $\mathcal V_k^1$ for every $k=1,\dots,|V|$. In particular, it has to visit the set $\mathcal V_{k^*}^1$ at least once, where $k^*=\ell^*(\ell^*-1)+1$. We prove the lower bound given in \eqref{lowerboundaligniphipos} by showing that $H_\text{pos}(\mathscr F(\mathcal V^1_{k^*}))=4\ell^*-h(\ell^*(\ell^*-1)+1)+H_{\text{pos}}(\bold m)$. Note that from \eqref{rewritegap1pos}, we get that the presence of disagreeing edges increases the energy. Thus, 
in order to describe the bottom $\mathscr F(\mathcal V^1_{k^*})$ we have to consider those configurations in which the $\ell^*(\ell^*-1)+1$ spins $1$ belong to a unique cluster inside a homogenous sea of spin $m\in S\backslash\{1\}$. Hence, consider $\tilde\omega$ be the reference path of Definition \ref{refpathmiopos} whose configurations satisfy this characterization. Note that $\tilde\omega\cap\mathcal V_{k^*}^1=\{\tilde\omega_{k^*}\}$ with $\tilde\omega_{k^*}\in\bar B^1_{\ell^*-1,\ell^*}(m,1)$. In particular,
\begin{align}\label{Homegakstar}
H_{\text{pos}}(\tilde\omega_{k^*})-H_{\text{pos}}(\bold m)=4\ell^*-h(\ell^*(\ell^*-1)+1),
\end{align}
where $4\ell^*$ is the perimeter of the cluster of spins $1$ in  $\tilde\omega_{k^*}$. We want to show that it is not possible to have a configuration with $k^*$ spins $1$ in a cluster of perimeter smaller than $4\ell^*$. Since the perimeter is an even integer, we suppose that there exists a configuration belonging to $\mathcal V^1_{k^*}$ such that the $1$-cluster has perimeter $4\ell^*-2$. Since $4\ell^*-2<4\sqrt{k^*}$, where $\sqrt{k^*}$ is the side-length of the square $\sqrt{k^*}\times\sqrt{k^*}$ of minimal perimeter among those in $\mathbb R^2$ of area $k^*$, using that the square is the figure that minimizes the perimeter for a given area, we conclude that there is no configuration with $k^*$ spins $1$ in a cluster with perimeter strictly smaller than $4\ell^*$. Hence, $\tilde\omega_{k^*}\in\mathscr F(\mathcal V^1_{k^*})$ and \eqref{lowerboundaligniphipos} is satisfied thanks to \eqref{Homegakstar}.  $\qed$

\subsubsection{Proof of Lemma \ref{lemmatrepospre}}\label{appendixproofposlemma63}
Let $\bold m\in\mathcal X^m_\text{pos}$. In the proof of Proposition \ref{lowerboundpos} we noted that any path $\omega:\bold m\to\bold 1$ has to visit $\mathcal V_k^1$  at least once for every $k=0,\dots,|V|$. Consider $\mathcal V_{\ell^*(\ell^*-1)}^1$. In \cite[Theorem 2.6]{alonso1996three} the authors show that the unique configuration of minimal energy in $\mathcal V_{\ell^*(\ell^*-1)}^1$ is the one in which all spins are $m$ except those that are $1$ in a quasi-square $\ell^*\times(\ell^*-1)$. In particular, this configuration has energy $\Phi_{\text{pos}}(\bold m,\mathcal X^s_{\text{pos}})-(2-h)=4\ell^*-2-h\ell^*(\ell^*-1)+H_\text{pos}(\bold m)$. Note that $4\ell^*-2$ is the perimeter of its $1$-cluster. Since the perimeter is an even integer, we have that the other configurations belonging to $\mathcal V_{\ell^*(\ell^*-1)}^1$ have energy that is larger than or equal to $4\ell^*-h\ell^*(\ell^*-1)+H_\text{pos}(\bold m)$. Thus, they are not visited by any optimal path. Indeed, $4\ell^*-h\ell^*(\ell^*-1)+H_\text{pos}(\bold m)>\Phi_{\text{pos}}(\bold m,\mathcal X^s_{\text{pos}})$.
Thus, we conclude that any optimal path intersects $\mathcal V_{\ell^*(\ell^*-1)}^1$ in a configuration belonging to $\bar R_{\ell^*-1,\ell^*}(m,1)$.
$\qed$

\subsection{Additional material for Subsection \ref{subsectionnew}}
\subsubsection{Proof of Corollary \ref{Hmaxlemmametameta}}\label{appendixproofcormetameta}
Every path from $\bold m\in\mathcal X^m_{\text{pos}}$ to the other metastable configurations in $\mathcal X^m_\text{pos}\backslash\{\bold m\}$ has to pass through the set $\mathcal V_k^m:=\{\sigma\in\mathcal X:\ N_m(\sigma)=k\}$ for any $k=|V|,\dots,0$. In particular, given $k^*:=\ell^*(\ell^*-1)+1$, any $\omega=(\omega_0,\dots,\omega_n) \in\Omega_{\bold m,\mathcal X^m_\text{pos}\backslash\{\bold m\}}^{opt}$ visits at least once the set $\mathcal V^m_{|\Lambda|-k^*}\equiv\mathscr D_\text{pos}^m$.  Hence, there exists $i\in\{0,\dots n\}$ such that $\omega_i\in\mathscr D_{\text{pos}}^m$. Thanks to \eqref{HbottomP} and to \eqref{phimetametanonres} we have that the energy value of any configuration belonging to $\mathscr F(\mathscr D_{\text{pos}}^m)$ is equal to the min-max reached by any optimal path from $\bold m$ to $\mathcal X^m_\text{pos}\backslash\{\bold m\}$. Thus, we conclude that $\omega_i\in\mathscr F(\mathscr D_{\text{pos}}^m)$.  $\qed$ 

\subsubsection{Proof of Proposition \ref{propgateposmetameta}}\label{appendicepropgatemetameta}
 For any $m\in S$, $m\neq 1$, let $\tilde{\mathscr D}^m_{\text{pos}}$ and $\hat{\mathscr D}^m_{\text{pos}}$ be the subsets of $\mathscr F(\mathscr D_{\text{pos}}^m)$ defined as follows. $\tilde{\mathscr D}^m_{\text{pos}}$ is  the set of those configurations of $\mathscr F(\mathscr D_{\text{pos}}^m)$ in which the boundary of the polyomino $C^1(\sigma)$ intersects each side of the boundary of its smallest surrounding rectangle $R(C^1(\sigma))$ on a set of the dual lattice $\mathbb{Z}^2+(1/2,1/2)$ made by at least two consecutive unit segments, see Figure \ref{figureexamplebis}(a). On the other hand, $\hat{\mathscr D}^m_{\text{pos}}$ is the set of those configurations of $\mathscr F(\mathscr D_{\text{pos}}^m)$ in which the boundary of the polyomino $C^1(\sigma)$ intersects at least one side of the boundary of $R(C^1(\sigma))$ in a single unit segment, see Figure \ref{figureexamplebis}(b) and (c). In particular note that $\mathscr F(\mathscr D_{\text{pos}}^m)=\tilde{\mathscr D}^m_{\text{pos}}\cup\hat{\mathscr D}^m_{\text{pos}}$. The proof proceeds in five steps.\\

\textbf{Step 1}. Our first aim is to prove that 
\begin{align}\label{Dhatmunion}
\hat{\mathscr D}^m_{\text{pos}}=\mathcal W_{\text{pos}}(\bold m,\mathcal X^s_{\text{pos}})\cup\mathcal W_{\text{pos}}'(\bold m,\mathcal X^s_{\text{pos}}).
\end{align}  
From \eqref{gatexm1} we have  $\mathcal W_{\text{pos}}(\bold m,\mathcal X^s_{\text{pos}})\cup\mathcal W_{\text{pos}}'(\bold m,\mathcal X^s_{\text{pos}})\subseteq\hat{\mathscr D}^m_{\text{pos}}$. Thus we reduce our proof to show that $\sigma\in\hat{\mathscr D}^m_{\text{pos}}$ implies $\sigma\in\mathcal W_{\text{pos}}(\bold m,\mathcal X^s_{\text{pos}})\cup\mathcal W_{\text{pos}}'(\bold m,\mathcal X^s_{\text{pos}})$. Note that this implication is not straightforward, since given $\sigma\in\hat{\mathscr D}^m_{\text{pos}}$, the boundary of the polyomino $C^1(\sigma)$ could intersect the other three sides of the boundary of its smallest surrounding rectangle $R(C^1(\sigma))$ in a proper subsets of the sides itself, see Figure \ref{figureexamplebis}(d) for an illustration of this hypothetical case. Hence, consider $\sigma\in\hat{\mathscr D}^m_{\text{pos}}$ and let $R(C^1(\sigma))=R_{(\ell^*+a)\times(\ell^*+b)}$ with $a,b \in\mathbb{Z}$. 
In view of the proof of Lemma \ref{bottomDP} we have that $C^1(\sigma)$ is a minimal polyomino and by \cite[Lemma 6.16]{cirillo2013relaxation} it is also convex and monotone, i.e., its perimeter of value $4\ell^*$ is equal to the one of $R(C^1(\sigma))$. Hence, the following equality holds
\begin{align}\label{equalityst}
4\ell^*=4\ell^*+2(a+b).
\end{align}
In particular, \eqref{equalityst} is satisfied only by $a=-b$. Now, let $\tilde R$ be the smallest rectangle surrounding the polyomino, say $\tilde C^1(\sigma)$, obtained by removing the unit protuberance from $C^1(\sigma)$. If $C^1(\sigma)$ has the unit protuberance adjacent to a side of length $\ell^*+a$, then $\tilde R$ is a rectangle $(\ell^*+a)\times(\ell^*-a-1)$. Note that $\tilde R$ must have an area larger than or equal to the number of spins $1$ of the polyomino $\tilde C^1(\sigma)$, that is $\ell^*(\ell^*-1)$. Thus, we have
\begin{align}\label{areartildepos1}
\text{Area}(\tilde R)=(\ell^*+a)(\ell^*-a-1)=\ell^*(\ell^*-1)-a^2-a\ge \ell^*(\ell^*-1) \iff -a^2-a\ge 0.
\end{align}
Since $a\in\mathbb{Z}$, $-a^2-a\ge 0$ is satisfied only if either $a=0$ or $a=-1$. 
Otherwise, if $C^1(\sigma)$ has the unit protuberance adjacent to a side of length $\ell^*-a$, then $\tilde R$ is a rectangle $(\ell^*+a-1)\times(\ell^*-a)$ and
\begin{align}\label{areartildepos2}
\text{Area}(\tilde R)=(\ell^*+a-1)(\ell^*-a)=\ell^*(\ell^*-1)-a^2+a\ge \ell^*(\ell^*-1) \iff -a^2+a\ge 0.
\end{align}
Since $a\in\mathbb{Z}$, $-a^2+a\ge 0$ is satisfied only if either $a=0$ or $a=1$. In both cases we get that $\tilde R$ is a rectangle of side lengths $\ell^*$ and $\ell^*-1$. Thus, if the protuberance is attached to one of the longest sides of $\tilde R$, then $\sigma\in\mathcal W_{\text{pos}}(\bold m,\mathcal X^s_{\text{pos}})$, otherwise $\sigma\in\mathcal W'_{\text{pos}}(\bold m,\mathcal X^s_{\text{pos}})$. In any case we conclude that \eqref{Dhatmunion} is satisfied.\\

\textbf{Step 2}. For any $\bold m\in\mathcal X^m_{\text{pos}}$ and for any path $\omega=(\omega_0,\dots,\omega_n)\in\Omega_{\bold m,\mathcal X^m_{\text{pos}}\backslash\{\bold m\}}^{opt}$, let 
\begin{align}\label{defgmomegametameta}
g_m(\omega):=\{i\in\mathbb{N}: \omega_i\in\mathscr F(\mathscr D_{\text{pos}}^m),\  N_1(\omega_{i-1})=\ell^*(\ell^*-1),\ N_m(\omega_{i-1})=|\Lambda|-\ell^*(\ell^*-1)\}.
\end{align}
We claim that $g_m(\omega)\neq\varnothing$. Let $\omega=(\omega_0,\dots,\omega_n)\in\Omega_{\bold m,\mathcal X^m_{\text{pos}}\backslash\{\bold m\}}^{opt}$ and let
$j^*\le n$ be the smallest integer such that after $j^*$ the path leaves $\mathscr D^{m,+}_\text{pos}$, i.e., $(\omega_{j^*},\dots,\omega_n)\cap\mathscr D^{m,+}_\text{pos}=\varnothing$.
 Since $\omega_{j^*-1}$ is the last configuration in $\mathscr D^{m,+}_\text{pos}$, it follows  that $\omega_{j^*}\in\mathscr D^m_\text{pos}$ and, by the proof of Corollary \ref{Hmaxlemmametameta}, we have that $\omega_{j^*}\in\mathscr F(\mathscr D^m_\text{pos})$. Moreover, since $\omega_{j^*-1}$ is the last configuration in $\mathscr D^{m,+}_\text{pos}$, we have that $N_m(\omega_{j^*-1})=|\Lambda|-\ell^*(\ell^*-1)$ and $\omega_{j^*}$ is obtained by $\omega_{j^*-1}$ by flipping a spin $m$ from $m$ to $s\neq m$. Note that $N_m(\omega_{j^*-1})=|\Lambda|-\ell^*(\ell^*-1)$ implies $N_s(\omega_{j^*-1})\le\ell^*(\ell^*-1)$ for any $s\in S\backslash\{m\}$. By Lemma \ref{bottomDP}, $\omega_{j^*}\in\mathscr F(\mathscr D^m_\text{pos})$ implies  $N_1(\omega_{j^*})=\ell^*(\ell^*-1)+1$, thus $N_1(\omega_{j^*-1})<\ell^*(\ell^*-1)$ is not feasible since $\omega_{j^*}$ and $\omega_{j^*-1}$ differ by a single spin update which increases the number of spins $1$ of at most one. Then, $j^*\in g_m(\omega)$ and the claim is proved.\\
 
\textbf{Step 3}. We claim that for any path $\omega\in\Omega_{\bold m,\mathcal X^m_{\text{pos}}\backslash\{\bold m\}}^{opt}$ one has $\omega_i\in\hat{\mathscr D}^m_{\text{pos}}$  for any $i\in g_m(\omega)$.
We argue by contradiction. Assume that there exists $i\in g_m(\omega)$ such that  $\omega_i\notin\hat{\mathscr D}^m_{\text{pos}}$ and $\omega_i\in\tilde{\mathscr D}^m_{\text{pos}}$.  Since $\omega_{i-1}$ is obtained from $\omega_i$ by flipping a spin $1$ to $m$ and since any configuration belonging to $\tilde{\mathscr D}^m_{\text{pos}}$ has all the spins $1$ with at least two nearest neighbors with spin $1$, using \eqref{energydifference3} we have
\begin{align}\label{disequalityHimetameta}
H_{\text{pos}}(\omega_{i-1})-H_{\text{pos}}(\omega_i)\ge(2-2)+h=h>0.
\end{align} 
In particular, from \eqref{disequalityHimetameta} we get a contradiction. Indeed, 
\begin{align}\label{absurdgate31metameta}
\Phi_\omega^\text{pos}\ge H_{\text{pos}}(\omega_{i-1})>H_{\text{pos}}(\omega_i)=H_{\text{pos}}(\bold m)+\Gamma_{\text{pos}}(\bold m,\mathcal X^m_{\text{pos}}\backslash\{\bold m\})=\Phi_{\text{pos}}(\bold m,\mathcal X^m_{\text{pos}}\backslash\{\bold m\}),
\end{align}
where the equality follows by \eqref{HbottomP}. Thus by \eqref{absurdgate31metameta} $\omega$ is not an optimal path, which is a contradiction, the claim is proved and we conclude the proof of Step 3.\\

\textbf{Step 4}. Now we claim that for any $\bold m\in\mathcal X^m_{\text{pos}}\backslash\{\bold m\}$ and for any path $\omega\in\Omega_{\bold m,\mathcal X^m_{\text{pos}}\backslash\{\bold m\}}^{opt}$,
\begin{align}
\omega_i\in\mathscr F(\mathscr D_{\text{pos}}^m)\implies \omega_{i-1}, \omega_{i+1}\notin\mathscr D_{\text{pos}}^m.
\end{align}
Using Corollary \ref{Hmaxlemmametameta}, for any $\bold m\in\mathcal X^m_\text{pos}$ and any path $\omega\in\Omega_{\bold m,\mathcal X^m_{\text{pos}}\backslash\{\bold m\}}^{opt}$ there exists an integer $i$ such that $\omega_i\in\mathscr F(\mathscr D_{\text{pos}}^m)$.  Assume by contradiction that $\omega_{i+1}\in\mathscr D_{\text{pos}}^m$. In particular,  since $\omega_i$ and $\omega_{i+1}$ have the same number of spins $m$, note that $\omega_{i+1}$ is obtained by flipping a spin $1$ from $1$ to $t\neq 1$. Since $\omega_i(v)\neq t$ for every $v\in V$,  the above flip increases the energy, i.e., $H_{\text{pos}}(\omega_{i+1})>H_{\text{pos}}(\omega_i)$. Hence, using this inequality and \eqref{HbottomP}, we have 
\begin{align}
\Phi_\omega^\text{pos}\ge H_{\text{pos}}(\omega_{i+1})>H_{\text{pos}}(\omega_i)=H_{\text{pos}}(\bold m)+\Gamma_{\text{pos}}(\bold m,\mathcal X^m_{\text{pos}}\backslash\{\bold m\})=\Phi_{\text{pos}}(\bold m,\mathcal X^m_{\text{pos}}\backslash\{\bold m\}),
\end{align}
which implies the contradiction because $\omega$ is not optimal. Thus $\omega_{i+1}\notin\mathscr D_{\text{pos}}^m$ and similarly we show that also $\omega_{i-1}\notin\mathscr D_{\text{pos}}^m$.\\

\textbf{Step 5}. In this last step of the proof we claim that for any $\bold m\in\mathcal X^m_{\text{pos}}$ and for any path $\omega\in\Omega_{\bold m,\mathcal X^m_{\text{pos}}\backslash\{\bold m\}}^{opt}$ there exists a positive integer $i$ such that $\omega_i\in\mathcal W_{\text{pos}}(\bold m,\mathcal X^s_{\text{pos}})$. Arguing by contradiction, assume that there exists $\omega\in\Omega_{\bold m,\mathcal X^m_{\text{pos}}\backslash\{\bold m\}}^{opt}$ such that $\omega\cap\mathcal W_{\text{pos}}(\bold m,\mathcal X^s_{\text{pos}})=\varnothing$. Thanks to Corollary \ref{Hmaxlemmametameta}, we know that $\omega$ visits $\mathscr F(\mathscr D^m_{\text{pos}})$ and thanks to Step 4 we have that the configurations along $\omega$ belonging to $\mathscr F(\mathscr D^m_{\text{pos}})$ are not consecutive. More precisely, they are linked by a sub-path that belongs either to $\mathscr D_{\text{pos}}^{m,+}$ or $\mathscr D_{\text{pos}}^{m,-}$. If $n$ is the length of $\omega$, then let $j\le n$ be the smallest integer such that $\omega_j\in\mathscr F(\mathscr D^m_{\text{pos}})$ and such that $(\omega_j,\dots,\omega_n)\cap\mathscr D_{\text{pos}}^{m,+}=\varnothing$, thus, $j\in g_m(\omega)$ since $j$ plays the same role of $j^*$ in the proof of Step 2. Using \eqref{Dhatmunion}, Step 3 and the assumption $\omega\cap\mathcal W_{\text{pos}}(\bold m,\mathcal X^s_{\text{pos}})=\varnothing$, it follows that $\omega_j\in\mathcal W_{\text{pos}}'(\bold m,\mathcal X^s_{\text{pos}})$. Moreover, starting from $\omega_j\in\mathscr F(\mathscr D^m_{\text{pos}})$ the energy along the path decreases only by either
\begin{itemize}
\item[(i)] flipping the spin in the unit protuberance from $1$ to $m$, or
\item[(ii)] flipping a spin, with two nearest neighbors with spin $1$, from $m$ to $1$.
\end{itemize}
Since by the definition of $j$ we have that $\omega_{j-1}$ is the last that visits $\mathscr D_{\text{pos}}^{m,+}$, $\omega_{j+1}\notin\mathscr D_{\text{pos}}^{m,+}$, (i) is not feasible. Considering (ii), we have $H_{\text{pos}}(\omega_{j+1})=H_{\text{pos}}(\bold m)+\Gamma_{\text{pos}}(\bold m,\mathcal X^m_{\text{pos}}\backslash\{\bold m\})-h$. Starting from $\omega_{j+1}$ we consider only moves which imply either a decrease of energy or an increase by at most $h$. Since $C^1(\omega_{j+1})$ is a polyomino $\ell^*\times(\ell^*-1)$ with a bar made of two adjacent unit squares on a shortest side, the only feasible moves are
\begin{itemize}
\item[(iii)] flipping a spin, with two nearest neighbors with spin $m$, from $m$ to $1$,
\item[(iv)] flipping  a spin, with two nearest neighbors with spin $1$, from $1$ to $m$.
\end{itemize}
By means of the moves (iii) and (iv), the process reaches a configuration $\sigma$ in which all the spins are equal to $m$ except those, that are $1$, in a connected polyomino $C^1(\sigma)$ that is convex and such that $R(C^1(\sigma))=R_{(\ell^*+1)\times(\ell^*-1)}$.
We cannot repeat the move (iv) otherwise we get a configuration that does not belong to $\mathscr D^m_\text{pos}$. While applying one time (iv) and iteratively (iii), until we fill the rectangle $R_{(\ell^*+1)\times(\ell^*-1)}$ with spins $1$, we get a set of configurations in which the one with the smallest energy is $\sigma$ such that $C^1(\sigma)\equiv R(C^1(\sigma))$. Moreover, from any configuration in this set, a possible move
is reached by flipping from $m$ to $1$ a spin $m$ with three nearest neighbors with spin $m$ that implies to enlarge the circumscribed rectangle. This spin-flip increases the energy by $2-h$. Thus, we obtain
\begin{align}\label{absurdfinallygate3}
\Phi_\omega^\text{pos}&\ge 4\ell^*-h(\ell^*+1)(\ell^*-1)+2-h+H_{\text{pos}}(\bold m)\notag\\
&=4\ell^*-h(\ell^*)^2+2+H_{\text{pos}}(\bold m)\notag\\
&>\Gamma_{\text{pos}}(\bold m,\mathcal X^m_{\text{pos}}\backslash\{\bold m\})+H_{\text{pos}}(\bold m)=\Phi_\text{pos}(\bold m,\mathcal X^m_{\text{pos}}\backslash\{\bold m\}),
\end{align}
which is a contradiction by the definition of an optimal path. Note that the last inequality follows by $2>h(\ell^*-1)$ since $0<h<1$, see Assumption \ref{remarkconditionpos}. It follows that it is not possible to have $\omega\cap\mathcal W_{\text{pos}}(\bold m,\mathcal X^s_{\text{pos}})=\varnothing$ for any $\omega\in\Omega_{\bold m,\mathcal X^m_{\text{pos}}\backslash\{\bold m\}}^{opt}$, namely $\mathcal W_{\text{pos}}(\bold m,\mathcal X^s_{\text{pos}})$ is a gate for this type of transition. $\qed$
\bibliographystyle{abbrv}
\bibliography{mybib}

\begin{thebibliography}{10}

\bibitem{alonso1996three}
L.~Alonso and R.~Cerf.
\newblock The three dimensional polyominoes of minimal area.
\newblock {\em The Electronic Journal of Combinatorics}, 3(1):R27, 1996.

\bibitem{ananikyan1995phase}
N.~Ananikyan and A.~Akheyan.
\newblock Phase transition mechanisms in the {P}otts model on a {B}ethe
  lattice.
\newblock {\em Journal of Experimental and Theoretical Physics},
  80(1):105--111, 1995.

\bibitem{apollonio2021metastability}
V.~Apollonio, V.~Jacquier, F.~R. Nardi, and A.~Troiani.
\newblock Metastability for the {I}sing model on the hexagonal lattice.
\newblock {\em arXiv:2101.11894}, 2021.

\bibitem{arous1996metastability}
G.~B. Arous and R.~Cerf.
\newblock Metastability of the three dimensional {I}sing model on a torus at
  very low temperatures.
\newblock {\em Electronic Journal of Probability}, 1, 1996.

\bibitem{bashiri2019on}
K.~Bashiri.
\newblock On the metastability in three modifications of the {I}sing model.
\newblock {\em View Journal Impact}, 25(3):483--532, 2019.

\bibitem{baxter1982critical}
R.~Baxter.
\newblock Critical antiferromagnetic square-lattice {P}otts model.
\newblock {\em Proceedings of the Royal Society of London. A. Mathematical and
  Physical Sciences}, 383(1784):43--54, 1982.

\bibitem{baxter1973potts}
R.~J. Baxter.
\newblock Potts model at the critical temperature.
\newblock {\em Journal of Physics C: Solid State Physics}, 6(23):L445, 1973.

\bibitem{baxter1978triangular}
R.~J. Baxter, H.~Temperley, and S.~E. Ashley.
\newblock Triangular {P}otts model at its transition temperature, and related
  models.
\newblock {\em Proceedings of the Royal Society of London. A. Mathematical and
  Physical Sciences}, 358(1695):535--559, 1978.

\bibitem{beltran2010tunneling}
J.~Beltran and C.~Landim.
\newblock Tunneling and metastability of continuous time {M}arkov chains.
\newblock {\em Journal of Statistical Physics}, 140(6):1065--1114, 2010.

\bibitem{beltran2012tunneling}
J.~Beltr{\'a}n and C.~Landim.
\newblock Tunneling and metastability of continuous time {M}arkov chains, the
  nonreversible case.
\newblock {\em Journal of Statistical Physics}, 149(4):598--618, 2012.

\bibitem{bet2021critical}
G.~Bet, A.~Gallo, and F.~R. Nardi.
\newblock Critical configurations and tube of typical trajectories for the
  {P}otts and {I}sing models with zero external field.
\newblock {\em arXiv:2102.06194}, 2021.

\bibitem{bet2021metastabilityneg}
G.~Bet, A.~Gallo, and F.~R. Nardi.
\newblock Metastability for the degenerate {P}otts {M}odel with negative
  external magnetic field under {G}lauber dynamics.
\newblock {\em arXiv:2105.14335}, 2021.

\bibitem{bet2020effect}
G.~Bet, V.~Jacquier, and F.~R. Nardi.
\newblock Effect of energy degeneracy on the transition time for a series of
  metastable states: application to probabilistic cellular automata.
\newblock {\em arXiv:2007.08342}, 2020.

\bibitem{bianchi2016metastable}
A.~Bianchi and A.~Gaudilliere.
\newblock Metastable states, quasi-stationary distributions and soft measures.
\newblock {\em Stochastic Processes and their Applications}, 126(6):1622--1680,
  2016.

\bibitem{bovier2016metastability}
A.~Bovier and F.~Den~Hollander.
\newblock {\em Metastability: a potential-theoretic approach}, volume 351.
\newblock Springer, 2016.

\bibitem{bovier2006sharp}
A.~Bovier, F.~den Hollander, and F.~R. Nardi.
\newblock Sharp asymptotics for {K}awasaki dynamics on a finite box with open
  boundary.
\newblock {\em Probability Theory and Related Fields}, 135(2):265--310, 2006.

\bibitem{bovier2002metastability}
A.~Bovier, M.~Eckhoff, V.~Gayrard, and M.~Klein.
\newblock Metastability and low lying spectral in reversible {M}arkov chains.
\newblock {\em Communications in Mathematical Physics}, 228(2):219--255, 2002.

\bibitem{bovier2004metastability}
A.~Bovier, M.~Eckhoff, V.~Gayrard, and M.~Klein.
\newblock Metastability in reversible diffusion processes {I}. {Sharp}
  asymptotics for capacities and exit times.
\newblock {\em Journal of the European Mathematical Society}, 2004.

\bibitem{boviermanzo2002metastability}
A.~Bovier and F.~Manzo.
\newblock Metastability in {{G}lauber} dynamics in the low-temperature limit:
  beyond exponential asymptotics.
\newblock {\em Journal of Statistical Physics}, 107(3-4):757--779, 2002.

\bibitem{cassandro1984metastable}
M.~Cassandro, A.~Galves, E.~Olivieri, and M.~E. Vares.
\newblock Metastable behavior of stochastic dynamics: a pathwise approach.
\newblock {\em Journal of Statistical Physics}, 35(5):603--634, 1984.

\bibitem{catoni1997exit}
O.~Catoni and R.~Cerf.
\newblock The exit path of a markov chain with rare transitions.
\newblock {\em ESAIM: Probability and Statistics}, 1:95--144, 1997.

\bibitem{cirillo1998metastability}
E.~N. Cirillo and J.~L. Lebowitz.
\newblock Metastability in the two-dimensional {I}sing model with free boundary
  conditions.
\newblock {\em Journal of Statistical Physics}, 90(1):211--226, 1998.

\bibitem{cirillo2003metastability}
E.~N. Cirillo and F.~R. Nardi.
\newblock Metastability for a stochastic dynamics with a parallel heat bath
  updating rule.
\newblock {\em Journal of Statistical Physics}, 110(1):183--217, 2003.

\bibitem{cirillo2013relaxation}
E.~N. Cirillo and F.~R. Nardi.
\newblock Relaxation height in energy landscapes: an application to multiple
  metastable states.
\newblock {\em Journal of Statistical Physics}, 150(6):1080--1114, 2013.

\bibitem{cirillo2015metastability}
E.~N. Cirillo, F.~R. Nardi, and J.~Sohier.
\newblock Metastability for general dynamics with rare transitions: escape time
  and critical configurations.
\newblock {\em Journal of Statistical Physics}, 161(2):365--403, 2015.

\bibitem{cirillo2008competitive}
E.~N. Cirillo, F.~R. Nardi, and C.~Spitoni.
\newblock Competitive nucleation in reversible probabilistic cellular automata.
\newblock {\em Physical Review E}, 78(4):040601, 2008.

\bibitem{cirillo2008metastability}
E.~N. Cirillo, F.~R. Nardi, and C.~Spitoni.
\newblock Metastability for reversible probabilistic cellular automata with
  self-interaction.
\newblock {\em Journal of Statistical Physics}, 132(3):431--471, 2008.

\bibitem{cirillo2017sum}
E.~N. Cirillo, F.~R. Nardi, and C.~Spitoni.
\newblock Sum of exit times in a series of two metastable states.
\newblock {\em The European Physical Journal Special Topics},
  226(10):2421--2438, 2017.

\bibitem{cirillo1996metastability}
E.~N. Cirillo and E.~Olivieri.
\newblock Metastability and nucleation for the {B}lume-{C}apel model.
  {D}ifferent mechanisms of transition.
\newblock {\em Journal of Statistical Physics}, 83(3):473--554, 1996.

\bibitem{costeniuc2005complete}
M.~Costeniuc, R.~S. Ellis, and H.~Touchette.
\newblock Complete analysis of phase transitions and ensemble equivalence for
  the {C}urie--{W}eiss--{P}otts model.
\newblock {\em Journal of Mathematical Physics}, 46(6):063301, 2005.

\bibitem{dai2015fast}
P.~Dai~Pra, B.~Scoppola, and E.~Scoppola.
\newblock Fast mixing for the low temperature 2d {I}sing model through
  irreversible parallel dynamics.
\newblock {\em Journal of Statistical Physics}, 159(1):1--20, 2015.

\bibitem{de1991metastability}
F.~de~Aguiar, L.~Bernardes, and S.~G. Rosa.
\newblock Metastability in the {P}otts model on the {C}ayley tree.
\newblock {\em Journal of Statistical Physics}, 64(3):673--682, 1991.

\bibitem{den2003droplet}
F.~Den~Hollander, F.~Nardi, E.~Olivieri, and E.~Scoppola.
\newblock Droplet growth for three-dimensional {K}awasaki dynamics.
\newblock {\em Probability Theory and Related Fields}, 125(2):153--194, 2003.

\bibitem{den2012metastability}
F.~den Hollander, F.~Nardi, and A.~Troiani.
\newblock Metastability for k{}awasaki dynamics at low temperature with two
  types of particles.
\newblock {\em Electronic Journal of Probability}, 17, 2012.

\bibitem{den2018metastability}
F.~den Hollander, F.~R. Nardi, and S.~Taati.
\newblock Metastability of hard-core dynamics on bipartite graphs.
\newblock {\em Electronic Journal of Probability}, 23, 2018.

\bibitem{di1987potts}
F.~di~Liberto, G.~Monroy, and F.~Peruggi.
\newblock The {P}otts model on {B}ethe lattices.
\newblock {\em Zeitschrift f{\"u}r Physik B Condensed Matter}, 66(3):379--385,
  1987.

\bibitem{ellis1990limit}
R.~S. Ellis and K.~Wang.
\newblock Limit theorems for the empirical vector of the
  {C}urie-{W}eiss-{P}otts model.
\newblock {\em Stochastic Processes and their Applications}, 35(1):59--79,
  1990.

\bibitem{ellis1992limit}
R.~S. Ellis and K.~Wang.
\newblock Limit theorems for maximum likelihood estimators in the
  {C}urie-{W}eiss-{P}otts model.
\newblock {\em Stochastic Processes and their Applications}, 40(2):251--288,
  1992.

\bibitem{enting1982triangular}
I.~Enting and F.~Wu.
\newblock Triangular lattice {P}otts models.
\newblock {\em Journal of Statistical Physics}, 28(2):351--373, 1982.

\bibitem{fernandez2015asymptotically}
R.~Fernandez, F.~Manzo, F.~Nardi, and E.~Scoppola.
\newblock Asymptotically exponential hitting times and metastability: a
  pathwise approach without reversibility.
\newblock {\em Electronic Journal of Probability}, 20, 2015.

\bibitem{fernandez2016conditioned}
R.~Fernandez, F.~Manzo, F.~Nardi, E.~Scoppola, and J.~Sohier.
\newblock Conditioned, quasi-stationary, restricted measures and escape from
  metastable states.
\newblock {\em Annals of Applied Probability}, 26(2):760--793, 2016.

\bibitem{gandolfo2010limit}
D.~Gandolfo, J.~Ruiz, and M.~Wouts.
\newblock Limit theorems and coexistence probabilities for the
  {C}urie-{W}eiss-{P}otts model with an external field.
\newblock {\em Stochastic Processes and their Applications}, 120(1):84--104,
  2010.

\bibitem{gaudillierelandim2014}
A.~Gaudilliere and C.~Landim.
\newblock {A Dirichlet principle for non reversible Markov chains and some
  recurrence theorems}.
\newblock {\em {Probability Theory and Related Fields}}, 158:55--89, 2014.

\bibitem{gaudilliere2005nucleation}
A.~Gaudilliere, E.~Olivieri, and E.~Scoppola.
\newblock Nucleation pattern at low temperature for local {K}awasaki dynamics
  in two dimensions.
\newblock {\em Markov Processes and Related Fields}, 11:553--628, 2005.

\bibitem{hollander2000metastability}
F.~d. Hollander, E.~Olivieri, and E.~Scoppola.
\newblock Metastability and nucleation for conservative dynamics.
\newblock {\em Journal of Mathematical Physics}, 41(3):1424--1498, 2000.

\bibitem{jovanovski2017metastability}
O.~Jovanovski.
\newblock Metastability for the {I}sing model on the hypercube.
\newblock {\em Journal of Statistical Physics}, 167(1):135--159, 2017.

\bibitem{kim2021metastability}
S.~Kim and I.~Seo.
\newblock Metastability of stochastic {I}sing and {P}otts models on lattices
  without external fields.
\newblock {\em arXiv:2102.05565}, 2021.

\bibitem{kotecky1994shapes}
R.~Koteck{\`y} and E.~Olivieri.
\newblock Shapes of growing droplets—a model of escape from a metastable
  phase.
\newblock {\em Journal of Statistical Physics}, 75(3):409--506, 1994.

\bibitem{manzo2004essential}
F.~Manzo, F.~R. Nardi, E.~Olivieri, and E.~Scoppola.
\newblock On the essential features of metastability: tunnelling time and
  critical configurations.
\newblock {\em Journal of Statistical Physics}, 115(1-2):591--642, 2004.

\bibitem{nardi2012sharp}
F.~Nardi and C.~Spitoni.
\newblock Sharp asymptotics for stochastic dynamics with parallel updating
  rule.
\newblock {\em Journal of Statistical Physics}, 146(4):701--718, 2012.

\bibitem{nardi1996low}
F.~R. Nardi and E.~Olivieri.
\newblock Low temperature stochastic dynamics for an {I}sing model with
  alternating field.
\newblock In {\em Markov Proc. Relat. Fields}, volume~2, pages 117--166, 1996.

\bibitem{nardi2019tunneling}
F.~R. Nardi and A.~Zocca.
\newblock Tunneling behavior of {I}sing and {P}otts models in the
  low-temperature regime.
\newblock {\em Stochastic Processes and their Applications},
  129(11):4556--4575, 2019.

\bibitem{nardi2016hitting}
F.~R. Nardi, A.~Zocca, and S.~C. Borst.
\newblock Hitting time asymptotics for hard-core interactions on grids.
\newblock {\em Journal of Statistical Physics}, 162(2):522--576, 2016.

\bibitem{neves1991critical}
E.~J. Neves and R.~H. Schonmann.
\newblock Critical droplets and metastability for a {G}lauber dynamics at very
  low temperatures.
\newblock {\em Communications in Mathematical Physics}, 137(2):209--230, 1991.

\bibitem{neves1992behavior}
E.~J. Neves and R.~H. Schonmann.
\newblock Behavior of droplets for a class of {G}lauber dynamics at very low
  temperature.
\newblock {\em Probability Theory and Related Fields}, 91(3-4):331--354, 1992.

\bibitem{olivieri1995markov}
E.~Olivieri and E.~Scoppola.
\newblock Markov chains with exponentially small transition probabilities:
  first exit problem from a general domain. the reversible case.
\newblock {\em Journal of Statistical Physics}, 79(3):613--647, 1995.

\bibitem{olivieri1996markov}
E.~Olivieri and E.~Scoppola.
\newblock Markov chains with exponentially small transition probabilities:
  first exit problem from a general domain. the general case.
\newblock {\em Journal of Statistical Physics}, 84(5):987--1041, 1996.

\bibitem{olivieri2005large}
E.~Olivieri and M.~E. Vares.
\newblock {\em Large deviations and metastability}, volume 100.
\newblock Cambridge University Press, 2005.

\bibitem{procacci2016probabilistic}
A.~Procacci, B.~Scoppola, and E.~Scoppola.
\newblock Probabilistic cellular automata for low-temperature 2-d {I}sing
  model.
\newblock {\em Journal of Statistical Physics}, 165(6):991--1005, 2016.

\bibitem{wang1994solutions}
K.~Wang.
\newblock Solutions of the variational problem in the {C}urie-{W}eiss-{P}otts
  model.
\newblock {\em Stochastic processes and their applications}, 50(2):245--252,
  1994.

\bibitem{zocca2019tunneling}
A.~Zocca.
\newblock Tunneling of the hard-core model on finite triangular lattices.
\newblock {\em Random Structures \& Algorithms}, 55(1):215--246, 2019.

\end{thebibliography}
\end{document}